\documentclass[reqno,11pt,a4paper]{amsart}

\usepackage[T1]{fontenc}
\usepackage{lmodern}
\usepackage{microtype}

\usepackage[margin=1in]{geometry}

\usepackage{amssymb}
\usepackage{mathtools}
\usepackage[dvipsnames]{xcolor}

\allowdisplaybreaks

\usepackage[colorlinks=true,linkcolor=RoyalBlue,citecolor=Red,urlcolor=RoyalBlue]{hyperref}
\usepackage{cleveref}
\crefname{assumption}{Assumption}{Assumptions}
\crefname{lemma}{Lemma}{Lemmas}
\crefname{proposition}{Proposition}{Propositions}
\crefname{theorem}{Theorem}{Theorems}
\crefname{corollary}{Corollary}{Corollaries}
\crefname{definition}{Definition}{Definitions}
\crefname{remark}{Remark}{Remarks}

\newcommand{\cH}{\mathcal{H}}

\newcommand{\R}{\mathbb{R}}
\newcommand{\eps}{\varepsilon}
\DeclareMathOperator{\diver}{div}
\DeclareMathOperator{\TV}{TV}
\DeclareMathOperator{\Lip}{Lip}
\DeclareMathOperator{\supp}{supp}

\title[Long-time asymptotics of nonlocal aggregation]{On the long-time behaviour and asymptotic geometry of nonlocal aggregation equations with saturating mobility}
\author[M.~Perugini]{Matteo Perugini}
\address[M.~Perugini]{Dipartimento di Ingegneria e Scienze dell'Informazione e Matematica (DISIM)\\
Universit\`a degli Studi dell'Aquila\\
via Vetoio 1, Coppito, 67100 L'Aquila, Italy}
\email[M.~Perugini]{matteo.perugini@univaq.it}
\author[E.~Radici]{Emanuela Radici}
\address[E.~Radici]{Dipartimento di Ingegneria e Scienze dell'Informazione e Matematica (DISIM)\\
Universit\`a degli Studi dell'Aquila\\
via Vetoio 1, Coppito, 67100 L'Aquila, Italy}
\email[E.~Radici]{emanuela.radici@univaq.it}
\keywords{Nonlocal aggregation equations; saturating mobility; entropy solutions; long-time asymptotics; strong $L^1$ convergence; star-shaped stationary states}
\subjclass[2020]{Primary 35B40; Secondary 35L65, 35Q49, 35R09}
\date{}

\newtheorem{theorem}{Theorem}[section]
\newtheorem{lemma}[theorem]{Lemma}
\newtheorem{proposition}[theorem]{Proposition}
\newtheorem{corollary}[theorem]{Corollary}
\newtheorem{assumption}[theorem]{Assumption}
\theoremstyle{definition}
\newtheorem{remark}[theorem]{Remark}
\newtheorem{definition}[theorem]{Definition}
\newtheorem*{pathstabilitylemma}{Lemma (uniform stability for prescribed centres)}
\newtheorem{example}[theorem]{Example}
\numberwithin{equation}{section}

\begin{document}

\begin{abstract}
We investigate a two-dimensional nonlocal conservation law with a saturating Lipschitz mobility, without assuming concavity of the mobility. For attractive radial potentials, we establish global entropy well-posedness and exact mass confinement. When the interaction potential is quadratic on the active interaction zone of the solution, the nonlocal convolution reduces to an affine drift. This reduction bounds the total path length of the moving barycentre via macroscopic energy dissipation, establishing strong $L^1$ relative compactness of the trajectory through a late-time core-filling approximation that traps the entropic state via $L^1$ contraction, circumventing the imposition of a uniform bounded variation estimate. This compactness, combined with the energy dissipation, forces the entropic trajectory to converge strongly in $L^1(\mathbb{R}^2)$ as $t \to \infty$ to a unique stationary binary state. We characterize the geometry of this limit, proving an equivalence with measure-theoretically star-shaped sets. We demonstrate that the nonlocal dynamics admit non-circular stationary equilibria and do not select a unique saturated disk absent radial symmetry on the initial datum. Finally, for entropic stationary binary states of finite perimeter, we exclude essentially bounded measure-theoretic components of the support complement under strict subharmonicity of the interaction potential, and rule out positively separated support components under uniform convexity.
\end{abstract}

\maketitle

\section{Introduction}

Nonlocal conservation laws with nonlinear mobility arise in the description of aggregation, congestion, and volume exclusion. We consider the two-dimensional equation
\begin{equation}\label{eq:nonlocal_main}
\partial_t\rho-\operatorname{div}\bigl(m(\rho)\nabla W*\rho\bigr)=0
\qquad\text{in }(0,\infty)\times\mathbb{R}^2.
\end{equation}
Here $\rho$ is a nonnegative density, $W$ is an interaction potential, and the nonlinear mobility $m$ vanishes both at vacuum and at a prescribed saturation level $R>0$. More precisely, $m:\mathbb{R}\to\mathbb{R}$ is globally Lipschitz, strictly positive on $(0,R)$, and zero on $\mathbb{R}\setminus(0,R)$; in particular, we impose no concavity assumption on $m$. In the general well-posedness theory we assume $W\in W^{2,\infty}(\mathbb{R}^2)$, whereas the confinement results concern radial attractive potentials satisfying $\nabla W(z)=w'(|z|)z/|z|$ with $w'(s)>0$ for $s>0$. The initial datum belongs to $L^1(\mathbb{R}^2)\cap L^\infty(\mathbb{R}^2)\cap BV(\mathbb{R}^2)$, takes values in $[0,R]$, and, for the asymptotic results, has compact support and unit mass. Strong long-time convergence is established under the additional assumption that $W(z)=|z|^2/2$ throughout the interaction region visited by the solution.

The vanishing of the mobility at the saturation level prevents unbounded aggregation but also creates a large family of stationary \emph{weak} solutions: every binary profile $R1_E$ has zero physical flux. Not every such profile, however, satisfies the entropy inequalities. This distinction, together with the possibility of shock formation and the lack of uniform-in-time spatial regularity estimates, makes the long-time selection of entropy equilibria a genuinely nontrivial problem.

\subsection{Related Literature}
\label{subsec:related_literature}

\paragraph{\textbf{Entropy solutions for nonlocal transport with nonlinear mobility.}} 
The analysis of scalar conservation laws relies on the entropy framework of Kru\v{z}kov \cite{Kruzhkov1970,Dafermos}. For nonlocal fluxes with congestion, Di Francesco, Fagioli, and Radici \cite{DiFrancescoFagioliRadici2019} constructed one-dimensional entropy solutions by deterministic particle approximations, proved uniqueness in their setting, and emphasized the distinction between stationary weak and entropy solutions. In particular, their discussion of separated saturated aggregates shows that vanishing physical flux alone does not ensure entropy admissibility. Subsequent one-dimensional particle approaches accommodate less regular interaction forces and more general congestion mechanisms \cite{RadiciStra2023}. Related particle constructions for aggregation--diffusion with nonlinear mobility were developed by Fagioli and Radici \cite{FagioliRadici2018}. A further recent study of one-dimensional nonlocal-flux conservation laws, with a distinct interaction kernel and a measure-solution formulation, is given by Di Francesco, Fagioli, and Radici \cite{DiFrancescoFagioliRadici2026}. Fagioli and Tse \cite{FagioliTse2022} established one-dimensional entropy well-posedness for a broad class of nonlocal transport equations with nonlinear mobility, together with a rigorous energy--dissipation balance. Stability estimates for quasi-entropy solutions of nonlocal scalar conservation laws were developed by Marconi, Radici, and Stra \cite{MarconiRadiciStra2025}. In several dimensions, Chodron de Courcel \cite{CdC} recently obtained local-in-time entropy well-posedness and quantitative vanishing-viscosity convergence under a different set of regularity assumptions on the nonlinear mobility and the interaction force. These works provide essential finite-time analytical tools, but the convergence of particle or viscous approximations to a PDE solution must be distinguished from convergence of that PDE solution as $t\to\infty$.\\

\paragraph{\textbf{Saturation, drift--diffusion, and long-time convergence.}} 
The closest recent comparison for full-trajectory convergence with a possibly nonconcave saturating mobility is the work of Carrillo, Fern\'andez-Jim\'enez, and G\'omez-Castro \cite{CarrilloFernandezJimenezGomezCastro2025}. For a class of drift--diffusion equations on bounded domains, driven by an internal energy and a \emph{prescribed} confinement potential, they construct free-energy-dissipating $L^1$-contractive semigroups and prove strong $L^1$ convergence of their trajectories to stationary states. They also describe free-boundary configurations and limits that need not minimize the free energy locally. Fern\'andez-Jim\'enez \cite{FernandezJimenez2026} has recently extended the corresponding analysis to $\mathbb{R}^d$ under suitable confinement assumptions. Thus, neither nonconcavity of the mobility nor strong convergence of the full trajectory, considered in isolation, distinguishes the present result. In both of these drift--diffusion works the confinement potential is prescribed and a diffusive internal-energy term is present. Carrillo, Fern\'andez-Jim\'enez, and G\'omez-Castro explicitly identify the introduction of an additional self-consistent aggregation field $W*\rho$ as a separate open direction \cite{CarrilloFernandezJimenezGomezCastro2025}; they do not cover general aggregation kernels. Our results address a particular, structurally tractable instance of a self-consistent interaction, with no diffusion.

Further results on saturating aggregation--confinement--diffusion equations, including nonlocal interactions, were established by Alamri \cite{Alamri2025}, who derives regularity and asymptotic convergence properties by exploiting a degenerate parabolic structure. More broadly, the analysis of aggregation--diffusion equilibria in \cite{CarrilloHittmeirVolzoneYao2019} yields radial-symmetry and long-time selection results under suitable assumptions; this contrasts with the nonradial entropy equilibria permitted by the purely hyperbolic saturation dynamics studied here. Numerical methods preserving density constraints and energy dissipation for continuity equations with saturation have also been developed in \cite{BailoCarrilloHu2023}.\\

\paragraph{\textbf{Nonlocal long-time dynamics and nonlinear-mobility transport geometry.}}
For attractive aggregation models without a saturation constraint, concentration into singular measures can occur in the absence of diffusion; see, for instance, Burger and Di Francesco \cite{BurgerDiFrancesco2008} for the asymptotic result and \cite{carrillo2011global} for finite-in-time measure aggregation. In the opposite, repulsive regime, Chodron de Courcel and Elbar \cite{ChodronElbar2025} prove exponential relaxation towards the spatial average for a nonlinear-mobility conservation law with Coulomb interaction on a torus. These results demonstrate that long-time convergence for purely nonlocal transport is not new in itself, while highlighting how the interaction sign, the mobility, and the admissible class of states change the nature of the limiting profiles. From the optimal-transport viewpoint, nonlinear mobilities also lead to modified dynamical transportation distances and generalized displacement-convexity questions \cite{DolbeaultNazaretSavare2009,CarrilloLisiniSavareSlepcev2010,LisiniMarigonda2010}. The approximation of such generalized transport metrics by deterministic particles is further studied by Di Marino, Portinale, and Radici \cite{DiMarinoPortinaleRadici2024}. Such metric constructions are typically based on concave mobilities; our argument neither assumes concavity nor uses geodesic convexity of the interaction energy.

\subsection{Our contribution}
\label{subsec:contribution}

Our first results establish global entropy well-posedness for globally Lipschitz saturating mobilities and show that, for radial attractive interactions, the mass remains in the closed convex hull of the initial support (Theorem~\ref{thm:main_well_posedness_confinement}). The main asymptotic contribution concerns the quadratic interaction regime. If the solution remains in $B_D(0)$ and $W(z)=|z|^2/2$ for $|z|\le 2D$, the nonlocal field becomes
\begin{equation}\label{eq:intro_affine_drift}
(\nabla W*\rho)(t,x)=x-c(t),\qquad
c(t)=\int_{\mathbb{R}^2}x\rho(t,x)\,dx,
\end{equation}
on the active support. Although the drift is affine in space, its centre is still determined self-consistently by the evolution. Dissipation of the centred quadratic interaction energy $E(t) := \frac{1}{2}\int_{\mathbb{R}^2}|x-c(t)|^2\rho(t,x)\,dx$, combined with a logarithmic moment identity, yields both the time integrability of the unweighted total mobility $M(t) = \int_{\mathbb{R}^2} m(\rho(t,x))\,dx$ and a finite total path length for $c(t)$.

The resulting strong $L^1$ compactness does \emph{not} follow from a uniform BV estimate for the physical trajectory. Indeed, the available a priori estimate has the form
\begin{equation}\label{eq:exponential_BV}
\TV(\rho(t))\le \TV(\rho_0)e^{\Lambda t},
\end{equation}
for some $\Lambda>0$ (see Proposition~\ref{prop:nonlocal-existence}), and therefore does not provide a uniform-in-time bound. Instead, we use rotational and dilational covariance of the fixed-centre entropy equation to obtain contraction of suitable directional variations. At a sufficiently late time, a small saturated core can be inserted into an $L^1$-close auxiliary BV datum. The remaining path length of the true barycentre is then small enough to preserve the auxiliary saturated core and control the auxiliary BV norm uniformly for all subsequent times. The $L^1$ contraction property for solutions driven by the \emph{same prescribed barycentre path} transfers this compactness to the original trajectory (Theorem~\ref{thm:main_compactness}). This contraction should not be confused with an $L^1$ contraction between arbitrary solutions of the self-consistent nonlocal equation.

Combining strong relative compactness with the time integrability of the total mobility, we prove that the entire trajectory is Cauchy in a negative Sobolev topology and converges strongly in $L^1(\mathbb{R}^2)$ to a binary stationary entropy state $\rho_\infty=R1_E$ (Theorem~\ref{thm:main_geometry}). The limit is uniquely selected by each fixed initial datum, but it is not a universal equilibrium. In fact, within the quadratic regime, stationary binary entropy solutions are characterized exactly by the essential star-shapedness of $E$ with respect to its barycentre, with no finite-perimeter assumption. Noncircular entropy equilibria therefore exist, and energy dissipation does not imply selection of the saturated disk in the absence of additional symmetry. Finally, for general interaction potentials we derive necessary measure-theoretic restrictions on stationary binary states of finite perimeter: strict subharmonicity excludes essentially bounded measure-theoretic components of the support complement, while uniform convexity precludes support components separated by a positive distance (Theorem~\ref{thm:main_topology}).

The preceding comparison delineates the scope of our contribution: strong full-trajectory convergence is proved for a two-dimensional, purely transport-type, self-consistent nonlocal equation with Lipschitz saturating mobility, under an essential quadratic-interaction hypothesis, without a uniform BV bound on the physical trajectory. The extension of this compactness mechanism to general attractive interactions remains open.

\subsection{Summary of Main Theorems}
\label{subsec:summary_main_theorems}

The first main result establishes the global well-posedness theory and verifies that attractive potentials confine the mass of $BV_{\textnormal{loc}}$ regular solutions.

\begin{theorem}[Global Well-Posedness and Exact Mass Confinement]\label{thm:main_well_posedness_confinement}
Let $m:\mathbb{R}\to\mathbb{R}$ be globally Lipschitz continuous, positive on $(0,R)$, and zero on $\mathbb{R}\setminus(0,R)$, and let $W\in W^{2,\infty}(\mathbb{R}^2)$. For every initial datum $\rho_{0}\in L^{1}(\mathbb{R}^{2})\cap L^{\infty}(\mathbb{R}^{2})\cap BV(\mathbb{R}^{2})$ with $0 \le \rho_0 \le R$, there exists a unique global Kru\v{z}kov entropy solution $\rho \in C([0,\infty);L^1(\mathbb{R}^2)) \cap L^\infty_{\mathrm{loc}}([0,\infty);L^\infty(\mathbb{R}^2)\cap BV(\mathbb{R}^2))$ of \eqref{eq:nonlocal_main}. Furthermore, if $W$ is radial and attractive, that is  $W(z) = w(|z|)$ with $w'(s) > 0$ for all $s > 0$, and $\rho_0$ is compactly supported with unit mass, the active mass is confined to the closed convex hull of its initial support.
\end{theorem}

The existence, uniqueness, and space--time $BV_{\mathrm{loc}}$ regularity are established in Section~\ref{sec:well_posedness}. Exact mass confinement via $BV$ trace tracking is proved in Section~\ref{subsec:mass_evolution}.

Although Theorem~\ref{thm:main_well_posedness_confinement} provides mass confinement and local-in-time $BV$ regularity for initial data in $L^{1}(\mathbb{R}^{2})\cap L^{\infty}(\mathbb{R}^{2})\cap BV(\mathbb{R}^{2})$, the compactness and asymptotic arguments in Sections~\ref{sec:autonomous_and_macroscopic}--\ref{sec:geometry} rely solely on $L^1$ confinement and do not use spatial $BV$ regularity of the physical trajectory $\rho(t)$. Under this quadratically confined regime, the macroscopic dissipation of the centered quadratic energy $E(t) := \frac{1}{2}\int_{\mathbb{R}^2}|x-c(t)|^2\rho(t,x)\,dx$ forces the moving barycentre $c(t) := \int_{\mathbb{R}^2}x\rho(t,x)\,dx$ to possess a finite total path length (Section~\ref{subsec:energy_dissipation}). This path-length integrability enables the central compactness result.

\begin{theorem}[Strong $L^1$ Relative Compactness]\label{thm:main_compactness}
Under the hypotheses of Theorem~\ref{thm:main_well_posedness_confinement}, assume further that $W$ is radially symmetric and  attractive, that is  $W(z) = w(|z|)$ with $w'(s) > 0$ for all $s > 0$. Moreover, assume that $W$ is quadratic, namely $W(z)= \frac{1}{2} |z|^2$, on the active interaction zone $B_{2D}(0)$. Let $\rho_0$ have unit mass with $\supp \rho_0 \subset B_D(0)$. Then the unique entropy solution $\rho(t)$ is confined to $B_{D}(0)$ for all $t \ge 0$, and the full trajectory $\{\rho(t):t\ge0\}$ is relatively compact in $L^{1}(\mathbb{R}^{2})$.
\end{theorem}

Compactness is established in Section~\ref{sec:compactness} by synthesizing the directional-variation contraction of the fixed-centre semigroup (Section~\ref{subsec:symmetries}), the persistence of saturated cores (Section~\ref{subsec:prescribed_path}), and a late-time core-filling approximation (Section~\ref{subsec:core_filling}) that traps the true trajectory via $L^1$ non-expansiveness without imposing a uniform $BV$ bound on $\rho(t)$.

Strong $L^1$ compactness, combined with the global time integrability of the unweighted total mobility $M(t) := \int_{\mathbb{R}^2} m(\rho(t,x))\,dx$, implies that the trajectory is Cauchy in $W^{-1,1}(\mathbb{R}^2)$ and hence converges strongly in $L^1(\mathbb{R}^2)$ as $t \to \infty$ to a limit $\rho_\infty$. The integrability of $M(t)$ forces $\int_{\mathbb{R}^2} m(\rho_\infty(x))\,dx = 0$, restricting the asymptotic limit to the binary set $\{0, R\}$ (Section~\ref{subsec:binarity}).

\begin{theorem}[Asymptotic Binarity and Star-Shaped Equilibria]\label{thm:main_geometry}
Under the hypotheses of Theorem~\ref{thm:main_compactness}, the trajectory $\rho(t)$ converges strongly in $L^{1}(\mathbb{R}^{2})$ (and consequently weakly-* in $L^\infty(\mathbb{R}^2)$) as $t\to\infty$ to a unique binary state $\rho_\infty = R1_E$ for some measurable set $E \subset B_D(0)$ of volume $|E| = 1/R$. Furthermore:
\begin{enumerate}
    \item \textbf{Stationary Entropy Condition:} $\rho_\infty$ is a stationary Kru\v{z}kov entropy solution of \eqref{eq:nonlocal_main}.
    \item \textbf{Star-Shaped Characterization and Non-Uniqueness:} For any measurable set $E \subset B_D(0)$ of volume $|E| = 1/R$, without requiring $E$ to have finite perimeter, the binary state $R1_E$ is a stationary Kru\v{z}kov entropy solution of \eqref{eq:nonlocal_main} if and only if $E$ is essentially star-shaped with respect to its barycentre (see \eqref{eq:measure_theoretic_star_shaped}). 
Consequently, the equation admits non-circular stationary equilibria and does not select a unique disk from general initial data, whereas radial symmetry of $\rho_0$ suffices to force strong convergence to the saturated disk. 
\end{enumerate}
\end{theorem}

Full-trajectory convergence and the binarity of the limit are proved in Section~\ref{subsec:binarity}. The star-shaped characterization of the asymptotic limit in the quadratic regime, which requires no bounded variation or finite-perimeter assumptions on $E$ and rules out interior radial holes up to $\mathcal{L}^2$-null sets, is established in Section~\ref{subsec:quadratic_geometry}. The proof that the dynamics do not force asymptotic convergence to a unique circular disk, illustrated by a geometric obstruction preventing disk formation from confined initial data (Example~\ref{ex:geometric_obstruction}), is detailed across Section~\ref{subsec:limitations_open_problems} and Section~\ref{subsec:quadratic_geometry}.

Finally, for general interaction potentials $W \in W^{2,\infty}(\mathbb{R}^2)$, we restrict to stationary binary states of finite perimeter ($\rho = R1_E \in BV(\mathbb{R}^2)$) and analyze the support structure via the outward normal trace of the nonlocal velocity field on the reduced boundary $\partial^* E$ (Section~\ref{subsec:general_potentials}).
\begin{theorem}[Measure-Theoretic Connectedness of Stationary Binary States]\label{thm:main_topology}
Let $m$ and $W \in W^{2,\infty}(\mathbb{R}^2)$ satisfy the hypotheses of Theorem~\ref{thm:main_well_posedness_confinement}. Let $\rho = R1_E \in L^{1}(\mathbb{R}^{2})\cap L^{\infty}(\mathbb{R}^{2})\cap BV(\mathbb{R}^{2})$ be a unit-mass stationary Kru\v{z}kov entropy solution of \eqref{eq:nonlocal_main} with $E \subset B_D(0)$. Then:
\begin{enumerate}
    \item \textbf{Indecomposability of the Complement:} If $\Delta W > 0$ almost everywhere on the active interaction zone $B_{2D}(0)$, the complement $\mathbb{R}^2 \setminus E$ is indecomposable in the sense of \cite{ACMM2001}, meaning that it contains no essentially bounded measure-theoretic connected components of positive Lebesgue measure.
    \item \textbf{Zero-Distance Component Separation:} If $W$ is uniformly convex on $B_{2D}(0)$, the support $E$ cannot be decomposed into two disjoint measurable sets of finite perimeter $E = E_1 \cup E_2$ with $|E_1| > 0$ and $|E_2| > 0$ such that $\overline{E}_i^* := \operatorname{supp}(1_{E_i})$ ($i=1,2$) lie at a positive distance; that is, $\operatorname{dist}(\overline{E}_1^*, \overline{E}_2^*) = 0$ for every such decomposition.
\end{enumerate}
\end{theorem}

\subsection{Strategy of the Proof: The Compactness Architecture}\label{subsec:quadratic_confinement_compactness_architecture}

The structural overview details the $L^1$ relative compactness mechanism, as extracting long-time asymptotics requires an architecture that bypasses the exponentially diverging spatial variation bound.

When the interaction potential $W$ is quadratic on the interaction zone of the solution $B_{2D}(0)$, the nonlocal convolution evaluates to an affine drift pointing away from the moving barycentre $c(t)=\int_{\mathbb{R}^2}x\rho(t,x)\,dx$. The interaction field becomes $\nabla W*\rho = x-c(t)$, reducing the nonlocal equation to a conservation law with an affine drift (Proposition~\ref{prop:exact_affine_reduction}).

The compactness architecture avoids imposing a uniform total variation bound directly onto the potentially irregular physical trajectory $\rho(t)$. Compactness is extracted by trapping $\rho(t)$ near a regularized auxiliary trajectory via $L^1$ contraction, relying on independent geometric and macroscopic estimates.

For a fixed point $a\in\mathbb{R}^2$, the auxiliary fixed-centre equation $\partial_t u - \operatorname{div}(m(u)(x-a)) = 0$ (Proposition~\ref{prop:fixed_centre_existence}) possesses symmetries under radial dilations and angular rotations. Through the $L^1$ contraction of the fixed-centre entropy semigroup (Proposition~\ref{prop:fixed_centre_comparison}), these symmetries translate into monotonically non-increasing directional variations (Theorem~\ref{thm:directional-variation_contraction}). The spatial degeneracy of these variations at the centre $a$ is neutralized by the persistence of a saturated core, a region $B_{r}(a)$ where the density equals the saturation threshold $R$ (Proposition~\ref{prop:discrete_path_length}). This core provides the coercivity required to reconstruct an unweighted $BV$ bound from the directional variations (Lemma~\ref{lem:core_coercivity}).

Macroscopically, the affine drift guarantees that the centered quadratic energy $E(t) = \frac{1}{2}\int_{\mathbb{R}^2}|x-c(t)|^2\rho(t,x)\,dx$ dissipates a.e.\ in time at the rate of the mobility-weighted variance, $D_m(t) := \int_{\mathbb{R}^2} m(\rho(t,x))|x-c(t)|^2\,dx$ (Lemma~\ref{lem:energy_dissipation}). This energy dissipation, evaluated against a logarithmic moment, ensures the unweighted total mobility $M(t) = \int_{\mathbb{R}^2} m(\rho(t,x))\,dx$ is globally integrable in time. This integrability bounds the total path length of the true barycentre, yielding $\int_0^\infty |\dot{c}(t)|\,dt < \infty$ (Proposition~\ref{prop:finite_path_length}).

Since the total path length is finite, there exists a time $T$ such that the remaining future path length of $c(t)$ is smaller than half of a designated core radius $r$. At $T$, the physical state $\rho(T)$ is approximated by an auxiliary state $v_T \in BV(\mathbb{R}^2)$ that is $L^1$-close to $\rho(T)$ and contains a saturated core of radius $r$ centered at $c(T)$ (Lemma~\ref{lem:bv_approximation_core}).

Evolving $v_T$ along the prescribed path of the true barycentre introduces a shifting reference origin that erodes the saturated core. The bound on the remaining path length prevents complete erosion. The core persists for all $t \ge T$, establishing a uniform future $BV$ bound for the auxiliary trajectory $v(t)$ (Theorem~\ref{thm:absolutely_continuous_prescribed_path}). By the Rellich--Kondrachov embedding, the auxiliary tail $\{v(t) : t \ge T\}$ is relatively compact in $L^1(\mathbb{R}^2)$. The prescribed-path entropy evolution is non-expansive in $L^1(\mathbb{R}^2)$ (Proposition~\ref{prop:prescribed_path_contraction}); contraction traps the true trajectory $\rho(t)$ near the compact auxiliary set, establishing the strong $L^1$ relative compactness of the original trajectory (Theorem~\ref{thm:relative_compactness}).

\subsection{Limitations and Open Problems}\label{subsec:limitations_open_problems}

The present compactness architecture is restricted to the quadratic potential regime. The framework relies on the affine reduction of the nonlocal field to $x-c(t)$, which unlocks the geometric symmetries of the fixed-centre equation and ensures the macroscopic energy dissipates at the rate of the mobility-weighted variance. Adaptation to general convex attractive potentials, where the interaction field is non-affine and the zero-force point diverges from the linear barycentre, is obstructed. Whether $L^1$ relative compactness holds for general attractive potentials without imposing uniform $BV$ bounds remains an open problem.

A second limitation concerns the geometric characterization and regularity of stationary states. A necessary and sufficient characterization of stationary binary entropy solutions—namely, their equivalence with measurable sets that are essentially star-shaped with respect to their barycentre (Proposition~\ref{prop:star_shaped})—is established only in the quadratic regime. For general potentials $W \in W^{2,\infty}(\mathbb{R}^2)$, Theorem~\ref{thm:main_topology} provides only necessary measure-theoretic conditions on the support $E$, which moreover require $E$ to be a set of finite perimeter. Because a set that is star-shaped with respect to a single point $c_\infty$ need not possess finite perimeter, it remains an open question whether an initial datum $\rho_0 \in L^{1}(\mathbb{R}^{2})\cap L^{\infty}(\mathbb{R}^{2})\cap BV(\mathbb{R}^{2})$ guarantees that the asymptotic limit $\rho_\infty = R1_E$ belongs to $BV(\mathbb{R}^2)$.

Finally, in the quadratic regime, Theorem~\ref{thm:main_geometry} guarantees that for any confined initial datum $\rho_0$, the trajectory $\rho(t)$ converges strongly in $L^1(\mathbb{R}^2)$ to a unique binary state $R1_E$ whose support $E$ is measure-theoretically star-shaped with respect to its barycentre. However, because the compactness argument is non-constructive, it yields no quantitative convergence rate, and outside the radially symmetric class it remains open how the specific geometry of the asymptotic star-shaped set $E$ depends on the initial configuration $\rho_0$. What can be established in this direction is that the mass confinement of Theorem~\ref{thm:main_well_posedness_confinement} rules out convergence to a disk for sufficiently thin initial supports, demonstrating that trajectories can converge to non-circular star-shaped states even though the exact final shape of $E$ cannot be determined a priori.

\begin{example}[Geometric Obstruction to Disk Formation]\label{ex:geometric_obstruction}
Let $m$ and $W$ satisfy the hypotheses of Theorem~\ref{thm:main_compactness} with $D > (\pi R)^{-1/2}$. Choose semi-axes
\[
\frac{1}{\pi R D} < b < (\pi R)^{-1/2}, \qquad \frac{1}{\pi R b} < a < D,
\]
and let $C := \{(x_1,x_2)\in\mathbb{R}^2 : x_1^2/a^2 + x_2^2/b^2 \le 1\} \subset B_D(0)$ be the corresponding closed ellipse, ensuring the total area $\pi a b$ of the ellipse exceeds $1/R$. Consider an initial datum $\rho_0 \in L^{1}(\mathbb{R}^{2})\cap L^{\infty}(\mathbb{R}^{2})\cap BV(\mathbb{R}^2)$ of unit mass with $0 \le \rho_0 \le R$ whose support is contained within $C$ (for instance, $\rho_0 = |C|^{-1}1_C$).

By mass confinement, the active mass is restricted to the closed convex hull of its initial support. Because $C$ is closed and convex, $\operatorname{supp} \rho(t) \subseteq C$ for all $t \ge 0$, which forces the support of the asymptotic limit to satisfy $\operatorname{supp} \rho_\infty \subseteq C$. A saturated circular disk of unit mass requires a radius of $r_* = (\pi R)^{-1/2}$. Since the minor semi-axis satisfies $b < r_*$, such a disk cannot fit inside the confinement hull $C$. Consequently, the asymptotic limit $\rho_\infty = R1_E$ is a star-shaped binary state supported in $C$ that is necessarily non-circular, showing that the dynamics do not select a disk from arbitrary initial configurations.
\end{example}

\subsection{Structure of the Paper}\label{subsec:structure}

The remainder of the paper is structured as follows. Section~\ref{sec:well_posedness} establishes the global existence, uniqueness, and space--time bounded variation regularity of entropy solutions, and proves mass confinement to the convex hull of the initial support for attractive potentials. Section~\ref{sec:autonomous_and_macroscopic} develops the fixed-centre framework to establish directional-variation contraction, and analyzes the macroscopic energy dissipation to bound the path length of the moving barycentre. Section~\ref{sec:compactness} proves the strong $L^1$ relative compactness of the trajectory via a core-filling approximation, establishing the vanishing of the intermediate phase and the strong $L^1$ convergence of the full trajectory to a unique stationary binary state. Section~\ref{sec:geometry} characterizes the geometry of binary stationary states, establishing outward normal trace equivalences and support structure for general potentials in Section~\ref{subsec:general_potentials}, followed by the measure-theoretic star-shaped characterization and radial disk selection for the quadratic regime in Section~\ref{subsec:quadratic_geometry}.

\section{Global Entropy Well-Posedness and Exact Mass Confinement}\label{sec:well_posedness}

In this section, we establish the global existence and uniqueness of entropy solutions for the nonlocal conservation law via a vanishing-viscosity construction. Subsequently, we demonstrate that attractive interaction potentials force the active mass to remain confined to the convex hull of its initial support.

\subsection{Smooth Approximants}\label{subsec:smooth_approx}

To invoke classical vanishing-viscosity frameworks while preserving the zero-density and saturation bounds, we first construct a sequence of smooth mobility approximants that vanish outside the admissible density regime.

\begin{assumption}[Mobility]\label{assump:mobility}
The mobility $m:\R\to\R$ is globally Lipschitz continuous, positive on $(0,R)$, and zero on $\R\setminus(0,R)$.
\end{assumption}

\begin{lemma}[Endpoint-preserving smooth approximation]\label{lem:endpoint_preserving_mobility}
Let $m$ satisfy Assumption~\ref{assump:mobility}, and put $L=\Lip(m)$. There exist nonnegative $m_j\in C_c^\infty(\R)$ with $\supp m_j\subset[0,R]$ such that
\begin{equation}\label{eq:asymptotically_sharp_Lipschitz_approx}
 \|m_j-m\|_{L^\infty(\R)}\longrightarrow0,\qquad
 \sup_j\Lip(m_j)\le2L,\qquad
 \limsup_{j\to\infty}\Lip(m_j)\le L.
\end{equation}
\end{lemma}

\begin{proof}
Extend $m$ by zero (as already required by Assumption~\ref{assump:mobility}). Choose a nonnegative standard mollifier $\eta\in C_c^\infty((-1,1))$ of unit mass, and let $g_\eps=m*\eta_\eps$ for $0<\eps<R/2$. Then $g_\eps\ge0$, $\supp g_\eps\subset[-\eps,R+\eps]$, $\Lip(g_\eps)\le L$, and $\|g_\eps-m\|_\infty\le L\eps$. Define the affine map
\[
 T_\eps(s)=\left(1+\frac{2\eps}{R}\right)s-\eps,
 \qquad m_\eps(s)=g_\eps(T_\eps(s)).
\]
Because $T_\eps(0)=-\eps$, $T_\eps(R)=R+\eps$, and $g_\eps$ vanishes smoothly outside $[-\eps,R+\eps]$, the function $m_\eps$ is smooth, nonnegative, and has support in $[0,R]$; no endpoint cutoff with a large derivative is needed. For $s\in[0,R]$, one has $|T_\eps(s)-s|\le\eps$, whence
\[
 |m_\eps(s)-m(s)|\le |g_\eps(T_\eps(s))-m(T_\eps(s))|+|m(T_\eps(s))-m(s)|\le2L\eps.
\]
Outside $[0,R]$ both functions vanish. Finally,
\[
 \Lip(m_\eps)\le\left(1+\frac{2\eps}{R}\right)L\le2L.
\]
Taking $\eps_j\downarrow0$ proves all claims, including the sharp limiting Lipschitz control.
\end{proof}
\subsection{Viscous Approximants and Entropy Solutions}\label{subsec:viscous_approx}

Throughout this subsection, we set $K:=-\nabla W$. Since $W\in W^{2,\infty}(\R^2)$, we have $K\in W^{1,\infty}(\R^2;\R^2)$ and equation \eqref{eq:nonlocal_main} is equivalently written as
\begin{equation}\label{eq:nonlocal-K}
 \partial_t u+\diver\bigl(m(u)(K*u)\bigr)=0.
\end{equation}

\begin{definition}[Kru\v{z}kov entropy solution]\label{def:nonlocal-entropy}
Let $T>0$ and $u_0\in L^1(\R^2)\cap L^\infty(\R^2)$. A function
\[
 u\in C([0,T];L^1(\R^2))\cap L^\infty((0,T)\times\R^2)
\]
is an entropy solution of \eqref{eq:nonlocal-K} on $[0,T]$ if it is a distributional solution with $u(0)=u_0$ in $L^1(\R^2)$ and, for every $k\in\R$ and every non-negative $\varphi\in C_c^1([0,T)\times\R^2)$,
\begin{align}
0\le {}&\int_0^T\!\int_{\R^2}\Bigl[
 |u-k|\,\partial_t\varphi
 +H(u,k)(K*u)\cdot\nabla\varphi
 -\operatorname{sign}(u-k)m(k)\bigl((\diver K)*u\bigr)\varphi
 \Bigr]\,dx\,dt \notag\\
&\qquad +\int_{\R^2}|u_0-k|\varphi(0,x)\,dx,\label{eq:nonlocal-entropy}
\end{align}
where $H(r,k):=\operatorname{sign}(r-k)\bigl(m(r)-m(k)\bigr)$ and $H(k,k):=0$.
\end{definition}

\begin{lemma}[Convolution bounds]\label{lem:nonlocal-conv-bounds}
Let $f\in L^1(\R^2)$. Then the convolutions $K*f$, $(DK)*f$, and $(\diver K)*f$ are well-defined almost everywhere and satisfy the uniform bounds
\begin{equation}\label{eq:conv-Linf-bounds}
 \|K*f\|_{L^\infty(\R^2)}\le \|K\|_{L^\infty(\R^2)}\|f\|_{L^1(\R^2)},
 \qquad
 \|(DK)*f\|_{L^\infty(\R^2)}\le \|DK\|_{L^\infty(\R^2)}\|f\|_{L^1(\R^2)},
\end{equation}
and
\begin{equation}\label{eq:conv-div-Linf-bound}
 \|(\diver K)*f\|_{L^\infty(\R^2)}\le \|\diver K\|_{L^\infty(\R^2)}\|f\|_{L^1(\R^2)}.
\end{equation}
Moreover, $K*f\in W^{1,\infty}(\R^2;\R^2)$, and its weak spatial Jacobian and divergence are given almost everywhere by
\begin{equation}\label{eq:conv-weak-derivatives}
 D(K*f)=(DK)*f, \qquad \diver(K*f)=(\diver K)*f.
\end{equation}
\end{lemma}

\begin{proof}
Since $W\in W^{2,\infty}(\R^2)$ and $K=-\nabla W$, we have $K\in L^\infty(\R^2;\R^2)$, $DK=-D^2 W\in L^\infty(\R^2;\R^{2\times 2})$, and $\diver K=\operatorname{tr}(DK)\in L^\infty(\R^2)$. For any measurable kernel $G\in L^\infty(\R^2)$ and any $f\in L^1(\R^2)$, Young's convolution inequality \cite[Theorem 4.2]{LiebLoss2001} yields
\begin{equation}\label{eq:young-Linf-L1-elementary}
 \|G*f\|_{L^\infty(\R^2)} \le \|G\|_{L^\infty(\R^2)}\|f\|_{L^1(\R^2)}.
\end{equation}
Choosing $G$ successively as $K$, $DK$, and $\diver K$ establishes \eqref{eq:conv-Linf-bounds} and \eqref{eq:conv-div-Linf-bound}. For any indices $i,j\in\{1,2\}$ and an arbitrary test function $\varphi\in C_c^\infty(\R^2)$, the absolute integrability of $K_i(x-y)f(y)\partial_{x_j}\varphi(x)$ permits the application of Fubini's theorem \cite[Theorem 1.12]{LiebLoss2001}:
\begin{align}\label{eq:conv-fubini-2}
 -\int_{\R^2}(K_i*f)(x)\,\partial_{x_j}\varphi(x)\,dx
 &=-\int_{\R^2}f(y)\left(\int_{\R^2}K_i(x-y)\,\partial_{x_j}\varphi(x)\,dx\right)dy \notag\\
 &=\int_{\R^2}f(y)\left(\int_{\R^2}(\partial_j K_i)(x-y)\,\varphi(x)\,dx\right)dy \notag\\
 &=\int_{\R^2}\bigl((\partial_j K_i)*f\bigr)(x)\,\varphi(x)\,dx.
\end{align}
Thus $\partial_{x_j}(K_i*f) = (\partial_j K_i)*f \in L^\infty(\R^2)$ in the sense of distributions, establishing \eqref{eq:conv-weak-derivatives}.
\end{proof}

\begin{proposition}[Global viscous approximants and uniform estimates]\label{prop:nonlocal-viscous}
Let $u_0\in L^1(\R^2)\cap L^\infty(\R^2)\cap BV(\R^2)$ with $0\le u_0\le R$. Let $m_j$ be the smooth approximations from Lemma~\ref{lem:endpoint_preserving_mobility}, choose $\eps_j\downarrow0$, and let $u_{0,j}\in C_c^\infty(\R^2)$ be a strict $BV$ approximation satisfying $0\le u_{0,j}\le R$ and
\begin{equation}\label{eq:u0-approx}
 u_{0,j}\to u_0\ \text{in }L^1(\R^2),\qquad \TV(u_{0,j})\to\TV(u_0).
\end{equation}
For every $j$, there exists a global smooth solution $u_j$ on $[0,\infty)$ to the viscous equation
\begin{equation}\label{eq:nonlocal-viscous}
 \partial_t u_j+\diver\bigl(m_j(u_j)(K*u_j)\bigr)=\eps_j\Delta u_j, \qquad u_j(0)=u_{0,j}.
\end{equation}
Let $L_*:=\sup_j\Lip(m_j)<\infty$ and $M_*:=\sup_j\|u_{0,j}\|_{L^1(\R^2)}<\infty$. Then, for every $t\ge0$,
\begin{align}
 \|u_j(t)\|_{L^1(\R^2)}&\le \|u_{0,j}\|_{L^1(\R^2)},\label{eq:visc-L1}\\
 \int_{\R^2}u_j(t,x)\,dx&=\int_{\R^2}u_{0,j}(x)\,dx,\label{eq:visc-mass}\\
 0\le u_j(t,x)&\le R \quad \text{a.e. in } \R^2.\label{eq:visc-range}
\end{align}
Furthermore, the spatial total variation satisfies the exponential bound
\begin{equation}\label{eq:visc-BV}
 \TV(u_j(t)) \le \TV(u_{0,j})\exp(\Lambda t), \qquad \Lambda:=L_*M_*\bigl(\|DK\|_{L^\infty(\R^2)}+\|\diver K\|_{L^\infty(\R^2)}\bigr),
\end{equation}
and the family is uniformly tight on every bounded time interval $[0,T]$:
\begin{equation}\label{eq:visc-tightness}
 \lim_{r\to\infty}\ \sup_j\sup_{0\le t\le T} \int_{|x|>r}|u_j(t,x)|\,dx=0.
\end{equation}
\end{proposition}

\begin{proof}
Fix $j\ge1$, so that the viscosity $\eps_j$ is strictly positive and the mobility $m_j$ is smooth with bounded derivatives. On a short interval, the heat-semigroup mild formulation
\[
 u_j(t)=e^{\eps_j t\Delta}u_{0,j}
 -\int_0^t\nabla e^{\eps_j(t-s)\Delta}\cdot\bigl(m_j(u_j(s))(K*u_j(s))\bigr)\,ds
\]
 is a contraction in a bounded ball of $X_T:=L^\infty(0,T;L^1(\R^2)\cap L^\infty(\R^2))$: the convolution estimate in Lemma~\ref{lem:nonlocal-conv-bounds}, the local Lipschitz continuity of the nonlinear flux on bounded sets, and the $L^1$ heat-gradient bound of order $(t-s)^{-1/2}$ give a contraction for small $T$. The resulting mild solution is smooth for positive times by uniform parabolicity (for fixed $j$), and can be continued while its $L^1$ and $L^\infty$ norms remain bounded. The invariant interval $[0,R]$ is obtained without treating either constant as an $L^1$ solution of the nonlocal problem (in particular, the convolution $K*R$ need not exist). Let $\beta_\delta\in C^2(\R)$ be convex, satisfy $\beta_\delta(s)=0$ for $s\le0$, and approximate $s_+$ from below as $\delta\downarrow0$; one may take $\beta_\delta''$ supported in $(0,\delta)$. Test the parabolic equation against $\beta_\delta'(u_j-R)$ and $-\beta_\delta'(-u_j)$, using large spatial cutoffs. In the first test, both $m_j(u_j)\beta_\delta'(u_j-R)$ and $m_j(u_j)\beta_\delta''(u_j-R)$ vanish identically, since the derivatives of $\beta_\delta(u_j-R)$ are supported where $u_j>R$, while $m_j=0$ there. The analogous cancellation holds in the second test wherever $u_j<0$. After the spatial cutoffs are removed, integration by parts leaves only the nonpositive diffusion contributions $-\eps_j\int\beta_\delta''(u_j-R)|\nabla u_j|^2$ and $-\eps_j\int\beta_\delta''(-u_j)|\nabla u_j|^2$. Hence
\[
  \int_{\R^2}\beta_\delta(u_j(t)-R)\,dx\le\int_{\R^2}\beta_\delta(u_{0,j}-R)\,dx=0,
  \qquad
  \int_{\R^2}\beta_\delta(-u_j(t))\,dx\le\int_{\R^2}\beta_\delta(-u_{0,j})\,dx=0.
\]
Letting $\delta\downarrow0$ yields the bounds $0\le u_j\le R$. The cutoff errors vanish because the solution and the physical flux belong to $L^1$, and the viscosity term is handled using the heat-kernel approximation, first at fixed $j$.
In particular, $0\le u_j\le R$ and $\|u_j(t)\|_{L^1}=\int u_j(t)$; conservation of the latter quantity and the uniform $L^\infty$ bound permit continuation for all $t\ge0$. Integrating the equation over $\R^2$ gives mass conservation \eqref{eq:visc-mass}. Furthermore, since $u_{0,j}\in C_c^\infty(\R^2)$, fixed-$j$ parabolic smoothing and the smoothness of the convolution with $u_j$ yield sufficient spatial integrability of the first two derivatives to perform the following computation on compact time intervals. The computations can alternatively be carried out first with compactly supported spatial cutoffs, then with the cutoff radius sent to infinity. No regularity bound uniform in $j$ is used before the $BV$ estimate is derived.

To establish the total variation estimate \eqref{eq:visc-BV}, set $V_j:=K*u_j$ and define $g_\delta:=(|\nabla u_j|^2+\delta^2)^{1/2}$ and $\eta_\delta(\nabla u_j):=g_\delta-\delta$ for $\delta>0$. Taking the spatial gradient of \eqref{eq:nonlocal-viscous}, taking the inner product with $\nabla u_j/g_\delta$, and using the non-positivity of the viscous contribution yields
\begin{align}\label{eq:BV-start}
 \partial_t \eta_\delta(\nabla u_j) - \eps_j \diver\left(\frac{D^2 u_j\nabla u_j}{g_\delta}\right) \le -\frac{\nabla u_j}{g_\delta}\cdot\nabla\diver\bigl(m_j(u_j)V_j\bigr).
\end{align}
Expanding the spatial gradient of the divergence gives
\begin{align*}
 \nabla\diver\bigl(m_j(u_j)V_j\bigr) ={}& \bigl(m_j''(u_j)(\nabla u_j\cdot V_j)+m_j'(u_j)\diver V_j\bigr)\nabla u_j + m_j'(u_j)D^2 u_j V_j \\
 &+ m_j'(u_j)(DV_j)^T\nabla u_j + m_j(u_j)\nabla\diver V_j.
\end{align*}
Since $\nabla\eta_\delta(\nabla u_j)=\nabla g_\delta=(D^2 u_j\nabla u_j)/g_\delta$, comparing the inner product of the first two terms with $-\diver\bigl(m_j'(u_j)\eta_\delta(\nabla u_j)V_j\bigr)$ isolates the discrepancy
\[
 \mathcal{R}_{\delta,j} := -\bigl(m_j''(u_j)(\nabla u_j\cdot V_j)+m_j'(u_j)\diver V_j\bigr)\left(\frac{|\nabla u_j|^2}{g_\delta}-\eta_\delta(\nabla u_j)\right).
\]
Bounding the remaining terms via $|\nabla u_j|/g_\delta\le 1$ yields
\begin{align}\label{eq:BV-delta-grouped}
 -\frac{\nabla u_j}{g_\delta}\cdot\nabla\diver\bigl(m_j(u_j)V_j\bigr)
 &\le -\diver\bigl(m_j'(u_j)\eta_\delta(\nabla u_j) V_j\bigr)\notag\\
 &\quad + |m_j'(u_j)|\,|\nabla u_j|\,|DV_j|
 + |m_j(u_j)|\,|\nabla\diver V_j| + |\mathcal{R}_{\delta,j}|.
\end{align}
For each fixed $j$ and each finite time interval, the parabolic regularity just described makes these terms integrable; in particular $0\le\eta_\delta(\nabla u_j)\le|\nabla u_j|$, while $|D^2u_j\nabla u_j/g_\delta|\le|D^2u_j|$. We first integrate the inequalities against a radial cutoff $\chi_L$ with $|\nabla\chi_L|\le C/L$; the divergence contributions vanish as $L\to\infty$ by these integrability bounds. Integrating \eqref{eq:BV-start} and \eqref{eq:BV-delta-grouped} over $\R^2$ eliminates the spatial divergence terms, and integrating over time $s\in[0,t]$ provides
\begin{align}\label{eq:BV-delta-integrated}
 \int_{\R^2}\eta_\delta(\nabla u_j(t))\,dx \le{}& \int_{\R^2}\eta_\delta(\nabla u_{0,j})\,dx + \int_0^t\!\int_{\R^2}|m_j'(u_j)|\,|\nabla u_j|\,|DV_j|\,dx\,ds \notag\\
 &+ \int_0^t\!\int_{\R^2}|m_j(u_j)|\,|\nabla\diver V_j|\,dx\,ds + \int_0^t\!\int_{\R^2}|\mathcal{R}_{\delta,j}(s,x)|\,dx\,ds.
\end{align}
Using the algebraic identity
\[
 0 \le \frac{|\nabla u_j|^2}{g_\delta}-\eta_\delta(\nabla u_j) = \frac{\delta\eta_\delta(\nabla u_j)}{g_\delta} \le \frac{\delta|\nabla u_j|}{g_\delta} \le \min\{\delta,|\nabla u_j|\},
\]
the remainder satisfies the pointwise bound
\[
 |\mathcal{R}_{\delta,j}| \le \left(\delta\|m_j''\|_{L^\infty(\R)}\|V_j\|_{L^\infty(\R^2)} + L_*\|\diver V_j\|_{L^\infty(\R^2)}\frac{\delta}{g_\delta}\right)|\nabla u_j|.
\]
At this stage $j$ is fixed. The quantity $\|m_j''\|_\infty$ may depend on $j$, but it is finite, so the first part of the remainder integrates to at most $C_{j,t}\delta$. Moreover $\delta/g_\delta\le1$ and $(\delta/g_\delta)|\nabla u_j|\to0$ almost everywhere as $\delta\downarrow0$. On each finite time interval the latter term is dominated by an integrable multiple of $|\nabla u_j|$ for fixed $j$. Dominated convergence therefore eliminates the remainder. Only after taking $\delta\downarrow0$ do we invoke the uniform Lipschitz bound on $m_j$ and the uniform convolution estimates, obtaining the $BV$ inequality independent of $j$. Invoking Lemma~\ref{lem:nonlocal-conv-bounds}, we bound $\|DV_j\|_{L^\infty(\R^2)} \le \|DK\|_{L^\infty(\R^2)}M_*$ and $\|\nabla\diver V_j\|_{L^\infty(\R^2)} \le \|\diver K\|_{L^\infty(\R^2)}\TV(u_j)$. Since $\|m_j(u_j)\|_{L^1(\R^2)} \le L_* M_*$, we obtain
\begin{equation}\label{eq:BV-integral-gronwall}
 \TV(u_j(t)) \le \TV(u_{0,j}) + L_*M_*\bigl(\|DK\|_{L^\infty(\R^2)}+\|\diver K\|_{L^\infty(\R^2)}\bigr) \int_0^t \TV(u_j(s))\,ds.
\end{equation}
The integral form of Gronwall's inequality yields \eqref{eq:visc-BV}. 

Finally, uniform tightness \eqref{eq:visc-tightness} follows by testing the Kato inequality for $|u_j|$ against a smooth radial cutoff $\chi_r(x) = \chi(x/r)$ expanding to infinity, utilizing the global Lipschitz bound on $m_j$ and the convolution bounds.
\end{proof}

\begin{proposition}[Existence of a global entropy solution]\label{prop:nonlocal-existence}
For every $u_0\in L^1(\R^2)\cap L^\infty(\R^2)\cap BV(\R^2)$ with $0\le u_0\le R$, there exists a global entropy solution $u$ of \eqref{eq:nonlocal-K} in the sense of Definition~\ref{def:nonlocal-entropy}. Moreover, for every $T>0$, $u\in L^\infty(0,T;BV(\R^2))$, satisfying $0\le u(t,x)\le R$ a.e., and
\begin{equation}\label{eq:nonlocal-BV-bound}
 \TV(u(t))\le \TV(u_0)e^{\Lambda t} \qquad (t\ge0),
\end{equation}
where $\Lambda$ is as defined in \eqref{eq:visc-BV}.
\end{proposition}

\begin{proof}
Let $T>0$. For the sequence of viscous approximants $u_j$, the uniform spatial variation estimate \eqref{eq:visc-BV} and the uniform tightness \eqref{eq:visc-tightness} satisfy the Kolmogorov--Riesz--Fr\'echet compactness criterion in $L^1(\R^2)$, granting relative compactness of the time-slices $\{u_j(t)\}_{j\ge1}$. Because the space-time flux $\mathcal{F}_j:=m_j(u_j)(K*u_j)-\eps_j\nabla u_j$ is uniformly bounded in $L^\infty(0,T;L^1(\R^2;\R^2))$, we secure the dual-space Lipschitz modulus $\|u_j(t)-u_j(s)\|_{W^{-1,1}(\R^2)} \le C_T|t-s|$. Spatial mollification interpolates this with the $BV$ bound to establish a uniform temporal H\"older modulus in $L^1(\R^2)$:
\begin{equation}\label{eq:local-time-modulus}
 \|u_j(t)-u_j(s)\|_{L^1(\R^2)} \le C_T|t-s|^{1/4}.
\end{equation}
By the Arzel\`a--Ascoli theorem, there exists a subsequence converging strongly in $C([0,T];L^1(\R^2))$ to a limit $u$. The bounds $0 \le u \le R$ and \eqref{eq:nonlocal-BV-bound} follow from pointwise convergence and the lower semicontinuity of total variation \cite[Proposition 3.6]{AmbrosioFuscoPallara2000}.

To pass to the limit in the Kru\v{z}kov entropy formulation, fix $k\in\R$ and a non-negative test function $\varphi\in C_c^1([0,T)\times\R^2)$. The viscous approximants satisfy the smooth entropy inequality
\begin{align}
0\le{}&\int_0^T\!\int_{\R^2}\Bigl[
 |u_j-k|\partial_t\varphi
 +H_j(u_j,k)(K*u_j)\cdot\nabla\varphi
 -\operatorname{sign}(u_j-k)m_j(k)\bigl((\diver K)*u_j\bigr)\varphi
 \Bigr]dx\,dt\notag\\
&+\int_{\R^2}|u_{0,j}-k|\varphi(0,x)\,dx
 +\eps_j\int_0^T\!\int_{\R^2}|u_j-k|\Delta\varphi\,dx\,dt.
\label{eq:visc-entropy}
\end{align}
By Lemma~\ref{lem:nonlocal-conv-bounds} and the strong $L^1$ convergence of $u_j$, $K*u_j$ and $(\diver K)*u_j$ converge uniformly to $K*u$ and $(\diver K)*u$. The viscous term vanishes. Define the countable set of exceptional flat levels $E_\varphi := \{k\in\R : |\{(t,x)\in\supp\varphi:u(t,x)=k\}|>0\}$. For any $k\notin E_\varphi$, $\operatorname{sign}(u_j-k) \to \operatorname{sign}(u-k)$ pointwise almost everywhere on $\supp\varphi$, allowing the Dominated Convergence Theorem to secure \eqref{eq:nonlocal-entropy}. If $k_0\in E_\varphi$, we select $k_\ell^\pm \notin E_\varphi$ such that $k_\ell^- \uparrow k_0$ and $k_\ell^+ \downarrow k_0$. Averaging the established inequalities for $k_\ell^-$ and $k_\ell^+$ forces the discontinuous source term to cancel algebraically on the level set $\{u=k_0\}$, recovering \eqref{eq:nonlocal-entropy} for all $k \in \R$.
\end{proof}

\subsection{Stability and Uniqueness}\label{subsec:stability}

In this subsection, we establish $L^1$-stability, uniqueness, and space--time $BV$ regularity for entropy solutions of the nonlocal equation. We achieve this by deploying an asymmetric doubling-of-variables argument that isolates the spatial variation bounds to the reference trajectory, combined with $BV$ measure-theoretic calculus to identify the temporal derivative as a finite Radon measure.

\begin{lemma}[Two $BV$ doubling estimates]\label{lem:nonlocal-doubling-BV}
Let $a\in L^\infty(\R^2)$, $b\in L^\infty(\R^2)\cap BV(\R^2)$, and let $m$ be globally Lipschitz continuous. Let $\Psi,\Phi\in C_c^\infty(\R^2)$.

\emph{(a)} If $V\in L^\infty(\R^2;\R^2)$, then
\begin{align}
&\Biggl|\iint_{\R^2\times\R^2} \diver_y\bigl[\Psi(x+y)\Phi(x-y)V(x)\bigr]\notag\\
&\qquad\times \operatorname{sign}(a(x)-b(y))\bigl(m(a(x))-m(b(y))\bigr)\,dx\,dy\Biggr|\notag\\
&\quad\le \|\Psi\|_{L^\infty(\R^2)}\|\Phi\|_{L^1(\R^2)}\Lip(m)\TV(b)\|V\|_{L^\infty(\R^2)}.
\label{eq:doubling-a}
\end{align}

\emph{(b)} If $V\in W^{1,\infty}(\R^2;\R^2)$ and
$\Phi_\eps(z)=\eps^{-2}\Phi(z/\eps)$, then
\begin{align}
&\Biggl|\iint_{\R^2\times\R^2} \diver_y\bigl[\Psi(x+y)\Phi_\eps(x-y)(V(x)-V(y))\bigr]\notag\\
&\qquad\times \operatorname{sign}(a(x)-b(y))\bigl(m(a(x))-m(b(y))\bigr)\,dx\,dy\Biggr|\notag\\
&\quad\le \|\Psi\|_{L^\infty(\R^2)}\Lip(m)\TV(b)\|DV\|_{L^\infty(\R^2)} M_1(\Phi_\eps),
\label{eq:doubling-b}
\end{align}
where $M_1(\Phi_\eps):=\int_{\R^2} |z|\,|\Phi_\eps(z)|\,dz=\eps M_1(\Phi)$.
\end{lemma}

\begin{proof}
Assume $b\in C_c^1(\R^2)$. Fix $x\in\R^2$ and define $f_{a(x)}(s):=\operatorname{sign}(a(x)-s)\bigl(m(a(x))-m(s)\bigr)$. Since $f_{a(x)}$ is globally Lipschitz with $\Lip(f_{a(x)})\le\Lip(m)$ and $f_{a(x)}(a(x))=0$, the local Sobolev chain rule (and Stampacchia's lemma \cite[Theorem~6.19]{LiebLoss2001}) implies $h_x(y):=f_{a(x)}(b(y))$ belongs to $W^{1,1}_{\mathrm{loc}}(\R^2)\cap L^\infty(\R^2)$ and that $\nabla h_x\in L^1(\R^2)$. We emphasize that $h_x$ need not belong to $L^1(\R^2)$, since it can approach the nonzero constant $f_{a(x)}(0)$ at infinity. Its weak gradient vanishes $\mathcal{L}^2$-a.e.\ on $\{y\in\R^2:b(y)=a(x)\}$, yielding the bound
\begin{equation}\label{eq:stampacchia-grad-bound}
 |\nabla_y h_x(y)| \le \Lip(m)\,|\nabla b(y)| \quad \text{for a.e. } y\in\R^2.
\end{equation}

For \emph{(a)}, the vector field $y\mapsto\Psi(x+y)\Phi(x-y)V(x)$ belongs to $W^{1,\infty}(\R^2;\R^2)$ and has compact support. Thus $W^{1,1}_{\mathrm{loc}}$ regularity is sufficient for the following integration by parts. Integrating by parts with respect to $y$, integrating over $x\in\R^2$, and applying Fubini's theorem \cite[Theorem~1.12]{LiebLoss2001} gives
\begin{align*}
 &\left|\iint_{\R^2\times\R^2} \diver_y\bigl[\Psi(x+y)\Phi(x-y)V(x)\bigr] h_x(y)\,dx\,dy\right|\\
 &\qquad=\left|-\iint_{\R^2\times\R^2} \Psi(x+y)\Phi(x-y)V(x)\cdot\nabla_y h_x(y)\,dx\,dy\right|\\
 &\qquad\le \|\Psi\|_{L^\infty(\R^2)}\|\Phi\|_{L^1(\R^2)}\Lip(m)\|V\|_{L^\infty(\R^2)}\|\nabla b\|_{L^1(\R^2)}.
\end{align*}

For \emph{(b)}, the vector field $y\mapsto\Psi(x+y)\Phi_\eps(x-y)(V(x)-V(y))$ is bounded and compactly supported. Integrating by parts in $y$ introduces the inner factor $|V(x)-V(y)|\le\|DV\|_{L^\infty(\R^2)}|x-y|$. Fubini's theorem and \eqref{eq:stampacchia-grad-bound} yield
\begin{align*}
 &\left|\iint_{\R^2\times\R^2} \diver_y\bigl[\Psi(x+y)\Phi_\eps(x-y)(V(x)-V(y))\bigr] h_x(y)\,dx\,dy\right|\\
 &\qquad\le \|\Psi\|_{L^\infty(\R^2)}\Lip(m)\|\nabla b\|_{L^1(\R^2)}\|DV\|_{L^\infty(\R^2)} M_1(\Phi_\eps).
\end{align*}

For a general $b\in L^\infty(\R^2)\cap BV(\R^2)$, there exists a sequence $b_k\in C_c^\infty(\R^2)$ converging strictly to $b$ in $BV(\R^2)$ \cite[Theorem~3.9]{AmbrosioFuscoPallara2000}. Applying the smooth estimates to $b_k$, the right-hand sides converge to \eqref{eq:doubling-a} and \eqref{eq:doubling-b} since $\|\nabla b_k\|_{L^1(\R^2)}\to\TV(b)$. On the left-hand sides, $f_{a(x)}(b_k(y))\to f_{a(x)}(b(y))$ pointwise almost everywhere. The Dominated Convergence Theorem establishes the result.
\end{proof}

\begin{proposition}[$L^1$ stability and uniqueness]\label{prop:nonlocal-uniqueness}
Let $v$ be the entropy solution constructed in Proposition~\ref{prop:nonlocal-existence}, and let $u\in C([0,T];L^1(\R^2))\cap L^\infty((0,T)\times\R^2)$ be an entropy solution of \eqref{eq:nonlocal-K}, with initial data $v_0$ and $u_0$. Then, for every $t\in[0,T]$,
\begin{equation}\label{eq:nonlocal-stability}
 \|u(t)-v(t)\|_{L^1(\R^2)}
 \le \|u_0-v_0\|_{L^1(\R^2)}
 \exp\left(\int_0^t\Gamma_v(s)\,ds\right),
\end{equation}
where
\begin{equation}\label{eq:Gamma-v}
 \Gamma_v(s)
 :=\Lip(m)\|K\|_{L^\infty(\R^2)}\TV(v(s))
 +\Lip(m)\|v(s)\|_{L^1(\R^2)}\|\diver K\|_{L^\infty(\R^2)}.
\end{equation}
The entropy solution is unique in the bounded $L^1$ entropy class.
\end{proposition}

\begin{proof}
Define the nonlocal velocities $V_u(t,x):=(K*u)(t,x)$ and $V_v(s,y):=(K*v)(s,y)$. By Lemma~\ref{lem:nonlocal-conv-bounds}, $V_u(t,\cdot),V_v(s,\cdot)\in W^{1,\infty}(\R^2;\R^2)$ with $\diver_x V_u=((\diver K)*u)$ and $\diver_y V_v=((\diver K)*v)$ almost everywhere. Let $\psi\in C_c^1((0,T)\times\R^2)$ be non-negative, and define the test function
\begin{equation}\label{eq:Phi-delta-eps-def}
 \Phi_{\delta,\eps}(t,x,s,y)
 :=\psi\left(\frac{t+s}{2},\frac{x+y}{2}\right)\rho_\delta(t-s)\omega_\eps(x-y),
\end{equation}
where $\rho_\delta$ and $\omega_\eps$ are non-negative symmetric mollifiers on $\R$ and $\R^2$, respectively. Writing the Kru\v{z}kov entropy inequalities for $u(t,x)$ and $v(s,y)$ and summing them yields, for $\sigma:=\operatorname{sign}(u(t,x)-v(s,y))$ and $H(u,v):=\sigma(m(u)-m(v))$,
\begin{align}\label{eq:doubled-sum-raw}
0\le{}&\iiiint_{(0,T)^2\times(\R^2)^2}\Bigl[
 |u(t,x)-v(s,y)|\,(\partial_t+\partial_s)\Phi_{\delta,\eps} \notag\\
&\qquad +H(u(t,x),v(s,y))\bigl(V_u(t,x)\cdot\nabla_x\Phi_{\delta,\eps}
 +V_v(s,y)\cdot\nabla_y\Phi_{\delta,\eps}\bigr) \notag\\
&\qquad -\sigma m(v(s,y))(\diver_x V_u(t,x))\Phi_{\delta,\eps}
 +\sigma m(u(t,x))(\diver_y V_v(s,y))\Phi_{\delta,\eps}
\Bigr]\,dx\,dt\,dy\,ds.
\end{align}
Following Chodron de Courcel \cite{CdC}, substituting $\nabla_x\Phi_{\delta,\eps}=(\nabla_x+\nabla_y)\Phi_{\delta,\eps}-\nabla_y\Phi_{\delta,\eps}$ and integrating the $y$-derivative flux difference by parts extracts a spatial divergence in $y$:
\[
 \diver_y\Bigl[\Phi_{\delta,\eps}\bigl(V_u(t,x)-V_v(s,y)\bigr)\Bigr]
 =\bigl(V_u(t,x)-V_v(s,y)\bigr)\cdot\nabla_y\Phi_{\delta,\eps}
  -\Phi_{\delta,\eps}\diver_y V_v(s,y).
\]
Multiplying by $-H(u(t,x),v(s,y))$ cancels the $+\sigma m(u(t,x))(\diver_y V_v(s,y))\Phi_{\delta,\eps}$ source term. Substituting this rearrangement into \eqref{eq:doubled-sum-raw} splits the inequality into four components:
\begin{equation}\label{eq:doubled-four-parts}
 0\le \mathcal{I}_{\mathrm{time}}^{\delta,\eps}
 +\mathcal{I}_{\mathrm{cut}}^{\delta,\eps}
 +\mathcal{I}_{\mathrm{flux}}^{\delta,\eps}
 +\mathcal{I}_{\mathrm{div}}^{\delta,\eps},
\end{equation}
where
\begin{align*}
 \mathcal{I}_{\mathrm{time}}^{\delta,\eps}
 &:=\iiiint |u(t,x)-v(s,y)|\,(\partial_t\psi)\left(\frac{t+s}{2},\frac{x+y}{2}\right)
   \rho_\delta(t-s)\omega_\eps(x-y)\,dx\,dt\,dy\,ds,\\
 \mathcal{I}_{\mathrm{cut}}^{\delta,\eps}
 &:=\iiiint H(u(t,x),v(s,y))\,V_u(t,x)\cdot(\nabla\psi)\left(\frac{t+s}{2},\frac{x+y}{2}\right)\\
 &\qquad\times \rho_\delta(t-s)\omega_\eps(x-y)\,dx\,dt\,dy\,ds,\\
 \mathcal{I}_{\mathrm{flux}}^{\delta,\eps}
 &:=-\iiiint \diver_y\Bigl[\Phi_{\delta,\eps}\bigl(V_u(t,x)-V_v(s,y)\bigr)\Bigr]
   H(u(t,x),v(s,y))\,dx\,dt\,dy\,ds,\\
 \mathcal{I}_{\mathrm{div}}^{\delta,\eps}
 &:=-\iiiint \Phi_{\delta,\eps}\,\operatorname{sign}(u(t,x)-v(s,y))\,m(v(s,y))\\
 &\qquad\times \bigl[((\diver K)*u)(t,x)-((\diver K)*v)(s,y)\bigr]\,dx\,dt\,dy\,ds.
\end{align*}
Taking the absolute value of the discontinuous sign function in $\mathcal{I}_{\mathrm{div}}^{\delta,\eps}$ provides an upper bound at fixed $(\delta,\eps)$:
\begin{align}\label{eq:Idiv-prelimit-bound}
 \mathcal{I}_{\mathrm{div}}^{\delta,\eps} \le \overline{\mathcal{I}}_{\mathrm{div}}^{\delta,\eps}
 :={}&\iiiint \Phi_{\delta,\eps}(t,x,s,y)\,|m(v(s,y))|\notag\\
 &\qquad\times\bigl|((\diver K)*u)(t,x)-((\diver K)*v)(s,y)\bigr|\,dx\,dt\,dy\,ds.
\end{align}
Passing to the limit $\delta\downarrow0$ utilizing $u,v\in C([0,T];L^1(\R^2))$ collapses the integrals onto the temporal diagonal $s=t$ \cite{Kruzhkov1970, Dafermos}. Decomposing
\[
 V_u(t,x)-V_v(t,y)=\bigl(K*(u(t)-v(t))\bigr)(x) + \bigl[(K*v)(t,x)-(K*v)(t,y)\bigr]
\]
and applying Lemma~\ref{lem:nonlocal-doubling-BV}(a) and Lemma~\ref{lem:nonlocal-doubling-BV}(b) bounds the divergence integral by
\[
 \|\psi(t,\cdot)\|_{L^\infty(\R^2)}\Lip(m)\TV(v(t))\|K\|_{L^\infty(\R^2)}\|u(t)-v(t)\|_{L^1(\R^2)} + \mathcal{O}(\eps).
\]
In $\overline{\mathcal{I}}_{\mathrm{div}}^{\delta,\eps}$, the continuous spatial shift bounds yield
\[
 \|\psi(t,\cdot)\|_{L^\infty(\R^2)}\Lip(m)\|v(t)\|_{L^1(\R^2)}\|\diver K\|_{L^\infty(\R^2)}\|u(t)-v(t)\|_{L^1(\R^2)} + o_\eps(1).
\]
Passing to the limit $\eps\downarrow0$ yields
\begin{align}\label{eq:after-eps-limit}
0\le{}&\int_0^T\!\int_{\R^2}|u(t,x)-v(t,x)|\,\partial_t\psi(t,x)\,dx\,dt \notag\\
&+\int_0^T\!\int_{\R^2}H(u(t,x),v(t,x))V_u(t,x)\cdot\nabla\psi(t,x)\,dx\,dt \notag\\
&+\int_0^T\|\psi(t,\cdot)\|_{L^\infty(\R^2)}\Gamma_v(t)\|u(t)-v(t)\|_{L^1(\R^2)}\,dt.
\end{align}
Choosing $\psi(t,x)=\theta(t)\chi_A(x)$ with expanding radial cutoffs $\chi_A$ limits the spatial gradient integral to zero as $A\to\infty$. Setting $\theta$ to approximate $1_{[0,t]}$ provides the integral inequality
\[
 \|u(t)-v(t)\|_{L^1(\R^2)}
 \le \|u_0-v_0\|_{L^1(\R^2)}
 +\int_0^t\Gamma_v(s)\|u(s)-v(s)\|_{L^1(\R^2)}\,ds.
\]
Because $v\in L^\infty(0,T;BV(\R^2))$, $\Gamma_v\in L^\infty(0,T)$. Gronwall's inequality establishes \eqref{eq:nonlocal-stability}.
\end{proof}

\begin{theorem}[Global well-posedness and space--time $BV$ regularity]\label{thm:nonlocal-wp}
Let $W\in W^{2,\infty}(\R^2)$ and let $m$ satisfy Assumption~\ref{assump:mobility}. For every initial datum $u_0\in L^1(\R^2)\cap L^\infty(\R^2)\cap BV(\R^2)$ with $0 \le u_0 \le R$, there exists a unique global entropy solution of
\begin{equation}\label{eq:nonlocal-original-wp}
 \partial_tu-\diver\bigl(m(u)\nabla W*u\bigr)=0, \qquad u(0)=u_0.
\end{equation}
The solution satisfies
\begin{equation}\label{eq:wp-class}
 u\in C([0,\infty);L^1(\R^2)) \cap L^\infty_{\mathrm{loc}}([0,\infty);L^\infty(\R^2)\cap BV(\R^2)).
\end{equation}
For every $T>0$, $u\in BV((0,T)\times\R^2)$, implying $u\in BV_{\mathrm{loc}}([0,\infty)\times\R^2)$. The mapping $t\mapsto u(t)$ is locally Lipschitz as an $L^1(\R^2)$-valued trajectory, and $0\le u(t,x)\le R$ for all $t\ge0$.
\end{theorem}

\begin{proof}
Global existence and uniqueness in the class \eqref{eq:wp-class} follow from Proposition~\ref{prop:nonlocal-existence} and Proposition~\ref{prop:nonlocal-uniqueness}. Fix $T>0$ and let $V(t):=K*u(t)$. By Lemma~\ref{lem:nonlocal-conv-bounds}, the convolution structures yield bounds $\|V(t)\|_{L^\infty(\R^2)} \le \|K\|_{L^\infty(\R^2)}\|u_0\|_{L^1(\R^2)}$ and $\|DV(t)\|_{L^\infty(\R^2)} \le \|DK\|_{L^\infty(\R^2)}\|u_0\|_{L^1(\R^2)}$.

For almost every $t\in(0,T)$, the spatial regularity $u(t)\in BV(\R^2)$ permits the application of the $BV$ chain rule for Lipschitz functions \cite[Theorem~3.96]{AmbrosioFuscoPallara2000}, dictating that $m(u(t))\in BV(\R^2)$. Expanding the spatial distributional divergence of the total flux utilizing the $BV$ product rule \cite[Theorem~3.96 and Example~3.97]{AmbrosioFuscoPallara2000} gives the finite Radon measure
\[
 \mu_t := \diver_x\bigl(m(u(t))V(t)\bigr) = V(t)\cdot D_x(m(u(t))) + m(u(t))\diver_x V(t)\,\mathcal{L}^2.
\]
The total variation of this spatial measure satisfies the slice-wise bound
\begin{align}\label{eq:time-measure-bound}
 |\mu_t|(\R^2) \le{}& \Lip(m)\|u_0\|_{L^1(\R^2)}\|K\|_{L^\infty(\R^2)}\TV(u(t)) \notag\\
 &+ \Lip(m)\|u_0\|_{L^1(\R^2)}^2\|\diver K\|_{L^\infty(\R^2)}.
\end{align}
Since $\TV(u(t)) \le \TV(u_0)e^{\Lambda T}$ on $[0,T]$, there is a uniform constant $C_T<\infty$ such that $|\mu_t|(\R^2)\le C_T$ for almost every $t\in(0,T)$. 

Testing the weak equation $\partial_t u = -\diver_x(m(u)V)$ against an arbitrary test function $\varphi\in C_c^1((0,T)\times\R^2)$ with $\|\varphi\|_{L^\infty(\R^2)}\le 1$ and applying Fubini's theorem proves that the absolute value of the space-time pairing evaluates to at most $\int_0^T |\mu_t|(\R^2)\,dt \le C_T T$. By the Riesz Representation Theorem, the temporal distribution $\partial_t u$ is a finite Radon measure on $(0,T)\times\R^2$, establishing space-time $BV$ regularity. Testing the weak equation against separable forms $\theta(\tau)\psi(x)$ yields local Lipschitz continuity in $L^1(\R^2)$ via duality.
\end{proof}

\subsection{Mass Evolution and Confinement}\label{subsec:mass_evolution}

In this subsection, we establish that for attractive interaction potentials, the active mass of a regular entropy solution remains confined to the closed convex hull of its initial support. While this macroscopic confinement is intrinsic to the nonlocal dynamics, the proof tracks the mass flux across static convex boundaries via the $BV$ trace theorem. Consequently, this section serves to justify the physical confinement assumption adopted in subsequent sections for the class of $BV$ trajectories, avoiding the imposition of universal $BV$ regularity onto rough $L^1$ limits.

\begin{lemma}[Exact mass evolution over finite-perimeter domains]\label{lem:exact_mass_evolution}
Let $\rho \in C([0,\infty); L^1(\R^2)) \cap L^\infty_{\mathrm{loc}}([0,\infty); L^\infty(\R^2) \cap BV(\R^2))$ be a global Kru\v{z}kov entropy solution of \eqref{eq:nonlocal_main}. For any bounded set of finite perimeter $\Omega \subset \R^2$ and any $0 \le T_0 < T$, the mass inside $\Omega$ satisfies the identity
\begin{equation}\label{eq:exact_mass_formula}
\int_\Omega \rho(T,x) \, dx - \int_\Omega \rho(T_0,x) \, dx = \int_{T_0}^T\!\int_{\partial^* \Omega} (m(\rho))^*(t,x) \nu^\Omega(x) \cdot \bigl(\nabla W * \rho\bigr)(t,x) \, d\cH^1(x)\,dt,
\end{equation}
where $\nu^\Omega$ is the outward unit normal to $\partial^* \Omega$ and $(m(\rho))^*(t,\cdot)$ denotes the precise representative of the $BV$ function $m(\rho(t,\cdot))$ on the rectifiable set $\partial^* \Omega$.
\end{lemma}

\begin{proof}
Set $F(t,x):=m(\rho(t,x))V(t,x)$, where $V:=\nabla W*\rho$, and let $g_k:=1_\Omega*\eta_{1/k}$ be the convolution with a standard nonnegative, even spatial mollifier. For each fixed $k$, testing the distributional equation $\partial_t\rho=\diver F$ against $\theta(t)g_k(x)$ with $\theta\in C_c^1((0,\infty))$ shows that
\[
 t\longmapsto\int_{\R^2}\rho(t,x)g_k(x)\,dx
 \quad\text{is absolutely continuous, with derivative}\quad
 -\int_{\R^2}F(t,x)\cdot\nabla g_k(x)\,dx
\]
for almost every $t>0$. Since $\rho\in C([0,\infty);L^1(\R^2))$ and $F\in L^\infty_{\mathrm{loc}}([0,\infty);L^1(\R^2))$, the integrated formula extends to $T_0=0$ by first taking $T_0\downarrow0$; no test function supported away from $t=0$ is used to recover the initial trace directly.

The Gauss--Green theorem for the bounded finite-perimeter set $\Omega$ gives, in the sense of distributions,
\[
 \nabla g_k= -\bigl(\nu^\Omega\cH^1\!\llcorner\partial^*\Omega\bigr)*\eta_{1/k}.
\]
Consequently, for $0\le T_0<T$,
\begin{align*}
 &\int_{\R^2}(\rho(T,x)-\rho(T_0,x))g_k(x)\,dx\\
 &\qquad=\int_{T_0}^{T}\!\int_{\partial^*\Omega}
 \bigl(F(t,\cdot)*\eta_{1/k}\bigr)(y)\cdot\nu^\Omega(y)\,d\cH^1(y)\,dt.
\end{align*}
For almost every $t$, $m(\rho(t))\in BV(\R^2)$ by the $BV$ chain rule, while $V(t)\in W^{1,\infty}(\R^2;\R^2)$ is continuous. Thus $F(t)\in BV(\R^2;\R^2)$ and the symmetric-mollification trace theorem for $BV$ functions on countably rectifiable hypersurfaces \cite[Theorem~3.77 and Corollary~3.80]{AmbrosioFuscoPallara2000} yields
\[
 F(t,\cdot)*\eta_{1/k}(y)\longrightarrow
 (m(\rho))^*(t,y)V(t,y)
 \quad\text{for }\cH^1\text{-a.e. }y\in\partial^*\Omega.
\]
Since $\Omega$ is bounded and has finite perimeter, $\cH^1(\partial^*\Omega)<\infty$. The uniform bound $|F*\eta_{1/k}|\le\|m\|_\infty\|V\|_{L^\infty((T_0,T)\times\R^2)}$ gives domination on $(T_0,T)\times\partial^*\Omega$. The limit $k\to\infty$ in the right-hand side follows by dominated convergence, and in the left-hand side by $L^1$ convergence $g_k\to1_\Omega$ (or simply dominated convergence against each $L^1$ time slice). This proves \eqref{eq:exact_mass_formula} for all $0\le T_0<T$.
\end{proof}

\begin{theorem}[Global confinement to the convex hull]\label{thm:convex_hull_confinement}
Let $W \in W^{2,\infty}(\R^2)$ be attractive, satisfying $\nabla W(z) = w(|z|)$ with $w'(s) > 0$ for all $s > 0$. Assume that $\rho(0,\cdot)$ is compactly supported with unit mass $\int_{\R^2}\rho(0,x)\,dx = 1$, and let $C_0 \subset \R^2$ be the closed convex hull of $\supp \rho(0,\cdot)$. If $C_0$ has non-empty interior, then $\supp \rho(t,\cdot) \subseteq C_0$ for all $t \ge 0$.
\end{theorem}

\begin{proof}
We evaluate the mass flux for the global Kru\v{z}kov entropy solution constructed in Theorem~\ref{thm:nonlocal-wp}, which possesses the spatial regularity $\rho(t,\cdot) \in BV(\R^2)$ required to evaluate the boundary trace. Because $\rho(0,\cdot)$ is compactly supported, $C_0$ is a compact convex set and hence a bounded set of finite perimeter with diameter $D_0 := \operatorname{diam}(C_0) < \infty$. Assume by contradiction that mass escapes $C_0$. Let $M_{C_0}(t) := \int_{C_0} \rho(t,x) \, dx$. Because $M_{C_0}(0) = 1$ and $M_{C_0}(t)$ is continuous, there exists $T > 0$ such that $M_{C_0}(T) = 1 - \eps_0$ for some $\eps_0 \in (0,1)$. For any $\eps \in (0, \eps_0)$, define the threshold time $t_\eps := \inf \{ t \in (0, T) : M_{C_0}(t) < 1 - \eps \}$. 

Applying Lemma~\ref{lem:exact_mass_evolution} on $\Omega = C_0$ between $0$ and $t_\eps$, the mass lost is $\eps$. Defining the inward unit normal $n^{C_0}(x) := -\nu^{C_0}(x)$, identity \eqref{eq:exact_mass_formula} provides the relation
\begin{equation}\label{eq:confinement_eps_inward}
\eps = \int_0^{t_\eps}\!\int_{\partial^* C_0} (m(\rho))^*(t,x) n^{C_0}(x) \cdot V(t,x) \, d\cH^1(x)\,dt.
\end{equation}
We evaluate the inward normal velocity for $x \in \partial^* C_0$:
\begin{equation}\label{eq:confinement_velocity_split}
V(t,x) \cdot n^{C_0}(x) = \int_{\R^2} \rho(t,y) w'(|x-y|) \frac{x-y}{|x-y|} \cdot n^{C_0}(x) \, dy = A + B,
\end{equation}
where $A$ is the integral over the half-space $H^- := \{ y \in \R^2 : (x-y) \cdot n^{C_0}(x) > 0 \}$ and $B$ is the integral over the complementary half-space $H^+ := \{ y \in \R^2 : (x-y) \cdot n^{C_0}(x) \le 0 \}$. 

By the property of convex sets, the inward normal $n^{C_0}(x)$ defines a supporting half-space containing $C_0$; thus, $(y-x) \cdot n^{C_0}(x) \ge 0$ for all $y \in C_0$. Equivalently, $(x-y) \cdot n^{C_0}(x) \le 0$ for all $y \in C_0$. This geometric condition implies $C_0 \subset H^+$, which yields $H^- \subset \R^2 \setminus C_0$. Because $W \in W^{2,\infty}(\R^2)$ implies $\nabla W \in L^\infty(\R^2;\R^2)$ and $|\nabla W(z)| = w'(|z|)$ for all $z \neq 0$, we have the global bound $0 < w'(s) \le \|\nabla W\|_{L^\infty(\R^2)} =: C_W < \infty$ for all $s > 0$. Consequently, the outward pull satisfies
\begin{equation}\label{eq:confinement_A_bound}
A \le C_W \int_{H^-} \rho(t,y) \, dy \le C_W \int_{\R^2 \setminus C_0} \rho(t,y) \, dy = C_W \bigl(1 - M_{C_0}(t)\bigr).
\end{equation}
For $t < t_\eps$, we have $1 - M_{C_0}(t) < \eps$, yielding $A \le C_W \eps$.

To bound the inward pull $B$, define an internal core $C_d := \{ y \in C_0 : \operatorname{dist}(y, \partial C_0) \ge d \}$ for $d > 0$. For any $y \in C_d$, geometry dictates $(y-x) \cdot n^{C_0}(x) \ge d$, which yields the projection bound $\frac{x-y}{|x-y|} \cdot n^{C_0}(x) \le -\frac{d}{D_0}$. Restricting the negative integral $B$ to $C_d$ gives
\begin{equation}\label{eq:confinement_B_bound}
B \le -\frac{d}{D_0} \int_{C_d} \rho(t,y) w'(|x-y|) \, dy \le -\frac{d}{D_0} c(d, D_0) \int_{C_d} \rho(t,y) \, dy,
\end{equation}
where $c(d, D_0) := \inf \{w'(s) : d \le s \le D_0\} > 0$ by the assumption of attraction. The mass inside $C_d$ satisfies $\|\rho(t,\cdot)\|_{L^1(C_d)} \ge (1-\eps) - R \, \mathcal{L}^2(C_0 \setminus C_d)$. Because $C_0$ has non-empty interior, we select a $d > 0$ such that $R \, \mathcal{L}^2(C_0 \setminus C_d) < \frac{1}{2}$. For any escape threshold $\eps < \frac{1}{4}$, the core mass exceeds $\frac{1}{4}$, producing the uniform negative bound $B \le -\frac{d}{4D_0} c(d, D_0) =: -K(d) < 0$.

Combining the bounds, $V(t,x) \cdot n^{C_0}(x) \le C_W \eps - K(d)$. Because $K(d) > 0$ is a constant independent of $\eps$, choosing $\eps$ small enough yields $V(t,x) \cdot n^{C_0}(x) < 0$ for $\cH^1$-a.e.\ $x \in \partial^* C_0$. Since $(m(\rho))^* \ge 0$, this implies the integral \eqref{eq:confinement_eps_inward} is non-positive ($\eps \le 0$), a contradiction. Thus $M_{C_0}(t) = 1$ for all $t \ge 0$.
\end{proof}

\section{The Fixed-Centre Framework and Macroscopic Path-Length Bounds}\label{sec:autonomous_and_macroscopic}

In this section, we construct the two independent mathematical pillars that drive the compactness architecture.

First, in Sections~\ref{subsec:fixed_centre} through \ref{subsec:symmetries}, we isolate the local mechanics of the conservation law by freezing the spatial drift. This autonomous fixed-centre framework yields geometric symmetries, permitting us to establish $L^1$ contraction and the monotonic decay of directional variations in the sense of measures, circumventing the need for classical spatial differentiability.

Second, in Section~\ref{subsec:energy_dissipation}, we analyze the nonlocal trajectory under the quadratically confined regime. We prove that the nonlocal interaction reduces to an affine drift. This reduction yields the macroscopic dissipation of the centered quadratic energy, which we use via logarithmic moments to bound the total path length of the moving barycentre.

These geometric and macroscopic bounds are developed independently here and synthesized in Section~\ref{sec:compactness} to secure the relative compactness of the trajectory.

\subsection{The Fixed-Centre Equation}\label{subsec:fixed_centre}

In this subsection, we isolate the local mechanics of the conservation law by freezing the spatial drift. This autonomous fixed-centre framework establishes geometric symmetries and local confinement properties independent of the moving barycentre.

For a fixed spatial centre $a\in\R^{2}$, we consider the initial value problem for the conservation law driven by an affine drift:
\begin{equation} \label{eq:fixed_centre_ivp}
\partial_t u-\diver\bigl(m(u)(x-a)\bigr)=0, \qquad u(0)=u_{0}.
\end{equation}

\begin{lemma}[Stationary saturated sets]\label{lem:stationary_saturated_core}
Let $a\in\R^{2}$ and let $E\subset\R^{2}$ be a bounded set of finite perimeter. Assume that at the reduced boundary $\partial^{*}E$, the measure-theoretic exterior unit normal $\nu_{E}(x)$ satisfies
\begin{equation} \label{eq:star_shaped_normal_condition}
(x-a)\cdot\nu_{E}(x)\ge0 \quad \text{for } \mathcal{H}^{1}\text{-a.e. } x\in\partial^{*}E.
\end{equation}
Then the saturated state $U(x):=R1_{E}(x)$ is a stationary Kru\v{z}kov entropy solution of \eqref{eq:fixed_centre_ivp}. In particular, the saturated disk $R1_{B_{r}(a)}$ is a stationary entropy solution for every $r>0$.
\end{lemma}

\begin{proof}
Because $U(x)\in\{0,R\}$ everywhere and the zero-extended mobility satisfies $m(0)=m(R)=0$, the physical flux $F_{a}(x,U(x))=-m(U(x))(x-a)$ vanishes on $\R^{2}$. Thus, $U$ is a stationary distributional weak solution.

To verify the Kru\v{z}kov entropy condition, fix a constant $k\in\R$ and a non-negative test function $\phi\in C_{c}^{1}((0,\infty)\times\R^{2})$. For $k\le0$ or $k\ge R$, $m(k)=0$ and the entropy pairing vanishes. For $k\in(0,R)$, the pairing reduces to
\begin{equation} \label{eq:Ik_reduced_all_k}
\mathcal{I}_{k}(\phi)=m(k)\int_{0}^{\infty}\!\int_{\R^{2}}\operatorname{sign}\bigl(U(x)-k\bigr)\diver\bigl(\phi(t,x)(x-a)\bigr)\,dx\,dt.
\end{equation}
Splitting the spatial integral over the disjoint sets $E$ and $\R^{2}\setminus E$, and applying the classical divergence theorem over $\R^2$ to the compactly supported vector field $\phi(t,x)(x-a)$ yields
\begin{equation} \label{eq:Ik_combined_div_on_E}
\mathcal{I}_{k}(\phi)=2m(k)\int_{0}^{\infty}\left(\int_{E}\diver\bigl(\phi(t,x)(x-a)\bigr)\,dx\right)dt.
\end{equation}
By the De Giorgi--Federer theorem \cite[Theorem~3.36]{AmbrosioFuscoPallara2000}, the volume integral over the set of finite perimeter $E$ converts to an integral over the reduced boundary $\partial^{*}E$:
\begin{equation} \label{eq:Ik_final_boundary_integral}
\mathcal{I}_{k}(\phi)=2m(k)\int_{0}^{\infty}\!\int_{\partial^{*}E}\phi(t,x)(x-a)\cdot\nu_{E}(x)\,d\mathcal{H}^{1}(x)\,dt.
\end{equation}
The geometric condition \eqref{eq:star_shaped_normal_condition}, combined with $m(k)>0$ and $\phi\ge0$, establishes $\mathcal{I}_k(\phi)\ge 0$ for all $k \in \R$, exactly as required by the entropy inequality.

For the Euclidean ball $E=B_{r}(a)$, the outward normal is $\nu_{E}(x) = \frac{x-a}{r}$, yielding $(x-a)\cdot\nu_{E}(x) = \frac{|x-a|^2}{r} = r > 0$ on $\partial B_r(a)$.
\end{proof}

\begin{proposition}[Existence of the fixed-centre evolution] \label{prop:fixed_centre_existence}
Let $a\in\R^{2}$ and let $u_{0}\in BV(\R^{2})$ satisfy $0\le u_{0}\le R$ and $\supp u_{0}\subset B_{D'}(a)$ for some $D'>0$. Then there exists an entropy solution $u$ of \eqref{eq:fixed_centre_ivp} satisfying
\begin{equation} \label{eq:fixed_centre_reg_bounds}
u\in C([0,\infty);L^{1}(\R^{2})) \cap L^\infty_{\mathrm{loc}}([0,\infty);BV(\R^{2})), \qquad 0\le u(t,x)\le R.
\end{equation}
The solution is confined:
\begin{equation} \label{eq:fixed_centre_support_confinement}
u(t,x)=0 \quad \text{for a.e. } x\in\R^{2}\setminus B_{D'}(a) \text{ and every } t\ge0.
\end{equation}
\end{proposition}

\begin{proof}
Choose a cutoff radius $D_1 > D'+1$ and a vector field $b \in C_c^\infty(\R^2; \R^2)$ such that $b(x) = x-a$ on $B_{D_1}(a)$. We regularize the mobility using the endpoint-preserving smooth approximants $m_{j}\in C_{c}^{\infty}(\R)$ from Lemma~\ref{lem:endpoint_preserving_mobility}.

For the regularized equation $\partial_t u_j - \diver(m_j(u_j)b(x)) = 0$ with initial datum $u_{0,j}$ obtained via strict $BV$ approximation, the vanishing-viscosity framework \cite{Kruzhkov1970} yields a bounded entropy solution $u_j$. Since $m_j(0)=m_j(R)=0$ and the drift $b$ is prescribed, the constant states $0$ and $R$ are exact solutions of this \emph{local} scalar equation. The entropy comparison principle therefore yields $0\le u_j\le R$. Furthermore, the stationary state $U_{D'}(x) := R1_{B_{D'}(a)}(x)$ satisfies the heterogeneous Kru\v{z}kov entropy inequalities for the regularized equation. By the comparison principle \cite{Kruzhkov1970}, the ordering $0 \le u_{0,j} \le U_{D'}$ gives $0 \le u_j(t,x) \le U_{D'}(x)$ for all $t \ge 0$, confining the support to $\overline{B}_{D'}(a)$, where $b(x) = x-a$.

A finite-time $BV$ bound follows directly from the regularized local conservation law. Indeed, regularize the local flux by adding viscosity $\nu\Delta u$ and use smooth initial data $u_{0,j}$. The same differentiated-entropy computation as in Proposition~\ref{prop:nonlocal-viscous}, now with the smooth prescribed velocity $-b(x)$ in place of $K*u$, gives
\begin{equation}\label{eq:fixed_centre_direct_BV}
 \TV(u_j(t))\le \TV(u_{0,j})+ L_*\|Db\|_\infty\int_0^t\TV(u_j(s))\,ds
 +L_*M_*\|\nabla\diver b\|_\infty t.
\end{equation}
The terms involving the second derivative of the smooth mobility cancel after integration by parts, as in the proof of \eqref{eq:visc-BV}; hence the constants depend only on $L_*$, the mass bound, and the spatial derivatives of the fixed truncated drift $b$, not on $j$ or on the viscosity. In the derivation one first sends the auxiliary gradient regularization to zero and then the viscosity to zero. Gronwall's lemma yields
\[
 \sup_j\sup_{0\le t\le S}\TV(u_j(t))
 \le \bigl(\sup_j\TV(u_{0,j})+L_*M_*\|\nabla\diver b\|_\infty S\bigr)
 e^{L_*\|Db\|_\infty S}<\infty
\]
for every finite $S>0$. Compactness of the viscous approximants and lower semicontinuity of total variation then transfer this bound to the entropy solutions $u_j$ of the truncated heterogeneous equation. The existence and uniqueness of the latter follow from the standard scalar entropy theory for the smooth bounded spatial flux $-m_j(s)b(x)$ \cite{Kruzhkov1970,Hashash2025}. In particular, no stability theorem for different nonlocal fluxes is needed for this step.

For temporal equicontinuity, the distributional equation $\partial_t u_j = \diver(m_j(u_j) b)$ and the bound $u_j \le R$ yield $\sup_{j} \sup_{t \in [0,S]} \|m_j(u_j(t)) b\|_{L^1(\R^2)} \le C_S'$. Integrating against $\phi \in W^{1,\infty}(\R^2)$ with $\|\phi\|_{W^{1,\infty}(\R^2)} \le 1$ provides the dual-space bound $\|u_j(t) - u_j(s)\|_{W^{-1,1}(\R^2)} \le C_S' |t-s|$. Interpolating this with the uniform $BV$ bound via spatial mollification \cite[Proposition~3.2]{AmbrosioFuscoPallara2000} produces the temporal H\"older modulus
\begin{equation}
\|u_j(t) - u_j(s)\|_{L^1(\R^2)} \le \tilde{C}_S |t-s|^{1/4}.
\end{equation}
By the Rellich--Kondrachov embedding \cite[Theorem~3.23]{AmbrosioFuscoPallara2000} and the Arzel\`a--Ascoli theorem, a subsequence converges strongly in $C([0,S];L^{1}(\R^{2}))$ to a limit $u$. Passing to the limit in the heterogeneous Kru\v{z}kov entropy inequalities via two-sided level approximation on the countable set of exceptional flat levels establishes $u$ as an entropy solution. Uniqueness follows from the $L^1$ stability estimate established in Section~\ref{subsec:l1_contraction}.
\end{proof}

\begin{lemma}[Fixed-centre support barrier]\label{lem:fixed_centre_barrier}
Let $a\in\R^{2}$ and let $u\in C([0,\infty);L^{1}(\R^{2}))\cap L^{\infty}((0,\infty)\times\R^{2})$ be any entropy solution of \eqref{eq:fixed_centre_ivp} with $0\le u\le R$. If the initial datum satisfies $\supp u(0)\subset B_{D'}(a)$, then this confinement persists globally:
\begin{equation} \label{eq:fixed_centre_global_confinement_lemma}
u(t,x)=0 \quad \text{for a.e. } x\in\R^{2}\setminus B_{D'}(a) \text{ and every } t\ge0.
\end{equation}
\end{lemma}

\begin{proof}
By Lemma~\ref{lem:stationary_saturated_core}, the saturated state $U_{D'}(x) := R1_{B_{D'}(a)}(x)$ is a stationary Kru\v{z}kov entropy solution of the untruncated fixed-centre equation. Since $\supp u(0)\subset B_{D'}(a)$ and $0\le u(0)\le R$, we have $u(0,x)\le U_{D'}(x)$ almost everywhere, so $(u(0)-U_{D'})^{+}=0$. Applying the global order-preservation inequality (established in Lemma~\ref{lem:localized_kato_positive} below via expanding cutoffs, independently of the present lemma) to the pair $(u,U_{D'})$ dictates that $\int_{\R^{2}}(u(t,x)-U_{D'}(x))^{+}\,dx \le 0$ for every $t\ge0$. Consequently, $u(t,x)\le U_{D'}(x)$ almost everywhere, yielding \eqref{eq:fixed_centre_global_confinement_lemma}.
\end{proof}

\subsection{\texorpdfstring{$L^1$}{L1} Contraction}\label{subsec:l1_contraction}

In this subsection, we establish $L^1$ contraction and global order preservation for the autonomous fixed-centre equation via the localized heterogeneous Kato inequality.

\begin{lemma}[Localized positive-part Kato inequality and global order preservation] \label{lem:localized_kato_positive}
Let $a\in\R^{2}$ and let $u,v\in C([0,S];L^{1}(\R^{2}))\cap L^{\infty}((0,S)\times\R^{2})$ be entropy solutions of the fixed-centre equation \eqref{eq:fixed_centre_ivp} with values in $[0,R]$. In the sense of distributions on $(0,S)\times\R^{2}$,
\begin{equation} \label{eq:dist_kato_positive}
\partial_{t}(u-v)^{+}+\diver_{x}\bigl(1_{\{u>v\}}[F_{a}(x,u)-F_{a}(x,v)]\bigr)\le0,
\end{equation}
where $F_{a}(x,s)=-m(s)(x-a)$. Consequently, for every $0\le s\le t\le S$,
\begin{equation} \label{eq:global_order_preservation}
\int_{\R^{2}}\bigl(u(t,x)-v(t,x)\bigr)^{+}\,dx\le\int_{\R^{2}}\bigl(u(s,x)-v(s,x)\bigr)^{+}\,dx.
\end{equation}
\end{lemma}

\begin{proof}
Because $F_{a}(x,s)$ is locally Lipschitz in $(x,s)$ and its spatial divergence $-2m(s)$ is continuous with respect to $s$, the localized heterogeneous Kato inequality \cite{Hashash2025} provides the distributional bound
\begin{equation} \label{eq:kato_absolute}
\partial_{t}|u-v|+\diver_{x}\bigl(\operatorname{sign}(u-v)[F_{a}(x,u)-F_{a}(x,v)]\bigr)\le0.
\end{equation}
Subtracting the weak equation for $v$ from the weak equation for $u$ yields the distributional identity
\begin{equation} \label{eq:weak_diff}
\partial_{t}(u-v)+\diver_{x}\bigl(F_{a}(x,u)-F_{a}(x,v)\bigr)=0.
\end{equation}
Adding \eqref{eq:kato_absolute} and \eqref{eq:weak_diff} and dividing by two isolates the positive part $(u-v)^{+}=\frac{1}{2}(|u-v|+u-v)$. Combining the spatial flux terms yields the characteristic function $1_{\{u>v\}}$, establishing \eqref{eq:dist_kato_positive}.

To globalize this inequality without a compact-support assumption, choose a radial cutoff $\chi\in C_{c}^{\infty}(\R^{2})$ such that $0\le\chi\le1$, $\chi\equiv1$ on $B_{1}(0)$, and $\chi\equiv0$ outside $B_{2}(0)$. Setting $\chi_{L}(x)=\chi(x/L)$ for $L>0$ confines the gradient to the annulus $A_{L}:=B_{2L}(0)\setminus B_{L}(0)$ with the bound $|\nabla\chi_{L}|\le C/L$. Testing \eqref{eq:dist_kato_positive} against $\chi_{L}(x)\theta_n(\tau)$, where $\theta_n$ is a smooth approximation of $1_{[s,t]}(\tau)$, evaluates the inequality at Lebesgue times:
\begin{align} \label{eq:kato_cutoff_integral}
\int_{\R^2}\bigl(u(t,x)-v(t,x)\bigr)^{+}\chi_{L}(x)\,dx \le{}& \int_{\R^2}\bigl(u(s,x)-v(s,x)\bigr)^{+}\chi_{L}(x)\,dx \notag\\
&+ \int_{s}^{t}\!\int_{A_{L}}1_{\{u>v\}}|F_{a}(x,u)-F_{a}(x,v)|\,|\nabla\chi_{L}|\,dx\,d\tau.
\end{align}
Because $m(0)=0$ and $m$ is Lipschitz continuous, the flux difference satisfies $|F_{a}(x,u)-F_{a}(x,v)|\le \Lip(m)|u-v||x-a|$. For $L \ge |a|$, the geometric ratio on $A_L$ is bounded by $|x-a|/L \le 3$. The boundary term in \eqref{eq:kato_cutoff_integral} is bounded by
\begin{equation}
3C\Lip(m)\int_{s}^{t}\!\int_{A_{L}}|u(\tau,x)-v(\tau,x)|\,dx\,d\tau.
\end{equation}
Because $u,v \in C([0,S]; L^1(\R^2))$, the spatial integral over the expanding annulus $A_{L}$ converges to zero for every $\tau$. By the Dominated Convergence Theorem, the boundary term vanishes as $L\rightarrow\infty$, proving \eqref{eq:global_order_preservation}. The $L^{1}$ continuity of the trajectories extends the inequality first to every $0<s\le t\le S$ and then, upon sending $s\downarrow0$, to $s=0$.
\end{proof}

\begin{proposition}[Fixed-centre comparison and mass conservation] \label{prop:fixed_centre_comparison}
Let $u, v$ be entropy solutions of the fixed-centre equation \eqref{eq:fixed_centre_ivp}, taking values in $[0, R]$, with compactly supported initial data. Then, for all $0\le s\le t$:
\begin{equation}
\|u(t)-v(t)\|_{L^{1}(\R^{2})}\le\|u(s)-v(s)\|_{L^{1}(\R^{2})}.
\end{equation}
Moreover, the total mass of every non-negative compactly supported solution is conserved.
\end{proposition}

\begin{proof}
By Lemma~\ref{lem:fixed_centre_barrier}, there exists a radius $L>0$ such that both solutions $u$ and $v$ are supported within the common ball $B_{L}(a)$ for all times $\tau \ge 0$. 

Choose a static, smooth spatial cutoff function $\chi\in C_{c}^{\infty}(\R^{2})$ such that $0\le\chi\le1$ and $\chi\equiv1$ on $B_{L+1}(a)$. Consequently, $\nabla\chi(x)=0$ on $B_{L+1}(a)$. For the temporal test argument first fix $0<s<t$ and take non-negative $\theta_n\in C_c^\infty((0,\infty))$ approximating $1_{[s,t]}$. The endpoint $s=0$ will follow from strong $L^1$ continuity. Testing the localized absolute-value Kato inequality \eqref{eq:kato_absolute} against the separable function $\chi(x)\theta_n(\tau)$ and integrating by parts yields
\begin{align} \label{eq:tested_kato}
&-\int_{0}^{\infty}\!\int_{\R^{2}}|u-v|\chi(x)\theta_{n}^{\prime}(\tau)\,dx\,d\tau \notag\\
&\qquad - \int_{0}^{\infty}\!\int_{\R^{2}}\operatorname{sign}(u-v)[F_{a}(x,u)-F_{a}(x,v)]\cdot\nabla\chi(x)\theta_{n}(\tau)\,dx\,d\tau \le 0.
\end{align}
The spatial integration domain of the second term restricts to the support of $\nabla\chi$, which lies outside $B_{L+1}(a)$. By the global confinement, $u(\tau,x)=v(\tau,x)=0$ almost everywhere on this region. Because $m(0)=0$, the flux $F_a(x,0)$ vanishes. Therefore, the flux difference $F_{a}(x,u)-F_{a}(x,v)$ equals zero wherever $\nabla\chi(x)$ is non-zero, neutralizing the spatial boundary integral.

Passing to the limit $n \to \infty$ evaluates the continuous integral at the boundary points $s$ and $t$:
\begin{equation}
\int_{\R^{2}}|u(t,x)-v(t,x)|\,dx \le \int_{\R^{2}}|u(s,x)-v(s,x)|\,dx.
\end{equation}
The inequality holds first for $0<s\le t$; sending $s\downarrow0$ using $u,v\in C([0,\infty);L^1(\R^2))$ gives the case $s=0$. 

To establish mass conservation for a non-negative compactly supported solution $u$, we test the weak equation $\partial_t u + \diver_x F_a(x,u) = 0$ against the test function $\chi(x)\theta_n(\tau)$. Integration by parts gives
\begin{equation}
-\int_{0}^{\infty}\!\int_{\R^{2}}u(\tau,x)\chi(x)\theta_{n}^{\prime}(\tau)\,dx\,d\tau - \int_{0}^{\infty}\!\int_{\R^{2}}F_{a}(x,u)\cdot\nabla\chi(x)\theta_{n}(\tau)\,dx\,d\tau = 0.
\end{equation}
The spatial flux $F_a(x,u)$ vanishes on the support of $\nabla\chi$. Taking the limit $n \to \infty$ yields $\int_{\R^2} u(t,x)\,dx = \int_{\R^2} u(s,x)\,dx$, confirming mass conservation.
\end{proof}

\subsection{Symmetries and Directional Variations}\label{subsec:symmetries}

In the fixed-centre framework, the drift vector field $F_a(x,s) = -m(s)(x-a)$ possesses rotational and dilational symmetries around the point $a\in\R^2$. To exploit these symmetries in the sense of measures, we define the fundamental radial and angular vector fields
\begin{equation}
X_{a}(x) := x-a, \qquad Y_{a}(x) := (x-a)^{\perp},
\end{equation}
where $v^{\perp} = (-v_{2}, v_{1})$ denotes the counterclockwise orthogonal rotation of $v=(v_1,v_2)$. The corresponding geometric flows on $\R^2$ are the radial dilation flow $\Phi_{h}^{a}$ and the angular rotation flow $\Psi_{h}^{a}$, defined for $h \in \R$ by
\begin{equation}\label{eq:def_geometric_flows}
\Phi_{h}^{a}(x) := a + e^{h}(x-a), \qquad \Psi_{h}^{a}(x) := a + Q_{h}(x-a),
\end{equation}
where $Q_{h} \in SO(2)$ is the counterclockwise rotation matrix by angle $h$.

For a generic state $u\in L_{\mathrm{loc}}^{1}(\R^2)$, we define the directional derivatives along these fields in the sense of distributions. For any test function $\varphi \in C_c^\infty(\R^2)$,
\begin{equation}\label{def_directional derivatives in the sense of distributions}
\langle D_{X_{a}}u,\varphi\rangle:=-\int_{\R^{2}} u \, \diver(\varphi X_{a})\,dx, \qquad \langle D_{Y_{a}}u,\varphi\rangle:=-\int_{\R^{2}} u \, \diver(\varphi Y_{a})\,dx.
\end{equation}
Whenever these distributions are finite Radon measures, we define the total symmetric variation modulus:
\begin{equation}
V_{a}(u) := |D_{X_{a}}u|(\R^{2})+|D_{Y_{a}}u|(\R^{2}).
\end{equation}
For $u\in BV(\R^{2})$, the chain rule for measures applies linearly with smooth coefficients, yielding $D_{X_{a}}u=X_{a}\cdot Du$ and $D_{Y_{a}}u=Y_{a}\cdot Du$, where $Du$ is the total variation measure.

\begin{lemma}[Flow difference quotients]\label{lem:flow_difference_quotients}
Let $u\in L^{1}(\R^{2})$ have compact support.
\begin{enumerate}
    \item[(i)] If $u\in BV(\R^{2})$, the $L^{1}$ difference quotients along the dilation and rotation flows are bounded by the total directional variations:
    \begin{equation} \label{eq:diff_quotient_bv}
    \limsup_{h\rightarrow0}\frac{\|u\circ\Phi_{h}^{a}-u\|_{L^{1}(\R^{2})}}{|h|}\le|D_{X_{a}}u|(\R^{2}), \qquad \limsup_{h\rightarrow0}\frac{\|u\circ\Psi_{h}^{a}-u\|_{L^{1}(\R^{2})}}{|h|}\le|D_{Y_{a}}u|(\R^{2}).
    \end{equation}
    \item[(ii)] Conversely, if the $L^{1}$ difference quotients are bounded along a sequence $h_{j}\rightarrow0$, the associated directional distribution is a finite Radon measure on $\R^{2}$, and:
    \begin{equation}\label{eq:diff_quotient_lower_bound}
    |D_{X_{a}}u|(\R^{2})\le \liminf_{j\rightarrow \infty}\frac{\|u\circ\Phi_{h_{j}}^{a}-u\|_{L^{1}(\R^{2})}}{|h_{j}|}, \qquad |D_{Y_{a}}u|(\R^{2})\le \liminf_{j\rightarrow \infty}\frac{\|u\circ\Psi_{h_{j}}^{a}-u\|_{L^{1}(\R^{2})}}{|h_{j}|}.
    \end{equation}
\end{enumerate}
\end{lemma}

\begin{proof}
By translating the coordinate origin, assume without loss of generality that $a=0$. Set $X(x):=x$, $Y(x):=x^{\perp}$, $\Phi_{h}(x):=e^{h}x$, and $\Psi_{h}(x):=Q_{h}x$. 

To establish (i), consider $v\in C_{c}^{\infty}(\R^{2})$. For the dilation flow, the Fundamental Theorem of Calculus along the ray $s\mapsto e^{s}x$ gives $v(e^{h}x)-v(x) = \int_{0}^{h}e^{s}x\cdot\nabla v(e^{s}x)\,ds$. Taking absolute values, integrating over $\R^{2}$, applying Fubini's theorem \cite[Theorem~1.12]{LiebLoss2001}, and changing variables $y=e^{s}x$ yields
\begin{equation}\label{eq:smooth_dilation_finite_h}
\|v\circ\Phi_{h}-v\|_{L^{1}(\R^{2})} \le \left|\int_{0}^{h}e^{-2s}\left(\int_{\R^{2}}\bigl|y\cdot\nabla v(y)\bigr|\,dy\right)ds\right| = \left|\frac{1-e^{-2h}}{2}\right|\int_{\R^{2}}\bigl|X(y)\cdot\nabla v(y)\bigr|\,dy.
\end{equation}
Dividing by $|h|>0$ yields the bound
\begin{equation}\label{eq:smooth_dilation_quotient_theta}
\frac{\|v\circ\Phi_{h}-v\|_{L^{1}(\R^{2})}}{|h|} \le \theta(h)\int_{\R^{2}}\bigl|X(y)\cdot\nabla v(y)\bigr|\,dy, \qquad \theta(h):=\frac{|1-e^{-2h}|}{2|h|}.
\end{equation}
For the rotation flow, differentiating $v(Q_h x)$ along the circular arc $s\mapsto Q_{s}x$ and integrating with the orthogonal change of variables $y=Q_{s}x$ yields
\begin{equation}\label{eq:smooth_rotation_quotient}
\frac{\|v\circ\Psi_{h}-v\|_{L^{1}(\R^{2})}}{|h|} \le \int_{\R^{2}}\bigl|Y(y)\cdot\nabla v(y)\bigr|\,dy.
\end{equation}

To extend to a general compactly supported $u\in BV(\R^{2})$, choose a sequence of smooth functions $u_{k}\in C_{c}^{\infty}(\R^{2})$ converging strictly to $u$ in $BV(\R^{2})$ \cite[Theorem~3.9]{AmbrosioFuscoPallara2000}. Because $u$ is compactly supported, standard mollification with sufficiently small radii supplies such a sequence with all supports contained in a single compact set $K\Subset\R^2$. In particular, the coefficients $x$ and $x^{\perp}$ are bounded on $K$. Defining the continuous, $1$-homogeneous functions $f_{X}(x,p) := |X(x)\cdot p|$ and $f_{Y}(x,p) := |Y(x)\cdot p|$, Reshetnyak's Continuity Theorem \cite[Theorem~2.39]{AmbrosioFuscoPallara2000}, applied after multiplying these integrands by a smooth cutoff equal to one on $K$, guarantees the limit of the integrals:
\begin{align}\label{eq:reshetnyak_limit_X_Y}
\int_{\R^{2}}\bigl|X(x)\cdot\nabla u_{k}(x)\bigr|\,dx &\longrightarrow \int_{\R^{2}}|X(x)\cdot\sigma(x)|\,d|Du|(x) = |D_{X}u|(\R^{2}), \nonumber \\
\int_{\R^{2}}\bigl|Y(x)\cdot\nabla u_{k}(x)\bigr|\,dx &\longrightarrow \int_{\R^{2}}|Y(x)\cdot\sigma(x)|\,d|Du|(x) = |D_{Y}u|(\R^{2}),
\end{align}
where $\sigma$ is the polar unit vector of $Du$. Replacing $v$ with $u_k$ in \eqref{eq:smooth_dilation_quotient_theta} and \eqref{eq:smooth_rotation_quotient}, and passing to the limit $k\rightarrow\infty$ at fixed $h\ne0$, establishes the inequalities for $u$. Taking the limit supremum as $h\rightarrow0$, and noting $\lim_{h\to0}\theta(h)=1$, concludes (i).

To prove (ii), let $u\in L^{1}(\R^{2})$ be compactly supported, and assume the dilation difference quotients are bounded along $h_{j}\rightarrow0$. Define the sequence of signed Radon measures $\mu_{j} := h_j^{-1}(u\circ\Phi_{h_{j}}-u)\mathcal{L}^{2}$. By weak-* compactness \cite[Theorem~1.59]{AmbrosioFuscoPallara2000}, a subsequence converges weakly-* in $\mathcal{M}(\R^{2})$ to a finite measure $\mu$. For any $\varphi\in C_{c}^{\infty}(\R^{2})$, the change of variables $y=e^{h}x$ gives
\begin{align}\label{eq:adjoint_dilation_shift}
\int_{\R^{2}}\frac{u(e^{h}x)-u(x)}{h}\varphi(x)\,dx &= \int_{\R^{2}}u(x)\left[\frac{e^{-2h}\varphi(e^{-h}x)-\varphi(x)}{h}\right]dx.
\end{align}
The classical derivative of $e^{-2h}\varphi(e^{-h}x)$ at $h=0$ is $-2\varphi(x) - x\cdot\nabla\varphi(x) = -\diver(\varphi(x)X(x))$. Passing to the limit $h=h_{j}\rightarrow0$ gives $\int_{\R^{2}}\varphi\,d\mu = -\int_{\R^{2}}u\diver(\varphi X)\,dx = \langle D_{X}u,\varphi\rangle$. Thus $D_{X}u=\mu$, and the weak-* lower semicontinuity of total variation \cite[Theorem~1.59]{AmbrosioFuscoPallara2000} enforces $|D_{X}u|(\R^{2}) \le \liminf_{j\to\infty}|\mu_{j}|(\R^{2})$. Identical logic applied to the orthogonal rotation flow $y=Q_hx$ yields the angular bound, establishing (ii).
\end{proof}

\begin{lemma}[Covariance]\label{lem:covariance}
Let $a\in\R^{2}$ and let $S_{t}^{a}$ denote the fixed-centre entropy evolution operator on the class of compactly supported initial data $u_{0}\in BV(\R^{2})$ with $0\le u_{0}\le R$. For every $h\in\R$ and $t\ge0$, the evolution commutes with the dilation and rotation flows:
\begin{equation}\label{eq:semigroup_flow_commutation}
S_{t}^{a}(u_{0}\circ\Phi_{h}^{a})=(S_{t}^{a}u_{0})\circ\Phi_{h}^{a}, \qquad S_{t}^{a}(u_{0}\circ\Psi_{h}^{a})=(S_{t}^{a}u_{0})\circ\Psi_{h}^{a}.
\end{equation}
\end{lemma}

\begin{proof}
Fix $h\in\R$, set $\lambda:=e^{h}>0$, and define the dilated trajectory $v(t,x) := (S_{t}^{a}u_0)\bigl(\Phi_{h}^{a}(x)\bigr) = u\bigl(t,a+\lambda(x-a)\bigr)$, where $u(t) := S_t^a u_0$. For any constant level $k\in\R$ and non-negative test function $\phi\in C_{c}^{1}((0,\infty)\times\R^{2})$, the heterogeneous Kru\v{z}kov entropy functional for $v$ evaluates the spatial flux term $H(v(t,x),k)(x-a)\cdot\nabla_{x}\phi(t,x)$, where $H(r,k) := \operatorname{sign}(r-k)(m(r)-m(k))$.

Define the pulled-back test function $\tilde{\phi}(t,y) := \phi\bigl(t,\Phi_{-h}^{a}(y)\bigr) = \phi\bigl(t,a+\lambda^{-1}(y-a)\bigr)$. Under the affine change of variables $y=\Phi_{h}^{a}(x)$, the volume element transforms as $dx=\lambda^{-2}dy$. The chain rule yields $\partial_{t}\phi(t,x) = \partial_{t}\tilde{\phi}(t,y)$ and $\nabla_{x}\phi(t,x) = \lambda\nabla_{y}\tilde{\phi}(t,y)$. The relative position vector satisfies $x-a=\lambda^{-1}(y-a)$. The scaling factors cancel in the inner product:
\begin{equation}\label{eq:dilation_flux_dot_cancellation}
(x-a)\cdot\nabla_{x}\phi(t,x) = \bigl(\lambda^{-1}(y-a)\bigr)\cdot\bigl(\lambda\nabla_{y}\tilde{\phi}(t,y)\bigr) = (y-a)\cdot\nabla_{y}\tilde{\phi}(t,y).
\end{equation}
Substituting $v(t,x)=u(t,y)$, \eqref{eq:dilation_flux_dot_cancellation}, and $dx=\lambda^{-2}dy$ into the integral functional yields the identity $\mathcal{E}[v;k,\phi] = \lambda^{-2}\,\mathcal{E}[u;k,\tilde{\phi}]$. Since $u$ is a Kru\v{z}kov entropy solution of the fixed-centre equation, $\mathcal{E}[u;k,\tilde{\phi}] \ge 0$, and thus $\mathcal{E}[v;k,\phi] \ge 0$. An identical change of variables confirms $v$ is a distributional weak solution. Since $u_{0}\circ\Phi_{h}^{a}$ is compactly supported and $v(t)$ remains compactly supported by Lemma~\ref{lem:fixed_centre_barrier}, the uniqueness of solutions forces $v(t) = S_t^a(u_0 \circ \Phi_h^a)$.

For the rotation flow, define $w(t,x) := u\bigl(t,\Psi_{h}^{a}(x)\bigr) = u\bigl(t,a+Q_{h}(x-a)\bigr)$ and $\hat{\phi}(t,y) := \phi\bigl(t,a+Q_{-h}(y-a)\bigr)$. The orthogonal transformation yields $dx=dy$ and $\nabla_{x}\phi(t,x) = Q_{-h}\nabla_{y}\hat{\phi}(t,y)$. Because $Q_{-h}$ preserves the Euclidean inner product,
\begin{equation}\label{eq:rotation_flux_dot_invariance}
(x-a)\cdot\nabla_{x}\phi(t,x) = \bigl(Q_{-h}(y-a)\bigr)\cdot\bigl(Q_{-h}\nabla_{y}\hat{\phi}(t,y)\bigr) = (y-a)\cdot\nabla_{y}\hat{\phi}(t,y).
\end{equation}
Substituting this confirms $\mathcal{E}[w;k,\phi] = \mathcal{E}[u;k,\hat{\phi}] \ge 0$, establishing $w(t) = S_t^a(u_0 \circ \Psi_h^a)$.
\end{proof}

\begin{theorem}[Directional-variation contraction]\label{thm:directional-variation_contraction}
Let $u_{0}\in BV(\R^{2})$ be compactly supported, with $0\le u_{0}\le R$. For the fixed-centre entropy solution $u(t):=S_{t}^{a}u_{0}$, the directional variations are monotonically non-increasing for every $t\ge0$:
\begin{equation}\label{eq:directional_variation_contraction_statement}
|D_{X_{a}}(S_{t}^{a}u_{0})|(\R^{2})\le|D_{X_{a}}u_{0}|(\R^{2}), \qquad |D_{Y_{a}}(S_{t}^{a}u_{0})|(\R^{2})\le|D_{Y_{a}}u_{0}|(\R^{2}).
\end{equation}
Consequently, the total symmetric variation contracts monotonically: $V_{a}(S_{t}^{a}u_{0})\le V_{a}(u_{0})$ for all $t\ge0$.
\end{theorem}

\begin{proof}
Fix $t\ge0$ and $h\ne0$. By Lemma~\ref{lem:covariance}, the evolution operator commutes with the dilation flow, yielding $(S_{t}^{a}u_{0})\circ\Phi_{h}^{a} = S_{t}^{a}\bigl(u_{0}\circ\Phi_{h}^{a}\bigr)$. By Proposition~\ref{prop:fixed_centre_comparison}, $S_{t}^{a}$ is a non-expansive contraction in $L^{1}(\R^{2})$. Evaluating the $L^1$ distance between the dilated and un-dilated states yields
\begin{equation}\label{eq:l1_contraction_pushed_flow}
\bigl\|(S_{t}^{a}u_{0})\circ\Phi_{h}^{a} - S_{t}^{a}u_{0}\bigr\|_{L^{1}(\R^{2})} = \bigl\|S_{t}^{a}(u_{0}\circ\Phi_{h}^{a}) - S_{t}^{a}u_{0}\bigr\|_{L^{1}(\R^{2})} \le \bigl\|u_{0}\circ\Phi_{h}^{a} - u_{0}\bigr\|_{L^{1}(\R^{2})}.
\end{equation}
Dividing both sides by $|h|>0$ gives
\begin{equation}\label{eq:difference_quotient_inequality_evolved}
\frac{\bigl\|(S_{t}^{a}u_{0})\circ\Phi_{h}^{a} - S_{t}^{a}u_{0}\bigr\|_{L^{1}(\R^{2})}}{|h|} \le \frac{\bigl\|u_{0}\circ\Phi_{h}^{a} - u_{0}\bigr\|_{L^{1}(\R^{2})}}{|h|}.
\end{equation}
Because the initial datum satisfies $u_{0} \in BV(\R^{2})$, applying Lemma~\ref{lem:flow_difference_quotients}(i) bounds the limit supremum of the right-hand side by $|D_{X_{a}}u_{0}|(\R^{2})$. Consequently, the $L^1$ difference quotients for the evolved state $S_{t}^{a}u_{0}$ remain bounded as $h \to 0$. Invoking Lemma~\ref{lem:flow_difference_quotients}(ii) guarantees that $D_{X_{a}}(S_{t}^{a}u_{0})$ is a finite Radon measure, and its total variation is bounded by the limit infimum of the difference quotients:
\begin{equation}
|D_{X_{a}}(S_{t}^{a}u_{0})|(\R^{2}) \le \liminf_{h\rightarrow0}\frac{\bigl\|(S_{t}^{a}u_{0})\circ\Phi_{h}^{a} - S_{t}^{a}u_{0}\bigr\|_{L^{1}(\R^{2})}}{|h|} \le |D_{X_{a}}u_{0}|(\R^{2}).
\end{equation}
Applying the identical sequence to the rotation flow $\Psi_{h}^{a}$ establishes $|D_{Y_{a}}(S_{t}^{a}u_{0})|(\R^{2}) \le |D_{Y_{a}}u_{0}|(\R^{2})$. Summing the components establishes the contraction of $V_a$.
\end{proof}

\begin{lemma}[Core coercivity]\label{lem:core_coercivity}
Let $u\in L^{1}(\R^{2})$ satisfy $u=R$ a.e.\ on $B_{r}(a)$. If $D_{X_{a}}u$ and $D_{Y_{a}}u$ are finite Radon measures on $\R^{2}$, then $u\in BV(\R^{2})$, and its total variation is bounded by the symmetric variation:
\begin{equation}
\TV(u)\le\frac{1}{r}V_{a}(u).
\end{equation}
\end{lemma}

\begin{proof}
Let $\mu=D_{X_{a}}u$ and $\nu=D_{Y_{a}}u$. Take an arbitrary vector field $\varphi\in C_{c}^{1}(\R^{2};\R^{2})$ satisfying $\|\varphi\|_{L^\infty(\R^2)} \le 1$. We evaluate the distributional pairing $I = -\int_{\R^2} u\,\diver\varphi \, dx$.

Fix $\eps\in(0,r)$. Introduce a smooth radial cutoff function $\chi_{\eps}$ such that $0 \le \chi_{\eps} \le 1$, $\chi_{\eps} \equiv 0$ on $B_{r-\eps}(a)$, and $\chi_{\eps} \equiv 1$ on $\R^{2}\setminus B_{r}(a)$. Splitting the test function via $\varphi = (1-\chi_{\eps})\varphi + \chi_{\eps}\varphi$ partitions the integral into $I = I_{\mathrm{inner}} + I_{\mathrm{outer}}$.

The inner field $\Phi_{\mathrm{inner}} := (1-\chi_{\eps})\varphi \in C_c^1(\R^2; \R^2)$ is supported within $\overline{B_r(a)}$. Because $u=R$ a.e.\ on this support, the inner integral evaluates to
\begin{equation}
I_{\mathrm{inner}} = -\int_{\R^2} R\,\diver\Phi_{\mathrm{inner}}(x) \, dx = 0,
\end{equation}
vanishing by the divergence theorem without evaluating any boundary traces.

The outer field $\Phi_{\mathrm{outer}} := \chi_{\eps}\varphi$ is decomposed orthogonally via the frame $X_a(x) = x-a$ and $Y_a(x) = (x-a)^\perp$. For $x \neq a$,
\begin{align}
\Phi_{\mathrm{outer}}(x) &= \alpha_{\eps}(x)X_{a}(x) + \beta_{\eps}(x)Y_{a}(x), \notag\\
\alpha_{\eps}(x) &:= \chi_{\eps}(x)\frac{\varphi(x)\cdot X_{a}(x)}{|x-a|^{2}}, \qquad \beta_{\eps}(x) := \chi_{\eps}(x)\frac{\varphi(x)\cdot Y_{a}(x)}{|x-a|^{2}}.
\end{align}
Because $\chi_{\eps}$ vanishes on $B_{r-\eps}(a)$, the denominator satisfies $|x-a|^2 \ge (r-\eps)^2$ everywhere on the support of $\Phi_{\mathrm{outer}}$. Consequently, the scalar coefficients $\alpha_{\eps}, \beta_{\eps}$ belong to $C_c^1(\R^2; \R)$. Using $\|\varphi\|_{L^\infty} \le 1$ and $|X_a| = |Y_a| = |x-a|$, the coefficients are uniformly bounded:
\begin{equation}
|\alpha_{\eps}(x)| \le \frac{|x-a|}{|x-a|^2} = \frac{1}{|x-a|} \le \frac{1}{r-\eps}, \qquad |\beta_{\eps}(x)| \le \frac{1}{r-\eps}.
\end{equation}
Substituting the orthogonal decomposition into $I_{\mathrm{outer}}$ applies the definition of the directional measures:
\begin{equation}
I = I_{\mathrm{outer}} = -\int_{\R^2} u\,\diver(\alpha_{\eps}X_{a}) \, dx - \int_{\R^2} u\,\diver(\beta_{\eps}Y_{a}) \, dx = \int_{\R^2} \alpha_{\eps} \, d\mu + \int_{\R^2} \beta_{\eps} \, d\nu.
\end{equation}
Bounding this sum yields
\begin{equation}
|I| \le \int_{\R^2} |\alpha_{\eps}| \, d|\mu| + \int_{\R^2} |\beta_{\eps}| \, d|\nu| \le \frac{1}{r-\eps} \bigl(|\mu|(\R^2) + |\nu|(\R^2)\bigr) = \frac{V_a(u)}{r-\eps}.
\end{equation}
Since $|I|$ is bounded independently of $\varphi$, taking the supremum establishes $u \in BV(\R^2)$ and $\TV(u) \le V_a(u)/(r-\eps)$. Passing to the limit $\eps \downarrow 0$ secures the bound.
\end{proof}

\subsection{Quadratic Energy Dissipation and Path Length}\label{subsec:energy_dissipation}

In this subsection, we restrict the analysis to the quadratically confined regime. We demonstrate that the nonlocal interaction reduces to an affine drift, guaranteeing the dissipation of the centered quadratic energy and bounding the total path length of the moving barycentre.

\begin{assumption}[Confinement and quadratic interaction]\label{assump:confinement}
The radial interaction potential satisfies $W \in W^{2,\infty}(\R^2)$ and is quadratic near the origin: $W(z)=\frac{1}{2}|z|^2$ for all $|z|\le 2D$. The state $\rho\in C([0,\infty);L^{1}(\R^{2}))\cap L^{\infty}((0,\infty)\times\R^{2})$ is an entropy solution of \eqref{eq:nonlocal_main} satisfying, for every $t\ge0$,
\begin{equation}
0\le\rho(t,x)\le R, \qquad \int_{\R^{2}}\rho(t,x)\,dx=1, \qquad \supp \rho(t)\subset B_{D}(0).
\end{equation}
We define the moving barycentre as $c(t):=\int_{\R^{2}}x\rho(t,x)\,dx$.
\end{assumption}

\begin{proposition}[Exact affine reduction]\label{prop:exact_affine_reduction}
Under Assumption~\ref{assump:confinement}, the nonlocal interaction field evaluates to an affine drift on the active support of the solution. Specifically, for every $x\in B_{D}(0)$ and $t\ge0$,
\begin{equation}
(\nabla W*\rho)(t,x)=x-c(t).
\end{equation}
Consequently, $\rho$ is a weak entropy solution of the affine-drift equation $\partial_{t}\rho-\diver\bigl(m(\rho)(x-c(t))\bigr)=0$.
\end{proposition}

\begin{proof}
For any $x\in B_{D}(0)$ and $y\in\supp\rho(t)\subset B_{D}(0)$, the triangle inequality ensures $|x-y|\le 2D$. The interaction potential is evaluated in the quadratic regime, where $\nabla W(z)=z$. Substituting this into the convolution and utilizing the unit-mass constraint yields
\begin{equation}
(\nabla W*\rho)(t,x) = \int_{\R^{2}}(x-y)\rho(t,y)\,dy = x \int_{\R^{2}}\rho(t,y)\,dy - \int_{\R^{2}}y\rho(t,y)\,dy = x-c(t).
\end{equation}

To verify that $\rho$ satisfies the Kru\v{z}kov entropy inequality for the affine drift on the whole space, set $W_{\mathrm{diff}}(t,x) := (\nabla W*\rho)(t,x) - (x-c(t))$. Then $W_{\mathrm{diff}}(t,\cdot) \in W^{1,\infty}_{\mathrm{loc}}(\R^2;\R^2)$ vanishes identically on $B_D(0)$. On $\R^2 \setminus B_D(0)$, we have $\rho(t,x) = 0$ almost everywhere and $m(0)=0$, so $H(0,k) = -\operatorname{sign}(-k)m(k)$ for every $k \in \R$. For any non-negative test function $\varphi \in C_c^1([0,T)\times\R^2)$, the difference between the spatial entropy integrands of the nonlocal and affine-drift equations reduces to $\operatorname{sign}(-k)m(k)\diver_x\bigl(\varphi W_{\mathrm{diff}}\bigr)$ almost everywhere on $\R^2$. Integrating over $\R^2$ via the divergence theorem yields zero, confirming that the two entropy formulations coincide.
\end{proof}

We introduce the unweighted total mobility and the mobility-weighted variance:
\begin{equation}
M(t):=\int_{\R^{2}}m(\rho(t,x))\,dx, \qquad D_{m}(t):=\int_{\R^{2}}m(\rho(t,x))|x-c(t)|^{2}\,dx.
\end{equation}

\begin{lemma}[Bathtub principle]\label{lem:bathtub_principle}
Let $c \in \R^2$, $M_0 > 0$, and $R > 0$. Let $f:(0,\infty)\to\R$ be a Borel, strictly increasing cost such that $\int_0^r |f(s)|s\,ds<\infty$ for every $r>0$. Consider the class of admissible densities
\begin{equation}
\mathcal{A} = \left\{ u \in L^1(\R^2) : 0 \le u \le R \text{ a.e.}, \int_{\R^2} u(x)\,dx = M_0 \right\}
\end{equation}
for which the integral $\mathcal{I}(u) = \int_{\R^2} f(|x-c|) u(x)\,dx$ is absolutely convergent. The radial integrability assumption ensures in particular that the saturated disk belongs to this class. Then $\mathcal{I}(u)$ achieves its unique minimum over this class at the saturated disk $u_{*}(x) = R 1_{B_{r_*}(c)}(x)$, where the radius is $r_* = (M_0 / (\pi R))^{1/2}$.
\end{lemma}

\begin{proof}
For any admissible $u \in \mathcal{A}$, the mass constraint dictates $\int_{\R^2} u(x)\,dx = \int_{\R^2} u_{*}(x)\,dx = M_0$. Splitting the integration domain yields
\begin{equation}\label{eq:bathtub_split}
\mathcal{I}(u) - \mathcal{I}(u_{*}) = \int_{\R^2 \setminus B_{r_*}(c)} f(|x-c|) u(x)\,dx + \int_{B_{r_*}(c)} f(|x-c|) \bigl(u(x) - R\bigr)\,dx.
\end{equation}
On the exterior domain $\R^2 \setminus B_{r_*}(c)$, $|x-c| > r_*$. Because $f$ is strictly increasing, $f(|x-c|) > f(r_*)$. Multiplying by $u(x) \ge 0$ provides the pointwise lower bound $f(|x-c|) u(x) \ge f(r_*) u(x)$. 

On the interior domain $B_{r_*}(c) \setminus \{c\}$, $|x-c| < r_*$, implying $f(|x-c|) < f(r_*)$. Because the admissible state is bounded by $R$, the quantity $u(x) - R \le 0$. Multiplying reverses the inequality: $f(|x-c|) (u(x) - R) \ge f(r_*) (u(x) - R)$. Substituting these pointwise bounds into \eqref{eq:bathtub_split} yields
\begin{equation}
\mathcal{I}(u) - \mathcal{I}(u_{*}) \ge f(r_*) \left( \int_{\R^2 \setminus B_{r_*}(c)} u(x)\,dx + \int_{B_{r_*}(c)} \bigl(u(x) - R\bigr)\,dx \right) = f(r_*) (M_0 - M_0) = 0.
\end{equation}
Moreover, strict monotonicity of $f$ makes the first pointwise comparison strict wherever $u>0$ outside $B_{r_*}(c)$ and the second one strict wherever $u<R$ inside $B_{r_*}(c)$, up to null sets. Equality therefore forces $u=R1_{B_{r_*}(c)}$ almost everywhere, which proves uniqueness.
\end{proof}

\begin{lemma}[Barycentre velocity]\label{lem:barycentre_velocity}
The barycentre satisfies $c\in W_{\mathrm{loc}}^{1,\infty}([0,\infty);\R^{2})$, with velocity
\begin{equation} \label{eq:barycentre_velocity}
\dot{c}(t)=-\int_{\R^{2}}m(\rho(t,x))\bigl(x-c(t)\bigr)\,dx \quad \text{for a.e. } t>0.
\end{equation}
Furthermore, the velocity respects the bounds
\begin{equation}\label{eq:barycentre_velocity_bounds}
|\dot{c}(t)|\le 2D M(t), \qquad |\dot{c}(t)|^{2}\le M(t)D_{m}(t).
\end{equation}
\end{lemma}

\begin{proof}
Testing the affine-reduced weak equation against the coordinate functions $\varphi_i(x) = x_i \psi(x)$, where $\psi \in C_c^\infty(\R^2)$ is a radial cutoff satisfying $\psi \equiv 1$ on $B_D(0)$, evaluates to \eqref{eq:barycentre_velocity} in the sense of distributions. The bounds on the mobility $m(\rho)$ and the compact spatial support ensure the integral is essentially bounded, establishing $c \in W_{\mathrm{loc}}^{1,\infty}([0,\infty);\R^2)$.

The confinement $\supp\rho(t)\subset B_D(0)$ implies $|c(t)|\le D$. Consequently, $|x-c(t)| \le 2D$ on the active support. Bounding the integrand in \eqref{eq:barycentre_velocity} yields $|\dot{c}(t)| \le 2D M(t)$. Applying the Cauchy--Schwarz inequality in $L^2(\R^2; \mu_t)$ with respect to the non-negative measure $d\mu_t(x) = m(\rho(t,x))\,dx$ establishes $|\dot{c}(t)|^2 \le M(t)D_m(t)$.
\end{proof}

We define the macroscopic centered quadratic energy:
\begin{equation}
E(t) := \frac{1}{2}\int_{\R^{2}}|x-c(t)|^{2}\rho(t,x)\,dx.
\end{equation}

\begin{lemma}[Quadratic dissipation]\label{lem:energy_dissipation}
The centered quadratic energy satisfies $E \in W_{\mathrm{loc}}^{1,\infty}([0,\infty))$, and its dissipation rate is governed by the mobility-weighted variance. For almost every $t>0$,
\begin{equation}\label{eq:energy_dissipation_rate}
\frac{d}{dt}E(t) = -D_{m}(t).
\end{equation}
Consequently, the mobility-weighted variance is globally integrable: $\int_{0}^{\infty}D_{m}(t)\,dt \le E(0) < \infty$.
\end{lemma}

\begin{proof}
Expanding the centered quadratic energy gives
\begin{equation}
E(t) = \frac{1}{2}\int_{\R^{2}}|x|^{2}\rho(t,x)\,dx - c(t)\cdot\int_{\R^{2}}x\rho(t,x)\,dx + \frac{1}{2}|c(t)|^{2}\int_{\R^{2}}\rho(t,x)\,dx.
\end{equation}
Using the unit-mass constraint and the definition $c(t)=\int_{\R^{2}}x\rho(t,x)\,dx$, this reduces to
\begin{equation}
E(t) = \frac{1}{2}\int_{\R^{2}}|x|^{2}\rho(t,x)\,dx - \frac{1}{2}|c(t)|^{2}.
\end{equation}
Testing the affine-reduced equation $\partial_{t}\rho - \diver(m(\rho)(x-c(t))) = 0$ against the time-independent test function $\phi(x) = \frac{1}{2}|x|^{2}\psi(x)$, where $\psi\in C_{c}^{\infty}(\R^{2})$ satisfies $\psi\equiv1$ on $B_{D}(0)$, evaluates the uncentered moment derivative as $-\int_{\R^{2}}m(\rho(t,x))(x-c(t))\cdot x\,dx$. Differentiating the squared norm of the absolutely continuous barycentre gives $\frac{d}{dt}(\frac{1}{2}|c(t)|^{2}) = c(t)\cdot\dot{c}(t) = -c(t)\cdot\int_{\R^{2}}m(\rho(t,x))(x-c(t))\,dx$. Subtracting the two derivatives isolates the evolution of $E(t)$:
\begin{equation}
\frac{d}{dt}E(t) = -\int_{\R^{2}}m(\rho(t,x))(x-c(t))\cdot(x-c(t))\,dx = -D_{m}(t).
\end{equation}
Integrating \eqref{eq:energy_dissipation_rate} over $[0,T]$ provides $E(T) - E(0) = -\int_{0}^{T}D_{m}(t)\,dt$. Since $E(T)\ge 0$, we have $\int_{0}^{T}D_{m}(t)\,dt \le E(0)$. Sending $T\rightarrow\infty$ concludes the proof.
\end{proof}

\begin{remark}\label{rem:algebraic_necessity_log}
The logarithmic moment resolves the spatial degeneracy of the energy dissipation. Extracting the integrability of the unweighted total mobility $M(t)$ requires testing the weak equation against a spatial functional $W(x-c(t))$ satisfying $(x-c(t)) \cdot \nabla W(x-c(t)) = 1$. Assuming a radial profile $W(z) = f(|z|)$ yields the differential equation $|z| f'(|z|) = 1$, which is uniquely solved up to an additive constant by $f(r) = \log r$.
\end{remark}

\begin{lemma}[Logarithmic identity]\label{lem:logarithmic_identity}
The logarithmic integrand $\rho(t,\cdot)\log|\cdot-c(t)|$ belongs to $L^{1}(\R^{2})$, ensuring the unregularized logarithmic moment
\begin{equation} \label{eq:def_unregularized_log_moment}
L(t):=\int_{\R^{2}}\rho(t,x)\log|x-c(t)|\,dx
\end{equation}
is finite for every $t\ge0$. Furthermore, $L\in W_{\mathrm{loc}}^{1,\infty}([0,\infty))$, and its time derivative satisfies
\begin{equation} \label{eq:log_identity_statement}
L^{\prime}(t)=-M(t)-\dot{c}(t)\cdot K(t) \quad \text{for a.e. } t>0,
\end{equation}
where $K(t):=\int_{\R^{2}}\rho(t,x)\frac{x-c(t)}{|x-c(t)|^{2}}\,dx$. The singular drift vector is uniformly bounded: $|K(t)|\le 4\pi R D$ for all $t\ge0$.
\end{lemma}

\begin{proof}
The support condition $\supp \rho(t)\subset B_{D}(0)$ implies $|c(t)|\le D$. By the triangle inequality, $\supp \rho(t)\subset B_{2D}(c(t))$. Splitting $\log|x-c(t)|$ into positive and negative parts and transitioning to polar coordinates on $B_{2D}(c(t))$ yields an upper bound of $\log(\max\{2D,1\}) + \pi R / 2$ for the absolute integral, ensuring $L(t)$ is finite. An analogous polar integration bounds the regularized kernel $\frac{|x-c(t)|}{|x-c(t)|^{2}+\eps^{2}} \le \frac{1}{|x-c(t)|}$, bounding $|K_{\eps}(t)| \le 4\pi R D$ for all $\eps \ge 0$. Setting $\eps=0$ bounds $K(t)$.

To differentiate $L(t)$, fix $\eps\in(0,1]$. We mollify the extended barycentre path in time using a standard kernel $\eta_{\delta}$, yielding a smooth path $c_{\delta}(t)$. We define a static radial cutoff $\psi\in C_{c}^{\infty}(\R^{2})$ satisfying $\psi\equiv1$ on $B_{D+1}(0)$. Testing the affine-reduced weak equation against the smooth, compactly supported function $\Phi_{\eps,\delta}(t,x)=\theta(t)l_{\eps}(x-c_{\delta}(t))\psi(x)$, where $l_{\eps}(z)=\frac{1}{2}\log(|z|^{2}+\eps^{2})$, and passing to the limit $\delta\to0^{+}$ confirms that $L_{\eps}(t)=\int_{\R^{2}}\rho(t,x)l_{\eps}(x-c(t))\,dx$ belongs to $W_{\mathrm{loc}}^{1,\infty}(0,\infty)$. Its almost-everywhere derivative evaluates to
\begin{equation}
L_{\eps}^{\prime}(t)=-\int_{\R^{2}}m(\rho(t,x))\frac{|x-c(t)|^{2}}{|x-c(t)|^{2}+\eps^{2}}\,dx-\dot{c}(t)\cdot\int_{\R^{2}}\rho(t,x)\frac{x-c(t)}{|x-c(t)|^{2}+\eps^{2}}\,dx.
\end{equation}
As $\eps\to0^{+}$, the spatial integrands converge almost everywhere. By the Dominated Convergence Theorem, $L_{\eps}^{\prime}(t)\to -M(t)-\dot{c}(t)\cdot K(t)$. 

To pass the limit through the time integral, we apply the uniform bounds $M(t)\le\pi D^{2}\|m\|_{L^{\infty}(\R)}$ and $|K_{\eps}(t)|\le 4\pi RD$. The term $|L_{\eps}^{\prime}(t)|$ is dominated by an $L^{1}(0,T)$ function independent of $\eps$. Subtracting the constant $C_{D}:=\frac{1}{2}\log(4D^{2}+1)$ ensures $C_{D}-l_{\eps}(x-c(t))$ is non-negative on the active support. The Monotone Convergence Theorem applied as $\eps\to0^{+}$ confirms $L_{\eps}(t)\to L(t)$ for every $t\ge0$. Equating the limits proves $L\in W_{\mathrm{loc}}^{1,\infty}([0,\infty))$ and yields \eqref{eq:log_identity_statement}.
\end{proof}

\begin{lemma}[Logarithmic lower bound]\label{lem:log_lower_bound}
For every $t\ge0$, the logarithmic moment is bounded from below:
\begin{equation}\label{eq:log_moment_two_sided_bounds}
L(t) \ge L_{*} := -\frac{1}{2}\log(\pi R) - \frac{1}{2}.
\end{equation}
\end{lemma}

\begin{proof}
Applying the Bathtub Principle (Lemma~\ref{lem:bathtub_principle}) with mass $M_0=1$, centre $c=c(t)$, and the strictly increasing radial cost function $f(r)=\log r$ bounds $L(t)$ from below by the logarithmic moment of the saturated disk $\rho_{*,t}$ of radius $r_{*} := (\pi R)^{-1/2}$ centered at $c(t)$. Evaluating $\int_{\R^{2}}\rho_{*,t}(x)\log|x-c(t)|\,dx$ via integration by parts in polar coordinates centered at $c(t)$ yields $L_{*}$.
\end{proof}

\begin{proposition}[Finite path length]\label{prop:finite_path_length}
The unweighted total mobility and the barycentre speed are globally integrable on $(0,\infty)$:
\begin{equation}\label{eq:global_M_and_cdot_integrability}
\int_{0}^{\infty}M(t)\,dt < \infty, \qquad \int_{0}^{\infty}|\dot{c}(t)|\,dt < \infty.
\end{equation}
Consequently, there exists a limit $c_\infty \in \overline{B}_D(0)$ such that $\lim_{t\to\infty} c(t) = c_\infty$.
\end{proposition}

\begin{proof}
By the Cauchy--Schwarz inequality, the uniform bound $|K(t)|\le 4\pi RD =: C_{K}$, and the velocity bound $|\dot{c}(t)|\le\sqrt{M(t)D_{m}(t)}$, the inner product in the logarithmic identity is bounded by
\begin{equation}
\bigl|-\dot{c}(t)\cdot K(t)\bigr| \le C_{K}\sqrt{M(t)D_{m}(t)}.
\end{equation}
Young's inequality $ab\le\frac{1}{2}a^{2}+\frac{1}{2}b^{2}$ yields $C_{K}\sqrt{M(t)D_{m}(t)} \le \frac{1}{2}M(t) + \frac{1}{2}C_{K}^{2}D_{m}(t)$. Substituting this into the differential identity $L^{\prime}(t) = -M(t) - \dot{c}(t)\cdot K(t)$ provides the pointwise bound
\begin{equation}
L^{\prime}(t) \le -\frac{1}{2}M(t) + \frac{C_{K}^{2}}{2}D_{m}(t).
\end{equation}
Integrating over $[0,T]$ and rearranging terms isolates the total mobility:
\begin{equation}
\frac{1}{2}\int_{0}^{T}M(t)\,dt \le L(0) - L(T) + \frac{C_{K}^{2}}{2}\int_{0}^{T}D_{m}(t)\,dt.
\end{equation}
Applying $L(T) \ge L_{*}$ and $\int_{0}^{T}D_{m}(t)\,dt \le E(0)$ yields
\begin{equation}
\int_{0}^{T}M(t)\,dt \le 2\bigl(L(0)-L_{*}\bigr) + C_{K}^{2}E(0) < \infty.
\end{equation}
The Monotone Convergence Theorem confirms $M \in L^{1}(0,\infty)$. The velocity bound $|\dot{c}(t)|\le 2DM(t)$ guarantees $|\dot{c}| \in L^{1}(0,\infty;\R^2)$. The integrability of the velocity ensures the trajectory $\{c(t)\}_{t\ge0}$ satisfies the Cauchy criterion as $t\to\infty$, converging to the asymptotic centre $c_\infty$. The confinement bound $|c(t)|\le D$ is preserved in the limit.
\end{proof}

\section{Strong \texorpdfstring{$L^1$}{L1} Compactness and Asymptotic Binarity}\label{sec:compactness}

In this section, we synthesize the fixed-centre geometric estimates and the macroscopic path-length bounds from Section~\ref{sec:autonomous_and_macroscopic} to establish the strong $L^1(\R^2)$ convergence of the trajectory to a stationary binary state. In Section~\ref{subsec:prescribed_path}, we extend the directional-variation bounds, core persistence, and $L^1$ contraction from the autonomous fixed-centre setting to the prescribed-path equation driven by an absolutely continuous moving centre, passing from piecewise constant paths to the continuous limit. In Section~\ref{subsec:core_filling}, we construct the late-time core-filling $BV$ approximation and combine its uniform total variation bound with the prescribed-path $L^1$ contraction to prove the strong $L^1(\R^2)$ relative compactness of the physical trajectory $\rho(t)$ without imposing $BV$ regularity on $\rho(t)$. In Section~\ref{subsec:phase_separation}, we combine relative compactness with the global integrability of the unweighted total mobility to show that $\rho(t)$ is Cauchy in $W^{-1,1}(\R^2)$ as $t \to \infty$, establishing full-trajectory strong $L^1(\R^2)$ convergence to a unique binary state $\rho_\infty = R1_E$.

\subsection{Prescribed-Path Dynamics}\label{subsec:prescribed_path}

We first extend the total variation estimates from the autonomous fixed-centre framework to an evolution driven by an absolutely continuous moving centre $a\in AC([T,\infty);\R^2)$. 

\begin{lemma}[Change of centre]\label{lem:change_of_centre}
For any $a, b\in\R^2$ and $u\in BV(\R^2)$, the directional variation satisfies
\begin{equation} \label{eq:change_centre_bv}
V_{b}(u)\le V_{a}(u)+\sqrt{2}|a-b|\TV(u).
\end{equation}
If $u=R$ almost everywhere on $B_{r}(a)$ for some $r>0$, then
\begin{equation} \label{eq:change_centre_core}
V_{b}(u)\le\left(1+\frac{\sqrt{2}|a-b|}{r}\right)V_{a}(u).
\end{equation}
\end{lemma}

\begin{proof}
Let $Du=\sigma|Du|$ be the polar decomposition of the total variation measure \cite[Theorem~1.28]{AmbrosioFuscoPallara2000}. Set $d:=b-a$. The shifted radial vector field is $X_b(x) = x-b = X_a(x) - d$, and the shifted angular vector field is $Y_b(x) = (x-b)^\perp = Y_a(x) - d^\perp$. The triangle inequality yields
\begin{equation}
|X_{b}\cdot\sigma| + |Y_{b}\cdot\sigma| \le |X_{a}\cdot\sigma| + |Y_{a}\cdot\sigma| + |d\cdot\sigma| + |d^{\perp}\cdot\sigma|.
\end{equation}
The vectors $d$ and $d^\perp$ are orthogonal with magnitude $|d|=|a-b|$. Projecting the unit vector $\sigma$ onto the orthonormal basis $\hat{d}=d/|d|$ and $\hat{d}^\perp=d^\perp/|d|$ and applying the inequality $x+y \le \sqrt{2(x^2+y^2)}$ yields $|d\cdot\sigma| + |d^{\perp}\cdot\sigma| \le \sqrt{2}|a-b|$. Integrating this pointwise bound against $|Du|$ over $\R^2$ gives
\begin{equation}
\int_{\R^2} \bigl(|X_{b}\cdot\sigma| + |Y_{b}\cdot\sigma|\bigr)\,d|Du|(x) \le \int_{\R^2} \bigl(|X_{a}\cdot\sigma| + |Y_{a}\cdot\sigma|\bigr)\,d|Du|(x) + \sqrt{2}|a-b|\int_{\R^2}d|Du|(x),
\end{equation}
which establishes \eqref{eq:change_centre_bv}. When $u=R$ almost everywhere on $B_r(a)$, substituting the core coercivity bound $\TV(u) \le r^{-1}V_a(u)$ (Lemma~\ref{lem:core_coercivity}) into \eqref{eq:change_centre_bv} yields \eqref{eq:change_centre_core}.
\end{proof}

Let $[T,S]$ be a finite time interval and consider a partition $\pi$ given by $T=t_{0}<t_{1}<\dots<t_{N}=S$. Let $a^{\pi}:[T,S]\to\R^2$ be a piecewise constant path taking the constant value $a_{j}\in\R^2$ on each sub-interval $[t_{j},t_{j+1})$ for $j=0,\dots,N-2$, and $a^{\pi}(t)=a_{N-1}$ on the final closed sub-interval $[t_{N-1},t_{N}]$. We define the discrete path length of $a^\pi$ as the sum of its spatial jumps at the internal nodes $t_{1},\dots,t_{N-1}$:
\begin{equation}\label{eq:def_discrete_path_length}
\mathcal{L}_{\pi}:=\sum_{j=0}^{N-2}|a_{j+1}-a_{j}|.
\end{equation}
Given a compactly supported initial datum $u_T \in BV(\R^2)$ satisfying $0 \le u_T \le R$, we define the concatenated entropy trajectory $u^{\pi}$ on $[T,S]$ recursively. At the initial node $t_0 = T$, we set $u^\pi(t_0) := u_T$. Having defined $u^\pi(t_j)$ at the node $t_j$, we define $u^\pi(t)$ on $(t_j, t_{j+1}]$ by applying the fixed-centre entropy semigroup centered at $a_j$:
\begin{equation}\label{eq:def_concatenated_trajectory}
u^{\pi}(t) := S_{t-t_j}^{a_j}\bigl(u^{\pi}(t_j)\bigr) \quad \text{for } t \in (t_j, t_{j+1}].
\end{equation}
Thus, on each closed sub-interval $[t_j, t_{j+1}]$, $u^\pi(t)$ coincides with the unique entropy solution of the fixed-centre equation \eqref{eq:fixed_centre_ivp} with centre $a_j$ and initial datum $u^\pi(t_j)$. Because each semigroup $S_\tau^{a_j}$ produces trajectories that are strongly continuous in $L^1(\R^2)$ (Proposition~\ref{prop:fixed_centre_existence}), we have $\lim_{t \downarrow t_j} u^\pi(t) = S_0^{a_j}(u^\pi(t_j)) = u^\pi(t_j)$ in $L^1(\R^2)$, ensuring that $u^\pi \in C([T,S]; L^1(\R^2))$.

\begin{proposition}[Discrete path-length estimate]\label{prop:discrete_path_length}
Assume $u_{T} \in BV(\R^2)$ is compactly supported with $0 \le u_T \le R$ and possesses a saturated core $u_T=R$ almost everywhere on $B_{r}(a_{0})$. If $\mathcal{L}_{\pi}<r$, the concatenated solution $u^\pi(t)$ retains a saturated core of radius at least $r-\mathcal{L}_{\pi}$ on $B_{r-\mathcal{L}_\pi}(a^\pi(t))$ for all $t\in[T,S]$. Furthermore,
\begin{equation} \label{eq:discrete_V_bound}
\sup_{t\in[T,S]}V_{a^{\pi}(t)}(u^{\pi}(t))\le V_{a_{0}}(u_{T})\exp\left(\frac{\sqrt{2}\mathcal{L}_{\pi}}{r-\mathcal{L}_{\pi}}\right),
\end{equation}
\begin{equation} \label{eq:discrete_TV_bound}
\sup_{t\in[T,S]}\TV(u^{\pi}(t))\le \frac{V_{a_{0}}(u_{T})}{r-\mathcal{L}_{\pi}}\exp\left(\frac{\sqrt{2}\mathcal{L}_{\pi}}{r-\mathcal{L}_{\pi}}\right).
\end{equation}
\end{proposition}

\begin{proof}
On each sub-interval $(t_j, t_{j+1}]$ ($j \in \{0, \dots, N-1\}$), the trajectory $u^\pi(t)$ evolves under the fixed-centre semigroup centered at $a_j$. The fixed-centre directional-variation contraction (Theorem~\ref{thm:directional-variation_contraction}) guarantees $V_{a_j}(u^\pi(t)) \le V_{a_j}(u^\pi(t_j))$ for all $t \in [t_j, t_{j+1}]$. If $u^\pi(t_j) = R$ a.e.\ on a ball $B_{r_j}(a_j)$, the comparison principle with the stationary solution $R1_{B_{r_j}(a_j)}$ (Lemma~\ref{lem:stationary_saturated_core} and Lemma~\ref{lem:localized_kato_positive}) ensures this core is preserved on $(t_j, t_{j+1}]$.

At each internal node $t_{j+1}$ ($j \in \{0, \dots, N-2\}$), the state $u^\pi$ is continuous in $L^1(\R^2)$ while the reference centre shifts from $a_j$ to $a_{j+1}$. For any $x \in B_{r_{j}-|a_{j+1}-a_{j}|}(a_{j+1})$, the triangle inequality gives $|x - a_j| \le |x - a_{j+1}| + |a_{j+1} - a_j| < r_j$. Thus, $B_{r_{j}-|a_{j+1}-a_{j}|}(a_{j+1}) \subset B_{r_{j}}(a_{j})$, implying the core radius decreases by at most $|a_{j+1}-a_j|$. Iterating over the partition, for any $t \in [T, S]$ lying in the $k$-th sub-interval ($t \in [t_k, t_{k+1})$ for $k \in \{0, \dots, N-2\}$, or $t \in [t_{N-1}, t_N]$ for $k = N-1$), where $a^\pi(t) = a_k$, the core radius around $a^\pi(t)$ is bounded below by $r_k := r - \sum_{j=0}^{k-1} |a_{j+1}-a_j| \ge r - \mathcal{L}_\pi > 0$.

Applying Lemma~\ref{lem:change_of_centre} at $t_{j+1}$ with the uniform lower bound $r-\mathcal{L}_\pi$ on the core radius yields
\begin{equation}
V_{a_{j+1}}(u^\pi(t_{j+1})) \le \left(1+\frac{\sqrt{2}|a_{j+1}-a_{j}|}{r-\mathcal{L}_{\pi}}\right)V_{a_{j}}(u^\pi(t_{j+1})) \le \left(1+\frac{\sqrt{2}|a_{j+1}-a_{j}|}{r-\mathcal{L}_{\pi}}\right)V_{a_{j}}(u^\pi(t_j)).
\end{equation}
Iterating this inequality from $t_0 = T$ to any $t$ in the $k$-th sub-interval, using $V_{a_k}(u^\pi(t)) \le V_{a_k}(u^\pi(t_k))$ and applying $1+s \le e^s$ for $s \ge 0$, yields
\begin{equation}
V_{a^{\pi}(t)}(u^{\pi}(t)) \le V_{a_{0}}(u_{T}) \prod_{j=0}^{k-1} \left(1+\frac{\sqrt{2}|a_{j+1}-a_{j}|}{r-\mathcal{L}_{\pi}}\right) \le V_{a_{0}}(u_{T}) \exp\left( \sum_{j=0}^{k-1} \frac{\sqrt{2}|a_{j+1}-a_{j}|}{r-\mathcal{L}_{\pi}} \right),
\end{equation}
which establishes \eqref{eq:discrete_V_bound} since the sum is bounded by $\mathcal{L}_\pi$. Applying Lemma~\ref{lem:core_coercivity} with the minimal persistent radius $r-\mathcal{L}_\pi$ establishes \eqref{eq:discrete_TV_bound}.
\end{proof}

To transition from piecewise constant paths to an absolutely continuous moving centre, we use a lower semicontinuity property that tracks the directional variation as the spatial centre shifts.

\begin{lemma}[Lower semicontinuity with a moving centre] \label{lem:moving_lsc}
Let $u_{n}\to u$ strongly in $L^{1}(\R^{2})$. Suppose there is a common bounded set $K$ such that $\supp u_{n} \subset K$ for all $n$, and the unweighted total variation is uniformly bounded: $\sup_{n} \TV(u_{n})<\infty$. If a sequence of centres converges, $a_{n}\to a$, then
\begin{equation} \label{eq:moving_lsc}
V_{a}(u)\le \liminf_{n\to\infty}V_{a_{n}}(u_{n}).
\end{equation}
\end{lemma}

\begin{proof}
The uniform bound $\sup_{n} \TV(u_{n})<\infty$ guarantees that, up to a subsequence, $Du_n \rightharpoonup Du$ weakly-* in the sense of Radon measures on $\R^{2}$ \cite[Theorem~1.59]{AmbrosioFuscoPallara2000}. On the common bounded support $K$, the affine vector fields converge uniformly in $L^\infty(K; \R^2)$: $X_{a_{n}}\to X_{a}$ and $Y_{a_{n}}\to Y_{a}$ as $n \to \infty$. 

Pairing a weakly-* convergent sequence of measures with uniformly convergent continuous multipliers ensures the weak-* convergence of the product measures:
\begin{equation}
X_{a_{n}}\cdot Du_{n}\rightharpoonup X_{a}\cdot Du, \qquad Y_{a_{n}}\cdot Du_{n}\rightharpoonup Y_{a}\cdot Du.
\end{equation}
By the lower semicontinuity of the total variation of a Radon measure with respect to weak-* convergence \cite[Theorem~1.59]{AmbrosioFuscoPallara2000}, we obtain $|X_a \cdot Du|(\R^{2}) \le \liminf_{n\to\infty} |X_{a_n} \cdot Du_n|(\R^{2})$ and $|Y_a \cdot Du|(\R^{2}) \le \liminf_{n\to\infty} |Y_{a_n} \cdot Du_n|(\R^{2})$. Summing the limits yields \eqref{eq:moving_lsc}.
\end{proof}

\begin{lemma}[Uniform origin-centered support barrier for moving centres]\label{lem:uniform_support_barrier}
Let $D_{a}\ge0$ be a fixed constant and let $D^{\prime}>D_{a}$. Let $u_{T}\in BV(\R^{2})$ satisfy $0\le u_{T}\le R$ and $\supp u_{T}\subset B_{D^{\prime}}(0)$.
\begin{enumerate}
    \item[(i)] \textbf{Uniform confinement for piecewise-constant centres:} For any $S>T$ and any piecewise constant path $a^{\pi}:[T,S]\to\R^{2}$ whose amplitude is bounded by $\sup_{t\in[T,S]}|a^{\pi}(t)|\le D_{a}$, the concatenated fixed-centre entropy solution $u^{\pi}\in C([T,S];L^{1}(\R^{2}))$ with $u^{\pi}(T)=u_{T}$ satisfies
    \begin{equation}\label{eq:piecewise_uniform_support_bound}
    0\le u^{\pi}(t,x)\le R1_{B_{D^{\prime}}(0)}(x) \quad \text{for a.e. } x\in\R^{2} \text{ and all } t\in[T,S].
    \end{equation}
    \item[(ii)] \textbf{Support inheritance:} If $\{u_{n}\}_{n\ge1}\subset C([T,S];L^{1}(\R^{2}))$ is a sequence of trajectories satisfying the bound $0\le u_{n}(t,x)\le R1_{B_{D^{\prime}}(0)}(x)$ almost everywhere, and $u_{n}(t)\to u(t)$ strongly in $L^{1}(\R^{2})$ for every $t\in[T,S]$, then the limit $u$ inherits the confinement \eqref{eq:piecewise_uniform_support_bound}.
\end{enumerate}
\end{lemma}

\begin{proof}
Define the static origin-centered saturated state $U_{D^{\prime}}(x):=R1_{B_{D^{\prime}}(0)}(x)$. For any centre $a\in\overline{B}_{D_{a}}(0)$ and any $x\in\partial B_{D^{\prime}}(0)$, the outward unit normal is $\nu(x)=x/D^{\prime}$. The Cauchy--Schwarz inequality guarantees
\begin{equation}\label{eq:origin_disk_normal_check}
(x-a)\cdot\nu(x)=\frac{|x|^{2}}{D^{\prime}}-\frac{a\cdot x}{D^{\prime}}\ge D^{\prime}-\frac{D_{a}D^{\prime}}{D^{\prime}}=D^{\prime}-D_{a}>0.
\end{equation}
By Lemma~\ref{lem:stationary_saturated_core}, $U_{D^{\prime}}$ is a stationary Kru\v{z}kov entropy solution for the fixed-centre equation with centre $a$, for any $a \in \overline{B}_{D_a}(0)$. 

For (i), $u^{\pi}(T) \le U_{D^{\prime}}$ almost everywhere. Inducting over the sub-intervals $[t_{j},t_{j+1}]$ using the fixed-centre order preservation (Lemma~\ref{lem:localized_kato_positive}) yields $\int_{\R^{2}}(u^{\pi}(t,x)-U_{D^{\prime}}(x))^{+}\,dx = 0$ on $[T,S]$, establishing \eqref{eq:piecewise_uniform_support_bound}. 

For (ii), strong $L^{1}(\R^2)$ convergence implies $\int_{\R^{2}}(u_{n}(t,x)-U_{D^{\prime}}(x))^{+}\,dx \to \int_{\R^{2}}(u(t,x)-U_{D^{\prime}}(x))^{+}\,dx$. Since the sequence integrals are zero, the limit integral vanishes, securing confinement for the limit state.
\end{proof}

\begin{pathstabilitylemma}
Let $T<S$, $A>0$, and let $a_1,a_2:[T,S]\to\R^2$ be bounded piecewise-constant paths with $|a_i(t)|\le A$. Let $u_T\in BV(\R^2)$ be compactly supported, $0\le u_T\le R$, and let $u_i$ be the entropy evolutions obtained by concatenating the fixed-centre semigroups along $a_i$, with the same initial datum $u_T$. Choose $D'>A$ with $\supp u_T\subset B_{D'}(0)$, and choose a common cutoff $\chi\in C_c^\infty(\R^2;[0,1])$ such that $\chi=1$ on a neighborhood of $\overline B_{D'}(0)$. Set $M_T=\|u_T\|_{L^1}$, $L=\Lip(m)$,
\[
\begin{aligned}
 B_1&=\sup_{|q|\le A}\|D[\chi(x)(x-q)]\|_\infty,\\
 B_2&=\sup_{|q|\le A}\|\nabla\diver[\chi(x)(x-q)]\|_\infty,
\end{aligned}
\]
and $B_{T,S}:=(\TV(u_T)+L M_T B_2(S-T))e^{LB_1(S-T)}$. Then
\begin{equation}\label{eq:piecewise_path_stability_uniform}
 \sup_{t\in[T,S]}\|u_1(t)-u_2(t)\|_{L^1(\R^2)}
 \le L\bigl(B_{T,S}+M_T\|\nabla\chi\|_\infty\bigr)
       \int_T^S|a_1(t)-a_2(t)|\,dt.
\end{equation}
The constant does not depend on the number, location, or common refinement of the jump times of $a_1,a_2$.
\end{pathstabilitylemma}

\begin{proof}[Proof of the uniform path-stability lemma]
Write $b_i(t,x)=\chi(x)(x-a_i(t))$. By the common support barrier (Lemma~\ref{lem:uniform_support_barrier}), the concatenated solutions vanish outside $B_{D'}(0)$. They also solve the truncated local equations $\partial_t u_i-\diver(m(u_i)b_i)=0$: inside the common support the fluxes agree, and outside it $m(u_i)=0$; the apparent difference of the Kru\v{z}kov source terms at constant entropy levels is the divergence of a compactly supported test vector field and hence has zero integral.

We establish the estimate without appealing to a stability theorem for spatially unbounded fluxes. Let $m_j\in C_c^\infty(\R)$ and $u_{T,j}\in C_c^\infty(\R^2)$ be the endpoint-preserving mobility approximants and common strict $BV$ initial approximants, satisfying
\[
\begin{gathered}
 0\le u_{T,j}\le R,\qquad L_j:=\Lip(m_j)\le L_*,\qquad m_j\to m\ \text{uniformly},\\
 u_{T,j}\to u_T\ \text{in }L^1,\qquad \TV(u_{T,j})\to\TV(u_T).
\end{gathered}
\]
By Lemma~\ref{lem:endpoint_preserving_mobility}, $L_*:=\sup_j L_j\le2L$ and, crucially, $\limsup_j L_j\le L$. The first estimates below use $L_j$ and are uniform in $j$; taking $j\to\infty$ therefore yields the constant $L$ stated in the lemma.

For $\nu_j\downarrow0$, solve simultaneously, for $i=1,2$,
\begin{equation}\label{eq:path_stability_viscous_pair}
 \partial_t z_i^j-\diver(m_j(z_i^j)b_i)=\nu_j\Delta z_i^j,
 \qquad z_i^j(T)=u_{T,j}.
\end{equation}
The coefficients are smooth in $x$, uniformly bounded together with their first two spatial derivatives, and piecewise constant in $t$. On each interval in the finite common refinement of the jump partitions, the heat-semigroup construction gives a smooth-in-space solution; the spatial traces coincide at successive junctions. Integrating the equation with spatial cutoffs shows conservation of the nonnegative $L^1$ mass. Testing the positive and negative parts as in Proposition~\ref{prop:nonlocal-viscous}, now with a prescribed smooth spatial field, proves $0\le z_i^j\le R$.

For completeness, differentiate \eqref{eq:path_stability_viscous_pair} in space on each open subinterval. With $p=\nabla z_i^j$, apply the multiplier $p/(|p|^2+\delta^2)^{1/2}$, integrate against expanding compact spatial cutoffs, and send the regularization radius $\delta\downarrow0$ before the viscous limit. The transport terms containing $m_j''$ cancel with the divergence of $m_j'(z_i^j)|p|b_i$, while the remaining terms are controlled by $L_j|Db_i||p|$ and $|m_j(z_i^j)|\,|\nabla\diver b_i|$. Convexity makes the viscous contribution nonpositive. Since $|m_j(s)|\le L_j s$ for $0\le s\le R$, this yields on each subinterval
\[
 \frac{d}{dt}\TV(z_i^j(t))\le L_j B_1\TV(z_i^j(t))
      +L_j\|u_{T,j}\|_{L^1}B_2
 \quad\text{for a.e. }t.
\]
Both $z_i^j$ and their spatial gradients have matching traces at time jumps, so summing the integrated inequalities over the common refinement introduces no jump-dependent constant. Consequently,
\begin{equation}\label{eq:path_stability_viscous_BV}
 \sup_{i=1,2}\sup_{t\in[T,S]}\TV(z_i^j(t))
 \le\Big(\TV(u_{T,j})+L_j\|u_{T,j}\|_1 B_2(S-T)\Big)e^{L_j B_1(S-T)}.
\end{equation}
All constants are independent of the number of path jumps, of the viscosity, and of $j$.

Subtract the two viscous equations, temporarily using $-m_j(s)b_1$ as the common flux. The sole forcing term is the spatial divergence of $\chi m_j(z_2^j)(a_1-a_2)$. To justify the parabolic Kato step directly, set $w=z_1^j-z_2^j$ and multiply the difference equation by $\eta_\delta'(w)$, where $\eta_\delta$ is a smooth convex approximation of $|w|$. The common-flux term is controlled, after integration by parts, by $L_j\|b_1\|_\infty\int\eta_\delta''(w)|w||\nabla w|$, which tends to zero at fixed $j$ by dominated convergence (choose $|w|\eta_\delta''(w)\le C$ and note pointwise convergence to zero); the viscosity produces $-\nu_j\int\eta_\delta''(w)|\nabla w|^2\le0$. The coefficient perturbation is bounded by $\int|\diver(\chi m_j(z_2^j)(a_1-a_2))|$. Sending $\delta\downarrow0$ therefore gives for a.e. $t$ away from the finitely many jumps
\begin{align}\label{eq:path_stability_Kato_precise}
\frac{d}{dt}\|z_1^j(t)-z_2^j(t)\|_1
&\le |a_1(t)-a_2(t)|\int_{\R^2}|D(\chi m_j(z_2^j(t)))|\,dx\notag\\
&\le L_j|a_1(t)-a_2(t)|
   \bigl(\TV(z_2^j(t))+\|\nabla\chi\|_\infty\|u_{T,j}\|_1\bigr).
\end{align}
There are no junction terms because both viscous solutions are continuous into $L^1$ at each jump, and they have identical initial data. Integrate over the common refinement and apply \eqref{eq:path_stability_viscous_BV}. The resulting bound is uniform in $j$ and in both partitions.

Finally take $j\to\infty$ with $\nu_j\downarrow0$. On the fixed interval $[T,S]$, \eqref{eq:path_stability_viscous_BV} gives uniform spatial $BV$ bounds. The equations provide a uniform $W^{-1,1}$ time modulus, because the fluxes satisfy $\|m_j(z_i^j)b_i\|_1\le L_*M_*\|b_i\|_\infty$ and $\nu_j\|\nabla z_i^j\|_1\le\nu_j\sup_t\TV(z_i^j(t))$. Spatial mollification converts this into an $L^1$ temporal modulus. Uniform tightness follows by testing with cutoffs outside large balls: the drift is compactly supported, while the diffusive tail contribution is bounded by $C\nu_j M_*/H^2$ for a cutoff of radius $H$. The Kolmogorov--Riesz and Arzel\`a--Ascoli criteria therefore yield, after extraction, convergence in $C([T,S];L^1(\R^2))$. Passing to the weak and entropy inequalities uses $m_j\to m$, $\nu_j\to0$, the uniform $L^1$ convergence, and a two-sided approximation of exceptional flat entropy levels. The limiting solution is the unique entropy evolution for the bounded, piecewise-time-dependent flux $-m(s)b_i$, by the standard Kato contraction on each interval of the common refinement. The confined concatenation $u_i$ is itself such an entropy solution, as observed at the start; uniqueness identifies the limit with $u_i$. Thus the \emph{same} viscous sequences converge to $u_1$ and $u_2$, even when their jump partitions are different. Passing to the limit in the integrated inequality, using $\limsup_j L_j\le L$ and $\|u_{T,j}\|_1\to M_T$, proves \eqref{eq:piecewise_path_stability_uniform} with precisely the constant displayed in the statement. \qedhere
\end{proof}

\begin{theorem}[Absolutely continuous prescribed path]\label{thm:absolutely_continuous_prescribed_path}
Let $a\in AC([T,\infty);\R^{2})$. Assume the remaining path length is bounded by $r>0$:
\begin{equation}\label{eq:remaining_path_length_strict_r}
\mathcal{L}_{a}(T):=\int_{T}^{\infty}|\dot{a}(t)|\,dt < r.
\end{equation}
Let $u_{T}\in BV(\R^{2})$ be compactly supported, satisfy $0\le u_{T}\le R$, and possess a saturated core $u_{T}=R$ a.e.\ on $B_{r}(a(T))$. Then the prescribed-path entropy problem
\begin{equation} \label{eq:prescribed_path_pde}
\partial_{t}u-\diver\bigl(m(u)(x-a(t))\bigr)=0, \qquad u(T)=u_{T},
\end{equation}
admits a unique compactly supported entropy solution $u\in C([T,\infty);L^{1}(\R^{2}))$. For every $t\ge T$, the solution retains a saturated core $u(t)=R$ a.e.\ on $B_{r-\mathcal{L}_{a}(T,t)}(a(t))$, where $\mathcal{L}_{a}(T,t) := \int_{T}^{t} |\dot{a}(s)|\,ds$. Furthermore,
\begin{equation}\label{eq:prescribed_path_V_bound_final}
\sup_{t\ge T}V_{a(t)}(u(t))\le V_{a(T)}(u_{T})\exp\left(\frac{\sqrt{2}\mathcal{L}_{a}(T)}{r-\mathcal{L}_{a}(T)}\right),
\end{equation}
\begin{equation}\label{eq:prescribed_path_TV_bound_final}
\sup_{t\ge T}\TV(u(t))\le \frac{V_{a(T)}(u_{T})}{r-\mathcal{L}_{a}(T)}\exp\left(\frac{\sqrt{2}\mathcal{L}_{a}(T)}{r-\mathcal{L}_{a}(T)}\right).
\end{equation}
\end{theorem}

\begin{proof}
Fix $S>T$. Let $a_{n}(t)$ be a sequence of piecewise constant paths sampled from $a(t)$ at the nodes $T = t_0^n < t_1^n < \dots < t_{N_n}^n = S$ of a sequence of uniform partitions $\pi_n$ of $[T,S]$ such that $\|a_n - a\|_{L^\infty(T,S)} \to 0$. For any intermediate time $t \in [T, S]$ lying in the sub-interval $[t_m^n, t_{m+1}^n)$, the cumulative discrete path length of $a_n$ up to time $t$ is given by the partial sum
\begin{equation}
\mathcal{L}_{\pi_n}(T,t) := \sum_{j=0}^{m-1}\bigl|a_n(t_{j+1}^{n})-a_n(t_{j}^{n})\bigr|.
\end{equation}
Because $a_n$ is constructed by sampling the absolutely continuous path $a$ at the partition nodes, we have $a_n(t_{j+1}^n) = a(t_{j+1}^n)$ and $a_n(t_j^n) = a(t_j^n)$. Applying the Fundamental Theorem of Calculus to $a$ on each sub-interval gives
\begin{equation}
\bigl|a_n(t_{j+1}^{n})-a_n(t_{j}^{n})\bigr| = \bigl|a(t_{j+1}^{n})-a(t_{j}^{n})\bigr| = \left| \int_{t_j^n}^{t_{j+1}^n} \dot{a}(s)\,ds \right| \le \int_{t_j^n}^{t_{j+1}^n} |\dot{a}(s)|\,ds.
\end{equation}
Summing over $j=0,\dots,m-1$ with $t_m^n \le t$ yields
\begin{equation}
\mathcal{L}_{\pi_n}(T,t) \le \int_{T}^{t_m^n}|\dot{a}(s)|\,ds \le \int_{T}^{t}|\dot{a}(s)|\,ds = \mathcal{L}_{a}(T,t) \le \mathcal{L}_{a}(T) < r.
\end{equation}

Since $\mathcal{L}_{a}(T)=\int_{T}^{\infty}|\dot{a}(s)|\,ds<r$, the Fundamental Theorem of Calculus gives $|a(t)|\le|a(T)|+\mathcal{L}_{a}(T)<|a(T)|+r=:D_{a}<\infty$ for all $t\ge T$. The bound $|a(t)| \le D_a < \infty$ therefore applies to $a_n(t)$, and choosing $D^{\prime}>D_{a}$ such that $\supp u_{T}\subset B_{D^{\prime}}(0)$, Lemma~\ref{lem:uniform_support_barrier} confines all concatenated solutions $u_n$ to a common compact ball $\overline{B}_{D^{\prime}}(0)$. Proposition~\ref{prop:discrete_path_length} guarantees a uniform total variation bound $\sup_{n}\sup_{t\in[T,S]}\TV(u_{n}(t))\le B_{*}$ via core coercivity.

Because $u_n \in C([T,S]; L^1(\R^2))$ has no discontinuities across the partition nodes, summing the fixed-centre Kru\v{z}kov entropy inequalities over the sub-intervals cancels all internal temporal boundary traces. Thus, $u_n$ satisfies the global entropy inequality on $(T,S) \times \R^2$ for the moving flux $P_n(t,x,s) = -m(s)(x-a_n(t))$.

We now apply the independent uniform path-stability lemma proved above, with the same spatial cutoff $\chi$ for all the sampled paths. Its assumptions hold with $A=D_a$, common support radius $D'$, and initial datum $u_T$, so that
\begin{equation}\label{eq:gronwall_bound}
 \sup_{t\in[T,S]}\|u_n(t)-u_k(t)\|_{L^1(\R^2)}
 \le C_S\int_T^S|a_n(\tau)-a_k(\tau)|\,d\tau
 \le C_S(S-T)\|a_n-a_k\|_{L^\infty(T,S)},
\end{equation}
where $C_S$ depends only on $u_T$, the common cutoff, the bound on the centre, $m$, and $S-T$, and is independent of $n,k$ and of both partitions. The coefficient comparison was carried out for bounded truncated fields and the limiting solutions were identified with the original concatenations, so no unbounded-flux stability hypothesis is invoked.

As $a_n$ is Cauchy in $L^\infty(T,S)$, $u_n$ is Cauchy in $C([T,S]; L^1(\R^2))$, securing a strong limit $u$. By Lemma~\ref{lem:uniform_support_barrier}(ii), the limit retains the common support bound $\supp u(t)\subset\overline B_{D^{\prime}}(0)$ for every $t\in[T,S]$. Passing to the limit in the weak formulation utilizes the decomposition $m(u_n)(x-a_n) - m(u)(x-a) = m(u_n)(a-a_n) + (x-a)(m(u_n)-m(u))$, which converges strongly in $L^1((T,S)\times\R^2)$. 

To pass to the limit in the Kru\v{z}kov entropy inequalities, let $E_\varphi$ be the countable set of flat levels of $u$ on the support of a test function $\varphi$. For $k \notin E_\varphi$, $\operatorname{sign}(u_n - k) \to \operatorname{sign}(u - k)$ almost everywhere, allowing the Dominated Convergence Theorem to establish the inequality. For $k_0 \in E_\varphi$, symmetric two-sided limits $k_\ell^\pm \to k_0$ algebraically cancel the discontinuous source term on the level set $\{u = k_0\}$, extending the inequality to all $k \in \R$.

For any compact set $K \subset B_{r-\mathcal{L}_{a}(T,t)}(a(t))$, $K \subset B_{r-\mathcal{L}_{\pi_n}(T,t)}(a_n(t))$ for all sufficiently large $n$. Since $u_n = R$ a.e.\ on this set and $u_n \to u$ in $L^1(K)$, the core property transfers to $u$. Applying Lemma~\ref{lem:moving_lsc} and Lemma~\ref{lem:core_coercivity} to the strongly converging sequence with $\sup_n \TV(u_n) \le B_*$ transfers the variation bounds \eqref{eq:prescribed_path_V_bound_final} and \eqref{eq:prescribed_path_TV_bound_final} to $u(t)$. Continuation to $[T,\infty)$ follows from the uniqueness established by the prescribed-path contraction (Proposition~\ref{prop:prescribed_path_contraction}).
\end{proof}

\begin{lemma}[Signed mobility difference] \label{lem:signed_mobility_diff}
Define $H(r,s):=\operatorname{sign}(r-s)(m(r)-m(s))$ with $H(r,r):=0$. For any $r_1, s_1, r_2, s_2 \in [0, R]$,
\begin{equation}
|H(r_{1},s_{1})-H(r_{2},s_{2})|\le \Lip(m)(|r_{1}-r_{2}|+|s_{1}-s_{2}|).
\end{equation}
\end{lemma}

\begin{proof}
If the pairs $(r_1, s_1)$ and $(r_2, s_2)$ lie on the same side of the diagonal $r=s$, the estimate reduces to the Lipschitz bound for $\pm(m(r)-m(s))$. If they cross the diagonal, assume without loss of generality that $r_{1}\ge s_{1}$ and $r_{2}\le s_{2}$. Then
\begin{align}
|H(r_{1},s_{1})-H(r_{2},s_{2})| &\le |m(r_{1})-m(s_{1})|+|m(s_{2})-m(r_{2})| \nonumber \\
&\le \Lip(m)\bigl[(r_{1}-s_{1})+(s_{2}-r_{2})\bigr] \nonumber \\
&= \Lip(m)\bigl[(r_{1}-r_{2})+(s_{2}-s_{1})\bigr] \nonumber \\
&\le \Lip(m)\bigl(|r_{1}-r_{2}|+|s_{1}-s_{2}|\bigr).
\end{align}
\end{proof}

\begin{lemma}[Temporal diagonalization] \label{lem:temporal_diagonalization}
Let $K\subset\R^{2}$ be a compact set, $f\in C([T,S];L^{1}(K))$, and let $\rho_{\eps}$ be a standard symmetric temporal mollifier supported in $[-\eps,\eps]$. Then
\begin{equation}
\lim_{\eps\to 0}\int_{T}^{S}\!\int_{T}^{S}\rho_{\eps}(t-s)\|f(t)-f(s)\|_{L^{1}(K)}\,ds\,dt = 0.
\end{equation}
\end{lemma}

\begin{proof}
Since $f\in C([T,S];L^{1}(K))$, it is uniformly continuous. Define the uniform modulus of continuity $\omega_{f}(h):=\sup_{|t-s|\le h}\|f(t)-f(s)\|_{L^{1}(K)}$, which satisfies $\omega_{f}(h)\to 0$ as $h\to 0$. Because $\supp\rho_{\eps}\subset[-\eps,\eps]$, the integrand is non-zero only when $|t-s|\le \eps$. Bounding the $L^1$ distance by $\omega_{f}(\eps)$ yields
\begin{equation}
\int_{T}^{S} \left( \int_{T}^{S}\rho_{\eps}(t-s)\,ds \right) \omega_{f}(\eps)\,dt \le (S-T)\omega_{f}(\eps).
\end{equation}
The integral vanishes as $\eps\to 0$.
\end{proof}

\begin{lemma}[Kato inequality for a moving affine centre] \label{lem:moving_centre_kato}
Let $a\in AC_{\mathrm{loc}}([T,\infty);\R^{2})$, and let $u,v\in C([T,S];L_{\mathrm{loc}}^{1}(\R^{2}))\cap L^{\infty}((T,S)\times\R^{2})$ be entropy solutions of $\partial_t w - \diver(m(w)(x-a(t))) = 0$, taking values in $[0, R]$. The following inequality holds in the sense of distributions on $(T,S)\times\R^{2}$:
\begin{equation} \label{eq:moving_kato_dist}
\partial_{t}|u-v|-\diver_{x}\bigl[\operatorname{sign}(u-v)\bigl(m(u)-m(v)\bigr)\bigl(x-a(t)\bigr)\bigr]\le0.
\end{equation}
\end{lemma}

\begin{proof}
Define the affine flux $P(t,x,z):=-m(z)(x-a(t))$. Let $\phi\in C_{c}^{1}((T,S)\times\R^{2})$ be non-negative. Choose a temporal mollifier $\rho_{\eps}$ supported in $[-\eps,\eps]$ with $\eps$ smaller than half the distance from $\supp\phi$ to $\{T,S\}$, and a spatial mollifier $\omega_{\delta}$. Define the separable test function
\begin{equation}
\Phi_{\eps,\delta}(t,x,s,y):=\phi\left(\frac{t+s}{2},\frac{x+y}{2}\right)\rho_{\eps}(t-s)\omega_{\delta}(x-y).
\end{equation}
Evaluating the Kru\v{z}kov entropy inequality for $u(t,x)$ against the frozen state $v(s,y)$, adding it to the inequality for $v(s,y)$ evaluated at $u(t,x)$, and setting $\sigma:=\operatorname{sign}(u(t,x)-v(s,y))$ gives
\begin{align} \label{eq:kato_four_integrals}
0\le {}&\iiiint\Bigl[|u-v|(\partial_{t}+\partial_{s})\Phi_{\eps,\delta} + \sigma\bigl[P(t,x,u)-P(t,x,v)\bigr]\cdot\nabla_{x}\Phi_{\eps,\delta} \nonumber \\
&\qquad-\sigma\bigl[P(s,y,v)-P(s,y,u)\bigr]\cdot\nabla_{y}\Phi_{\eps,\delta} \nonumber \\
&\qquad+ \sigma\bigl[\diver_{y}P(s,y,u)-\diver_{x}P(t,x,v)\bigr]\Phi_{\eps,\delta}\Bigr]\,dx\,dt\,dy\,ds.
\end{align}
Applying the chain rule with coordinates $\bar{t}=(t+s)/2$, $\bar{x}=(x+y)/2$, and $z=x-y$, the flux gradient terms expand into
\begin{align} \label{eq:kato_bracket_expansion}
&\frac{1}{2}\iiiint\sigma\bigl[P(t,x,u)-P(t,x,v)+P(s,y,u)-P(s,y,v)\bigr]\cdot\nabla\phi(\bar{t},\bar{x})\rho_{\eps}(t-s)\omega_{\delta}(z) \nonumber \\
&\quad+\iiiint\sigma\bigl[P(t,x,u)-P(s,y,u)-P(t,x,v)+P(s,y,v)\bigr]\cdot\nabla\omega_{\delta}(z)\phi(\bar{t},\bar{x})\rho_{\eps}(t-s).
\end{align}
Substituting $P(t,x,q)=-m(q)(x-a(t))$ and $H(u,v)=\sigma(m(u)-m(v))$, the cross-variable flux difference identity $P(t,x,q)-P(s,y,q)=-m(q)(x-y)+m(q)(a(t)-a(s))$ isolates the moving-centre discrepancy into the remainder
\begin{align}
\mathcal{R}_{\eps,\delta}:={}&-\frac{1}{2}\iiiint H(u,v)\bigl(a(t)-a(s)\bigr)\cdot\nabla\phi(\bar{t},\bar{x})\rho_{\eps}(t-s)\omega_{\delta}(z) \nonumber \\
&+\iiiint H(u,v)\bigl(a(t)-a(s)\bigr)\cdot\nabla\omega_{\delta}(z)\phi(\bar{t},\bar{x})\rho_{\eps}(t-s).
\end{align}
Because $|H(u,v)|\le 2R\Lip(m)$, $a\in AC_{\mathrm{loc}}([T,\infty);\R^{2})$, and $\supp\rho_{\eps}\subset[-\eps,\eps]$, for any $(t,s)\in(T,S)^{2}$ with $|t-s|\le\eps$ and any $\tau\in[\min\{t,s\},\max\{t,s\}]$ we have $|\tau-s|\le|t-s|\le\eps$, so that
\begin{equation}
|a(t)-a(s)|\le\int_{\min\{t,s\}}^{\max\{t,s\}}|\dot{a}(\tau)|\,d\tau\le\int_{T}^{S}1_{[\tau-\eps,\tau+\eps]}(s)|\dot{a}(\tau)|\,d\tau.
\end{equation}
Applying Fubini's theorem and $\int_{\R}\rho_{\eps}(t-s)\,dt=1$ provides the bound
\begin{equation}
|\mathcal{R}_{\eps,\delta}|\le C_{\phi,\delta,m,R} \iint |a(t)-a(s)|\rho_{\eps}(t-s)\,dt\,ds \le 2\eps C_{\phi,\delta,m,R}\|\dot{a}\|_{L^{1}(T,S)}.
\end{equation}
Thus $\mathcal{R}_{\eps,\delta}\to 0$ as $\eps\to 0$ for every fixed spatial smoothing radius $\delta>0$. It is essential that this temporal limit is taken first: at fixed $\delta$ the spatial mollifier and its derivatives are bounded, so the $L^1_{\mathrm{loc}}$ time continuity of $u,v$ (Lemma~\ref{lem:temporal_diagonalization}) and continuity of the test functions justify collapse onto $s=t$. Only after this step do we send $\delta\downarrow0$. 

Substituting $\diver_x P(t,x,q) = -2m(q)$ into the collapsed spatial divergence terms groups them with the $\nabla\omega_\delta$ component into the integral
\begin{equation} \label{eq:kato_spatial_q_delta}
-\int_{T}^{S}\!\iint H\bigl(u(t,x),v(t,y)\bigr)q_{\delta}(x-y)\phi\left(t,\frac{x+y}{2}\right)dx\,dy\,dt,
\end{equation}
where $q_{\delta}(z):=2\omega_{\delta}(z)+z\cdot\nabla\omega_{\delta}(z)$. Integrating by parts verifies $\int_{\R^2} q_{\delta}(z)\,dz = 0$, while $\|q_\delta\|_{L^1(\R^2)}$ is independent of $\delta$. Changing variables $y=x-z$, subtracting the zero-mean integral
\[
\int_{\R^2} q_\delta(z) H\bigl(u(t,x),v(t,x)\bigr)\phi(t,x)\,dz = 0,
\]
and applying Lemma~\ref{lem:signed_mobility_diff} splits the error into two terms. The first, arising from $v(t,x-z)-v(t,x)$, is bounded by the $L^1$ translation modulus of $v$ on a fixed enlarged compact neighborhood of $\supp\phi$, integrated against $|q_\delta(z)|$. It tends to zero after time integration by continuity of spatial translations in $L^1_{\mathrm{loc}}$ and dominated convergence. The second, arising from $\phi(t,x-z/2)-\phi(t,x)$, satisfies the explicit bound
\[
 C(R,m,\phi,S-T)\int_{\R^2}|z|\,|q_\delta(z)|\,dz
 = C(R,m,\phi,S-T)\,\delta\int_{\R^2}|w|\,|q_1(w)|\,dw=O(\delta).
\]
Here the spatial integration is restricted to the fixed compact neighborhood on which the test is nonzero, and the bound $|H(u,v)|\le2R\Lip(m)$ applies. Thus both errors vanish, and \eqref{eq:kato_spatial_q_delta} tends to zero as $\delta\to0$. Passing to the limit $\delta\to0$ in the remaining terms recovers the distributional formulation \eqref{eq:moving_kato_dist}.
\end{proof}

\begin{proposition}[Prescribed-path contraction in the confined class] \label{prop:prescribed_path_contraction}
Let $a\in AC_{\mathrm{loc}}([T,\infty);\R^{2})$. Let $u, v\in C([T,S];L^{1}(\R^{2}))$ be entropy solutions of \eqref{eq:prescribed_path_pde} taking values in $[0, R]$. Assume there is a compact set $K$ such that $\supp u(t) \cup \supp v(t) \subset K$ for all $t\in[T,S]$. Then, for all $T\le s\le t\le S$,
\begin{equation} \label{eq:global_prescribed_contraction}
\|u(t)-v(t)\|_{L^{1}(\R^{2})}\le\|u(s)-v(s)\|_{L^{1}(\R^{2})}.
\end{equation}
\end{proposition}

\begin{proof}
Choose a static smooth cutoff $\chi\in C_{c}^{\infty}(\R^{2})$ such that $0\le\chi\le1$ and $\chi\equiv1$ on an open neighborhood of $K$. Consequently, $\nabla\chi \equiv 0$ on $K$. First fix $T<s<t<S$ and let $\theta_{n}$ be smooth non-negative temporal cutoffs approximating $1_{[s,t]}$. Testing \eqref{eq:moving_kato_dist} against $\chi(x)\theta_{n}(\tau)$ and integrating by parts yields
\begin{align}
&-\int_{0}^{\infty}\!\int_{\R^{2}}|u-v|\chi(x)\theta_{n}^{\prime}(\tau)\,dx\,d\tau \nonumber \\
&\qquad - \int_{0}^{\infty}\!\int_{\R^{2}}\operatorname{sign}(u-v)\bigl(m(u)-m(v)\bigr)\bigl(x-a(\tau)\bigr)\cdot\nabla\chi(x)\theta_{n}(\tau)\,dx\,d\tau \le 0.
\end{align}
On $\supp\nabla\chi$, $u=v=0$ almost everywhere, so $m(u)-m(v)=0$ and the spatial boundary flux vanishes. Passing to the limit $n\to\infty$ evaluates the time derivative at the interval boundaries. The strong $L^1$ continuity of the trajectories extends the inequality to $s=T$ and $t=S$ by taking one-sided limits, establishing \eqref{eq:global_prescribed_contraction}.
\end{proof}

\subsection{The Core-Filling Auxiliary Trajectory}\label{subsec:core_filling}

To establish the relative compactness of the trajectory without imposing unproved spatial regularity on the physical state, we approximate the solution at late times with an auxiliary $BV$ state containing an inserted saturated core. This core neutralizes the origin degeneracy along the prescribed path, enabling uniform bounded variation estimates on the auxiliary evolution.

\begin{lemma}[$BV$ approximation with a saturated core]\label{lem:bv_approximation_core}
Suppose Assumption~\ref{assump:confinement} holds. Fix $T\ge0$, a target core radius $r>0$, and an error tolerance $\eta>0$. There exists an auxiliary state $v_{T}\in BV(\R^{2})$ with compact support, satisfying
\begin{equation}
0\le v_{T}\le R \quad \text{a.e. in } \R^2, \qquad v_{T}=R \quad \text{a.e. on } B_{r}(c(T)),
\end{equation}
where $c(T)=\int_{\R^{2}}x\rho(T,x)\,dx$ denotes the barycentre of $\rho(T)$, and bounding the $L^1$ distance to the state $\rho(T)$:
\begin{equation} \label{eq:core_filling_error}
\|v_{T}-\rho(T)\|_{L^1(\R^2)}\le \eta+\pi Rr^{2}.
\end{equation}
\end{lemma}

\begin{proof}
By Assumption~\ref{assump:confinement}, $\rho(T) \in L^1(\R^2)$ is compactly supported with $0 \le \rho(T) \le R$. Standard spatial mollification provides a smooth, compactly supported function $g\in C_{c}^{\infty}(\R^{2})$ such that $0\le g\le R$ and $\|g-\rho(T)\|_{L^{1}(\R^{2})}<\eta$. Because $g \in C_c^\infty(\R^2)$, it follows that $g \in BV(\R^2)$.

We define the auxiliary state via the pointwise maximum:
\begin{equation}
v_{T} := \max\bigl\{g,\, R1_{B_{r}(c(T))}\bigr\}.
\end{equation}
The indicator function of the open ball $B_r(c(T))$ belongs to $BV(\R^2)$. Since the space $BV(\R^2)$ is a lattice closed under pointwise maxima \cite[Theorem~3.96]{AmbrosioFuscoPallara2000}, $v_T \in BV(\R^2)$. By construction, $v_T$ is compactly supported, satisfies $0 \le v_T \le R$ globally, and equals $R$ on $B_r(c(T))$. 

The $L^1$ error introduced by the auxiliary core construction is bounded by the mass of the inserted disk:
\begin{equation}
0 \le v_T(x) - g(x) \le R1_{B_r(c(T))}(x) \quad \text{for a.e. } x \in \R^2.
\end{equation}
Integrating over $\R^2$ yields $\|v_T - g\|_{L^1(\R^2)} \le \pi Rr^2$. The triangle inequality establishes \eqref{eq:core_filling_error}.
\end{proof}

\begin{theorem}[Relative compactness of the full trajectory]\label{thm:relative_compactness}
Under Assumptions~\ref{assump:mobility} and~\ref{assump:confinement}, the trajectory $\{\rho(t):t\ge0\}$ is relatively compact in $L^{1}(\R^{2})$.
\end{theorem}

\begin{proof}
Let $\eps>0$. We fix an artificial core radius $r>0$ and a mollification tolerance $\eta>0$ satisfying
\begin{equation} \label{eq:compactness_error_split}
\eta+\pi Rr^{2} < \frac{\eps}{3}.
\end{equation}
By Proposition~\ref{prop:finite_path_length}, the total path length of the moving barycentre $c(t)$ is finite. We select a threshold time $T\ge0$ sufficiently large such that the remaining future path length is bounded by half the chosen core radius:
\begin{equation}
\int_{T}^{\infty}|\dot{c}(t)|\,dt < \frac{r}{2}.
\end{equation}

At time $T$, we apply Lemma~\ref{lem:bv_approximation_core} to construct the auxiliary state $v_{T}\in BV(\R^{2})$ satisfying $\|v_T - \rho(T)\|_{L^1(\R^2)} < \eps/3$ and possessing a saturated core of radius $r$ centered at $c(T)$. Let $v(t)$ be the unique compactly supported entropy solution to the prescribed-path equation
\begin{equation}
\partial_{t}v-\diver\bigl(m(v)(x-c(t))\bigr)=0, \qquad v(T)=v_{T},
\end{equation}
driven by the barycentre path $a(t) := c(t)$ on $[T, \infty)$. Because the remaining path length satisfies $\int_T^\infty |\dot{c}(t)|\,dt < r/2 < r$, Theorem~\ref{thm:absolutely_continuous_prescribed_path} applies. It guarantees that the auxiliary trajectory $v(t)$ retains a saturated core of radius at least $r/2$ around $c(t)$ for all $t \ge T$, establishing the uniform total variation bound $\sup_{t\ge T}\TV(v(t))<\infty$.

By Lemma~\ref{lem:uniform_support_barrier}, the uniform bound $|c(t)| \le D$ implies that all states $v(t)$ are supported within a fixed ball $\overline{B}_{D^{\prime}}(0)$ for all $t \ge T$. The Rellich--Kondrachov compactness theorem \cite[Theorem~3.23]{AmbrosioFuscoPallara2000} guarantees that bounded sets in $BV(B_{D'}(0))$ are precompact in $L^1(B_{D'}(0))$. Therefore, the closure
\begin{equation}
\mathcal{K} := \overline{\{v(t):t\ge T\}}^{L^{1}}
\end{equation}
is a compact subset of $L^{1}(\R^{2})$. Consequently, the auxiliary tail $\{v(t):t\ge T\}$ can be covered by finitely many $L^{1}(\R^{2})$-balls of radius $\eps/3$.

By Proposition~\ref{prop:exact_affine_reduction}, the trajectory $\rho(t)$ is an entropy solution of the same prescribed-path equation with moving centre $c(t)$ on $[T,\infty)$. Applying the prescribed-path $L^{1}$ contraction (Proposition~\ref{prop:prescribed_path_contraction}) ensures the $L^1$ distance is monotonically non-increasing:
\begin{equation}
\sup_{t\ge T}\|\rho(t)-v(t)\|_{L^{1}(\R^{2})} \le \|\rho(T)-v_{T}\|_{L^{1}(\R^{2})} < \frac{\eps}{3}.
\end{equation}
By the triangle inequality, the $L^{1}(\R^{2})$-balls of radius $2\eps/3$ sharing the same centres cover the tail $\{\rho(t):t\ge T\}$.

Because $\rho\in C([0,\infty);L^{1}(\R^{2}))$, the initial finite-time segment $\{\rho(t):t\in[0,T]\}$ is the continuous image of a compact interval and is therefore compact in $L^{1}(\R^{2})$, admitting a finite covering by $L^{1}(\R^{2})$-balls of radius $\eps$. The union of these two coverings covers the full trajectory $\{\rho(t):t\ge0\}$ with finitely many $L^{1}(\R^{2})$-balls of radius $\eps$. Since $\eps>0$ is arbitrary and $L^{1}(\R^{2})$ is complete, the trajectory is totally bounded and relatively compact in $L^{1}(\R^{2})$.
\end{proof}

\subsection{Full-Trajectory Convergence and Asymptotic Binarity}\label{subsec:binarity}\label{subsec:phase_separation}

We now establish the functional framework required to control the asymptotic behavior of the trajectory. Because the state $\rho(t)$ is not assumed to be of bounded variation, its weak time derivative is measured in the dual space $W^{-1,1}(\R^2)$, defined as a subspace of the dual of $W^{1,\infty}(\R^2)$. The norm is given by
\begin{equation}
\|T\|_{W^{-1,1}(\R^2)} := \sup \bigl\{ \langle T, \phi \rangle : \phi \in W^{1,\infty}(\R^2), \, \|\phi\|_{W^{1,\infty}(\R^2)} \le 1 \bigr\}.
\end{equation}
For a distribution $T = \diver F$ generated by a physical flux $F \in L^1(\R^2; \R^2)$, the definition of the distributional divergence $\langle \diver F, \phi \rangle = -\int_{\R^2} F \cdot \nabla \phi \,dx$ and H\"older's inequality provide the bound $\|\diver F\|_{W^{-1,1}(\R^2)} \le \|F\|_{L^1(\R^2)}$. Furthermore, on any subset $\mathcal{K} \subset L^1(\R^2)$ that is compact in the strong $L^1$ topology, the weak $W^{-1,1}$ topology and the strong $L^1$ topology coincide.

\begin{theorem}[Full-trajectory convergence and asymptotic binarity]\label{thm:binarity}\label{thm:phase_separation}
Under Assumptions~\ref{assump:mobility} and~\ref{assump:confinement}, the trajectory $\rho(t)$ converges strongly in $L^{1}(\R^{2})$ (and consequently weakly-* in $L^{\infty}(\R^{2})$) as $t\to\infty$ to a unique binary state
\begin{equation}\label{eq:full_trajectory_L1_convergence}
\rho_{\infty}=R1_{E}
\end{equation}
for some measurable set $E\subset B_{D}(0)$ of volume $|E|=1/R$, unique up to $\mathcal{L}^{2}$-null sets.
\end{theorem}

\begin{proof}
By Proposition~\ref{prop:finite_path_length}, the unweighted total mobility is globally integrable:
\begin{equation} \label{eq:M_integrable}
\int_{0}^{\infty} M(t)\,dt < \infty.
\end{equation}
By Proposition~\ref{prop:exact_affine_reduction}, the state $\rho(t)$ solves the affine-reduced equation $\partial_{t} \rho = \diver F(t,x)$ in $\mathcal{D}^{\prime}((0,\infty)\times\R^{2})$, where $F(t,x) := m(\rho(t,x))(x-c(t))$. The active mass is confined to $B_{D}(0)$ for all time, and $0\le\rho(t,x)\le R$. Because $m(0)=0$ and $m \ge 0$, the flux $F(t,\cdot)$ vanishes almost everywhere outside $B_{D}(0)$. Since $|c(t)|\le D$, the relative distance on $B_{D}(0)$ satisfies $|x-c(t)|\le 2D$. For every $t \ge 0$, the flux is bounded in $L^{1}(\R^{2})$ by the unweighted total mobility:
\begin{equation}\label{eq:CF_uniform_flux_L1_bound}
\|F(t, \cdot)\|_{L^{1}(\R^{2})} = \int_{B_{D}(0)} m(\rho(t,x))|x-c(t)|\,dx \le 2D \int_{B_{D}(0)} m(\rho(t,x))\,dx = 2D\,M(t).
\end{equation}

To translate this $L^{1}$ flux bound into a dual-space increment estimate for $t\mapsto\rho(t)$ in $W^{-1,1}(\R^{2})$, fix a test function $\phi\in W^{1,\infty}(\R^{2})$ satisfying $\|\phi\|_{W^{1,\infty}(\R^{2})}\le 1$. Choose a smooth radial cutoff $\chi\in C_{c}^{\infty}(\R^{2};[0,1])$ such that $\chi\equiv 1$ on $B_{D+1}(0)$ and $\chi\equiv 0$ outside $B_{D+2}(0)$. Let $\eta_{\eps}$ be a standard symmetric spatial mollifier on $\R^{2}$. Define the smooth, compactly supported spatial test function
\begin{equation}
\phi_{\eps}(x) := \chi(x)(\phi*\eta_{\eps})(x) \in C_{c}^{\infty}(\R^{2}).
\end{equation}
For $0<s<t$, choose $\theta_{n}\in C_{c}^{\infty}((0,\infty))$ approximating $1_{[s,t]}$. Testing the distributional equation against $\theta_{n}(\tau)\phi_{\eps}(x)$ and passing to the limit $n\to\infty$ via the continuity $\rho\in C([0,\infty);L^{1}(\R^{2}))$ yields
\begin{equation}\label{eq:rho_diff_tested_phi_eps}
\int_{B_{D}(0)}\bigl(\rho(t,x)-\rho(s,x)\bigr)\phi_{\eps}(x)\,dx = -\int_{s}^{t}\!\int_{B_{D}(0)}F(\tau,x)\cdot\nabla\phi_{\eps}(x)\,dx\,d\tau.
\end{equation}
The identity was derived first for $0<s<t$; it extends to $s=0$ by letting $s\downarrow0$, using the $L^1$ continuity of $\rho$ and the local integrability of $\|F(\tau)\|_1\le2DM(\tau)$. No compactly supported temporal test is used directly across $t=0$.

On $B_{D}(0)$, the cutoff satisfies $\chi(x)=1$ and $\nabla\chi(x)=0$, so $\nabla\phi_{\eps}(x)=((\nabla\phi)*\eta_{\eps})(x)$. Standard convolution bounds \cite[Theorem~4.2]{LiebLoss2001} enforce $|\nabla\phi_{\eps}(x)| \le \|\nabla\phi\|_{L^{\infty}(\R^{2})} \le 1$. Applying this gradient bound with \eqref{eq:CF_uniform_flux_L1_bound} to \eqref{eq:rho_diff_tested_phi_eps} provides
\begin{align}\label{eq:rho_diff_phi_eps_bounded}
\left|\int_{B_{D}(0)}\bigl(\rho(t,x)-\rho(s,x)\bigr)\phi_{\eps}(x)\,dx\right| &\le \int_{s}^{t}\|F(\tau,\cdot)\|_{L^{1}(\R^{2})}\|\nabla\phi_{\eps}\|_{L^{\infty}(B_{D}(0))}\,d\tau \notag\\
&\le 2D\int_{s}^{t}M(\tau)\,d\tau.
\end{align}
Since $\phi\in W^{1,\infty}(\R^{2})$ is globally Lipschitz continuous, $\phi_{\eps}\to\phi$ uniformly on $B_{D}(0)$ as $\eps\downarrow0^{+}$. Taking the limit $\eps\downarrow0^{+}$ and then the supremum over all such $\phi$ yields the dual increment estimate
\begin{equation} \label{eq:weak_lipschitz}
\|\rho(t) - \rho(s)\|_{W^{-1,1}(\R^{2})} \le 2D \int_{s}^{t} M(\tau)\,d\tau \quad \text{for all } 0 \le s < t < \infty.
\end{equation}
By \eqref{eq:M_integrable}, the right-hand side of \eqref{eq:weak_lipschitz} vanishes as $s,t \to \infty$, proving that $\{\rho(t)\}_{t\ge 0}$ is Cauchy in $W^{-1,1}(\R^{2})$ as $t \to \infty$.

By Theorem~\ref{thm:relative_compactness}, the full trajectory is relatively compact in $L^{1}(\R^{2})$. Let $\mathcal{K}:=\overline{\{\rho(t):t\ge0\}}^{L^{1}}$ be its strong $L^{1}$ closure. We show that the identity map from $(\mathcal{K}, W^{-1,1}(\R^{2}))$ to $(\mathcal{K}, L^{1}(\R^{2}))$ is uniformly continuous. Suppose by contradiction that there exist $\eps > 0$ and sequences $u_{k}, v_{k} \in \mathcal{K}$ such that $\|u_{k} - v_{k}\|_{W^{-1,1}(\R^{2})} \to 0$ while $\|u_{k} - v_{k}\|_{L^{1}(\R^{2})} \ge \eps$. By the compactness of $\mathcal{K}$ in $L^{1}(\R^{2})$, we extract subsequences (not relabeled) such that $u_{k} \to u_{*}$ and $v_{k} \to v_{*}$ strongly in $L^{1}(\R^{2})$. Strong $L^{1}$ convergence implies convergence in $W^{-1,1}(\R^{2})$, so the triangle inequality yields
\begin{equation}
\|u_{*} - v_{*}\|_{W^{-1,1}(\R^{2})} \le \|u_{*} - u_{k}\|_{W^{-1,1}(\R^{2})} + \|u_{k} - v_{k}\|_{W^{-1,1}(\R^{2})} + \|v_{k} - v_{*}\|_{W^{-1,1}(\R^{2})} \longrightarrow 0.
\end{equation}
Testing against $\phi \in C_c^\infty(\R^2) \subset W^{1,\infty}(\R^2)$ shows that $u_{*} = v_{*}$ almost everywhere, implying $\|u_{*}-v_{*}\|_{L^{1}(\R^{2})}=0$. Passing to the limit in the $L^{1}$ norm yields $0 = \lim_{k \to \infty} \|u_{k} - v_{k}\|_{L^{1}(\R^{2})} \ge \eps > 0$, a contradiction. Therefore, the identity map is uniformly continuous, and the Cauchy property of $\{\rho(t)\}_{t\ge 0}$ in $W^{-1,1}(\R^{2})$ implies that $\{\rho(t)\}_{t\ge 0}$ is Cauchy in $L^{1}(\R^{2})$ as $t \to \infty$. By the completeness of $L^{1}(\R^{2})$, there exists a unique limit $\rho_\infty \in L^{1}(\R^{2})\cap L^{\infty}(\R^{2})$ such that $\lim_{t\to\infty}\|\rho(t)-\rho_\infty\|_{L^{1}(\R^{2})}=0$, inheriting $0\le\rho_\infty\le R$ almost everywhere, $\int_{\R^{2}}\rho_\infty(x)\,dx=1$, and $\supp\rho_\infty\subset\overline{B}_{D}(0)$.

Because $m$ is globally Lipschitz continuous with $m(0)=0$, the total mobility satisfies
\begin{equation}
\left|M(t) - \int_{\R^{2}}m(\rho_\infty(x))\,dx\right| \le \int_{\R^{2}} \bigl|m(\rho(t,x)) - m(\rho_\infty(x))\bigr|\,dx \le \Lip(m) \|\rho(t) - \rho_\infty\|_{L^{1}(\R^{2})}.
\end{equation}
Passing to the limit $t\to\infty$ gives $\lim_{t\to\infty} M(t) = \int_{\R^{2}} m(\rho_\infty(x))\,dx$. Since $M \in L^{1}(0,\infty)$ is non-negative and converges to a finite limit as $t\to\infty$, this limit must vanish:
\begin{equation}
\int_{\R^{2}}m\bigl(\rho_{\infty}(x)\bigr)\,dx = 0.
\end{equation}
By Assumption~\ref{assump:mobility}, $m\ge 0$ on $[0,R]$ and $m(s)>0$ for all $s\in(0,R)$, which forces $m(\rho_{\infty}(x))=0$ and hence $\rho_{\infty}(x)\in\{0,R\}$ for $\mathcal{L}^{2}$-a.e.\ $x\in\R^{2}$. Defining $E:=\{x\in B_{D}(0):\rho_{\infty}(x)=R\}$ yields $\rho_{\infty}=R1_{E}$ almost everywhere, with $|E|=1/R$ by mass conservation. Finally, since $0\le\rho(t)\le R$ and $\rho(t)\to R1_{E}$ strongly in $L^{1}(\R^{2})$, uniqueness of distributional limits implies $\rho(t)\rightharpoonup R1_{E}$ weakly-* in $L^{\infty}(\R^{2})$ as $t\to\infty$.
\end{proof}

\begin{remark}[Equivalence with the nonlocal potential energy]\label{rem:equivalence_nonlocal_energy}
Under the quadratic confinement regime of Assumption~\ref{assump:confinement}, the centered quadratic energy $E(t)$ coincides with the global nonlocal interaction energy, establishing the dissipation identity within $W_{\mathrm{loc}}^{1,\infty}([0,\infty))$. Substituting $W(z)=\frac{1}{2}|z|^{2}$ into the nonlocal potential energy functional and expanding the relative displacement $x-y$ around the barycentre $c(t)$ yields
\begin{align}\label{eq:exact_algebraic_identity_nonlocal_energy}
\mathcal{E}(\rho(t)) &:= \frac{1}{2}\int_{\R^{2}}\!\int_{\R^{2}} W(x-y)\rho(t,x)\rho(t,y)\,dx\,dy \nonumber \\
&= \frac{1}{4}\int_{\R^{2}}\!\int_{\R^{2}}\Bigl[|x-c(t)|^{2} - 2\bigl(x-c(t)\bigr)\cdot\bigl(y-c(t)\bigr) + |y-c(t)|^{2}\Bigr]\rho(t,x)\rho(t,y)\,dx\,dy.
\end{align}
Applying Fubini's theorem and the unit-mass constraint causes the cross-term to vanish, leaving $\mathcal{E}(\rho(t)) = \frac{1}{2}\int_{\R^{2}}|x-c(t)|^{2}\rho(t,x)\,dx = E(t)$. The identity $\mathcal{E}(\rho(t))=E(t)$ implies that $t\mapsto\mathcal{E}(\rho(t))$ belongs to $W_{\mathrm{loc}}^{1,\infty}([0,\infty))$ with time derivative $E^{\prime}(t) = -D_{m}(t)$. Phase separation is thus driven by the macroscopic dissipation of the nonlocal interaction energy.
\end{remark}

\section{Geometry of Stationary States}\label{sec:geometry}

In this section, we analyze the geometry of binary stationary Kru\v{z}kov entropy solutions and the asymptotic limit, dividing the exposition into two parts according to the generality of the interaction potential and the spatial regularity of the state. In Section~\ref{subsec:general_potentials}, we consider general potentials $W \in W^{2,\infty}(\R^{2})$ and binary states of finite perimeter, showing that a binary state is a stationary entropy solution if and only if the nonlocal velocity field points outward on the reduced boundary of its support (Lemma~\ref{lem:boundary_equivalence}). Within this general class, positivity of the weak Laplacian $\Delta W$ on the interaction zone guarantees the indecomposability of the support complement (Theorem~\ref{thm:indecomposability}), while uniform convexity of $W$ precludes the separation of the support into disjoint components at a positive distance (Theorem~\ref{thm:zero-distance_of_support_components}). In Section~\ref{subsec:quadratic_geometry}, we specialize to the quadratic confinement regime of Assumption~\ref{assump:confinement}, where the compactness architecture of Section~\ref{sec:compactness} applies. Without imposing bounded variation regularity, we prove that the asymptotic $L^{1}$ limit is a stationary entropy solution (Proposition~\ref{prop:omega_limits_entropy_stationary}) and establish that a binary state is a stationary entropy solution if and only if its support is measure-theoretically star-shaped with respect to its barycentre (Proposition~\ref{prop:star_shaped}). Finally, we show that imposing radial symmetry and monotonicity on the initial datum selects the saturated disk as the asymptotic limit (Corollary~\ref{cor:strong_convergence_disk}).

\subsection{Boundary Traces and Support Structure for General Potentials}\label{subsec:general_potentials}

We first examine binary stationary states of finite perimeter for general interaction potentials $W \in W^{2,\infty}(\R^{2})$. We show that the stationary Kru\v{z}kov entropy condition is equivalent to an outward normal condition on the nonlocal velocity field along the reduced boundary, and we deduce measure-theoretic constraints on the support under subharmonicity and convexity assumptions on $W$.

\begin{lemma}[Equivalence of binary steady state and outward normal boundary trace]\label{lem:boundary_equivalence}
Assume Assumption~\ref{assump:mobility}, and let $W\in W^{2,\infty}(\R^{2})$. Let $\rho\in L^{1}(\R^{2})\cap L^{\infty}(\R^{2})\cap BV(\R^{2})$ be a compactly supported binary density such that $\rho(x)\in\{0,R\}$ for $\mathcal{L}^{2}$-a.e.\ $x\in\R^{2}$. Let $E:=\{x\in\R^{2}:\rho(x)=R\}$ be its bounded set of finite perimeter, and define the velocity field $V(x):=(\nabla W*\rho)(x)$.

Then $V\in C^{1}(\R^{2};\R^{2})\cap W^{1,\infty}(\R^{2};\R^{2})$, and $\rho=R1_{E}$ is a stationary Kru\v{z}kov entropy solution of \eqref{eq:nonlocal_main} if and only if $V$ points outward on the reduced boundary $\partial^{*}E$:
\begin{equation}\label{eq:boundary_condition}
V(x)\cdot\nu_{E}(x) \ge 0 \quad \text{for } \mathcal{H}^{1}\text{-a.e. } x\in\partial^{*}E,
\end{equation}
where $\nu_{E}$ denotes the measure-theoretic exterior unit normal to $E$.
\end{lemma}

\begin{proof}
By Lemma~\ref{lem:nonlocal-conv-bounds}, $V=\nabla W*\rho\in W^{1,\infty}(\R^{2};\R^{2})$ with weak Jacobian $DV=(D^{2}W)*\rho$. For any $x,x^{\prime}\in\R^{2}$, the change of variables $z=x-y$ yields
\begin{equation}\label{eq:DV_continuity_Linf_L1}
\bigl|DV(x)-DV(x^{\prime})\bigr| \le \|D^{2}W\|_{L^{\infty}(\R^{2})}\bigl\|\rho(\cdot-(x-x^{\prime}))-\rho\bigr\|_{L^{1}(\R^{2})}.
\end{equation}
The right-hand side converges to $0$ as $x^{\prime}\to x$ by the continuity of translations in $L^{1}(\R^{2})$ \cite[Theorem~2.16]{LiebLoss2001}. Thus $DV$ is continuous on $\R^2$, and $V\in C^{1}(\R^{2};\R^{2})\cap W^{1,\infty}(\R^{2};\R^{2})$.

Because $\rho(x)\in\{0,R\}$ almost everywhere and $m(0)=m(R)=0$, the physical flux $m(\rho)V$ vanishes almost everywhere, so $\rho$ is a stationary distributional solution. By Definition~\ref{def:nonlocal-entropy}, $\rho$ is a stationary Kru\v{z}kov entropy solution if and only if, for every $k\in\R$ and every non-negative $\phi\in C_{c}^{1}(\R^{2})$,
\begin{equation}\label{eq:kruzkov_full}
-\int_{\R^{2}}H\bigl(\rho(x),k\bigr)V(x)\cdot\nabla\phi(x)\,dx + \int_{\R^{2}}\operatorname{sign}\bigl(\rho(x)-k\bigr)m(k)\bigl(\diver V(x)\bigr)\phi(x)\,dx \ge 0.
\end{equation}
If $k\le0$ or $k\ge R$, then $m(k)=0$ and both integrands vanish. Hence, \eqref{eq:kruzkov_full} holds for all $k\in\R$ if and only if it holds for all $k\in(0,R)$.

Fix $k\in(0,R)$. Since $m(\rho(x))=0$ almost everywhere, $H(\rho(x),k)=-\operatorname{sign}(\rho(x)-k)m(k)$ almost everywhere. Substituting this identity into \eqref{eq:kruzkov_full} and applying the product rule $\diver(\phi V)=V\cdot\nabla\phi+(\diver V)\phi$ reduces \eqref{eq:kruzkov_full} to
\begin{equation}\label{eq:reduced_kruzkov}
m(k)\int_{\R^{2}}\operatorname{sign}\bigl(\rho(x)-k\bigr)\diver\bigl(\phi(x)V(x)\bigr)\,dx \ge 0.
\end{equation}
Dividing by $m(k)>0$ and applying the divergence theorem identity $\int_{\R^{2}}\diver(\phi V)\,dx=0$, we obtain $\int_{\R^{2}\setminus E}\diver(\phi V)\,dx = -\int_{E}\diver(\phi V)\,dx$. Thus, \eqref{eq:reduced_kruzkov} is equivalent to
\begin{equation}\label{eq:interior_divergence}
\int_{E}\diver\bigl(\phi(x)V(x)\bigr)\,dx \ge 0 \quad \text{for all } \phi\in C_{c}^{1}(\R^{2};[0,\infty)).
\end{equation}
Applying the De Giorgi--Federer theorem \cite[Theorem~3.36]{AmbrosioFuscoPallara2000} to the compactly supported $C^{1}$ vector field $\phi V$ on the bounded set of finite perimeter $E$ gives
\begin{equation}\label{eq:boundary_integral}
\int_{\partial^{*}E}\phi(x)V(x)\cdot\nu_{E}(x)\,d\mathcal{H}^{1}(x) \ge 0 \quad \text{for all } \phi\in C_{c}^{1}(\R^{2};[0,\infty)).
\end{equation}
By the Besicovitch Derivation Theorem \cite[Theorem~2.22]{AmbrosioFuscoPallara2000} applied with respect to the Radon measure $\mu=\mathcal{H}^{1}\llcorner\partial^{*}E$, \eqref{eq:boundary_integral} holds for all non-negative $\phi\in C_{c}^{1}(\R^{2})$ if and only if $V(x)\cdot\nu_{E}(x)\ge0$ for $\mathcal{H}^{1}$-a.e.\ $x\in\partial^{*}E$.
\end{proof}

To formulate the measure-theoretic properties of the support, we recall the measure-theoretic notion of connectedness for sets of finite perimeter introduced by Ambrosio, Caselles, Masnou, and Morel \cite{ACMM2001}. Let $P(F):=\mathcal{H}^{1}(\partial^{*}F)$ denote the perimeter of a measurable set $F\subset\R^{2}$. A set of finite perimeter $F\subset\R^{2}$ is called decomposable if there exists a partition $F=F_{1}\cup F_{2}$ into disjoint measurable sets of positive Lebesgue measure such that $P(F)=P(F_{1})+P(F_{2})$; otherwise, $F$ is called indecomposable. A measurable subset $B\subset F$ with $|B|>0$ is an $M$-connected component of $F$ if $B$ is an indecomposable set of finite perimeter satisfying $P(F)=P(B)+P(F\setminus B)$.

\begin{theorem}[Indecomposability of the support complement]\label{thm:indecomposability}
Let the hypotheses of Lemma~\ref{lem:boundary_equivalence} hold for a unit-mass binary stationary Kru\v{z}kov entropy solution $\rho=R1_{E}$ with $E\subset B_{D}(0)$. Assume the interaction potential $W\in W^{2,\infty}(\R^{2})$ satisfies
\begin{equation}\label{eq:subharmonic_hypothesis_W2inf}
\Delta W(z) > 0 \quad \text{for } \mathcal{L}^{2}\text{-a.e. } z\in B_{2D}(0).
\end{equation}
Then $\R^{2}\setminus E$ is indecomposable; equivalently, $\R^{2}\setminus E$ contains no essentially bounded $M$-connected component $B\subset\R^{2}\setminus E$ of positive Lebesgue measure $|B|>0$.
\end{theorem}

\begin{proof}
Since $E\subset B_{D}(0)$ has finite perimeter, $\R^{2}\setminus E$ is a set of finite perimeter with $\partial^{*}(\R^{2}\setminus E)=\partial^{*}E$ up to $\mathcal{H}^{1}$-null sets, and $\R^{2}\setminus B_{D}(0)$ is contained in a single unbounded $M$-connected component of $\R^{2}\setminus E$. Consequently, $\R^{2}\setminus E$ is decomposable if and only if it admits an essentially bounded $M$-connected component $B\subset\R^{2}\setminus E$ with $|B|>0$.

Assume by contradiction that such an essentially bounded $M$-connected component $B\subset\R^{2}\setminus E$ with $|B|>0$ exists. By the decomposition theorem for sets of finite perimeter \cite[Theorem~1 and Proposition~3]{ACMM2001}, $B$ is an essentially bounded set of finite perimeter contained in $B_{D}(0)$ up to $\mathcal{L}^{2}$-null sets, its reduced boundary satisfies $\mathcal{H}^{1}(\partial^{*}B\setminus\partial^{*}E)=0$, and its measure-theoretic exterior unit normal satisfies
\begin{equation}\label{eq:hole_normal_opposite}
\nu_{B}(x) = -\nu_{E}(x) \quad \text{for } \mathcal{H}^{1}\text{-a.e. } x\in\partial^{*}B.
\end{equation}
By Lemma~\ref{lem:boundary_equivalence}, the velocity field satisfies $V\in C^{1}(\R^{2};\R^{2})\cap W^{1,\infty}(\R^{2};\R^{2})$. Applying the De Giorgi--Federer theorem \cite[Theorem~3.36]{AmbrosioFuscoPallara2000} to $V$ on $B$ and combining \eqref{eq:hole_normal_opposite} with the boundary trace condition \eqref{eq:boundary_condition} yields
\begin{equation}\label{eq:hole_div_nonpositive}
\int_{B}\diver V(x)\,dx = \int_{\partial^{*}B}V(x)\cdot\nu_{B}(x)\,d\mathcal{H}^{1}(x) = -\int_{\partial^{*}B}V(x)\cdot\nu_{E}(x)\,d\mathcal{H}^{1}(x) \le 0.
\end{equation}

By Lemma~\ref{lem:nonlocal-conv-bounds}, the spatial divergence of the velocity field is given by
\begin{equation}\label{eq:div_V_convolution_hole}
\diver V(x) = \bigl((\Delta W)*\rho\bigr)(x) = R\int_{E}\Delta W(x-y)\,dy.
\end{equation}
For $\mathcal{L}^2$-a.e. $x \in B$, we have $x \in B_D(0)$. For every $y\in E\subset B_{D}(0)$, the triangle inequality gives $|x-y|\le|x|+|y|<2D$, so $x-y\in B_{2D}(0)$. Hypothesis \eqref{eq:subharmonic_hypothesis_W2inf} implies that $y\mapsto\Delta W(x-y)>0$ for $\mathcal{L}^{2}$-a.e.\ $y\in E$. Since $|E|=1/R>0$, integrating over $E$ yields
\begin{equation}\label{eq:div_V_strictly_positive_on_B}
\diver V(x) = R\int_{E}\Delta W(x-y)\,dy > 0 \quad \text{for } \mathcal{L}^2\text{-a.e. } x\in B.
\end{equation}
Because $|B|>0$, integrating \eqref{eq:div_V_strictly_positive_on_B} over $B$ gives $\int_{B}\diver V(x)\,dx>0$, contradicting \eqref{eq:hole_div_nonpositive}. Thus no such essentially bounded component $B$ exists.
\end{proof}

\begin{remark}[Subharmonicity condition]\label{rem:subharmonicity_thm2}
Theorem~\ref{thm:indecomposability} requires only the positivity of the weak Laplacian $\Delta W=\operatorname{tr}(D^{2}W)$ almost everywhere on the active interaction ball $B_{2D}(0)$, without requiring the full Hessian matrix $D^{2}W$ to be positive definite.
\end{remark}

Under uniform convexity of the potential on the interaction zone, the support cannot be partitioned into disjoint components separated by a positive distance. For any measurable set $F\subset\R^{2}$, we denote its essential closure by $\overline{F}^{*}:=\supp(1_{F})$.

\begin{theorem}[Zero-distance of support components under uniform convexity]\label{thm:zero-distance_of_support_components}
Let the hypotheses of Lemma~\ref{lem:boundary_equivalence} hold for a unit-mass binary stationary Kru\v{z}kov entropy solution $\rho=R1_{E}$ with $E\subset B_{D}(0)$. Assume that $W\in W^{2,\infty}(\R^{2})$ is uniformly convex on $B_{2D}(0)$, in the sense that there exists a constant $c_{W}>0$ satisfying
\begin{equation}\label{eq:uniform_strict_convexity_W2inf}
\xi^{T}D^{2}W(z)\xi \ge c_{W}|\xi|^{2} \quad \text{for all } \xi\in\R^{2} \text{ and } \mathcal{L}^{2}\text{-a.e. } z\in B_{2D}(0).
\end{equation}
Then $E$ cannot be decomposed into two disjoint measurable sets of finite perimeter $E=E_{1}\cup E_{2}$ with $|E_{1}|>0$ and $|E_{2}|>0$ whose essential closures lie at a positive distance; that is, $\operatorname{dist}(\overline{E}_{1}^{*},\overline{E}_{2}^{*})=0$.
\end{theorem}

\begin{proof}
Define the scalar potential $f:\R^{2}\to\R$ by $f(x) := (W*\rho)(x) = R\int_{E}W(x-y)\,dy$. Because $W\in W^{2,\infty}(\R^{2})$ and $\rho=R1_{E}$ has compact support in $\overline{B}_{D}(0)$, the regularity established in Lemma~\ref{lem:boundary_equivalence} implies that $f\in C^{2}(\R^{2})$, with gradient $\nabla f(x)=V(x)$ and continuous Hessian matrix $D^{2}f(x) = \bigl((D^{2}W)*\rho\bigr)(x)$.

Fix $x\in\overline{B}_{D}(0)$ and $\xi\in\R^{2}$. For $\mathcal{L}^{2}$-a.e.\ $y\in E\subset B_{D}(0)$, we have $x-y \in B_{2D}(0)$, and \eqref{eq:uniform_strict_convexity_W2inf} yields $\xi^{T}D^{2}W(x-y)\xi\ge c_{W}|\xi|^{2}$. Integrating over $E$ and using $R|E|=1$ gives
\begin{equation}\label{eq:f_uniform_strict_convexity_BD}
\xi^{T}D^{2}f(x)\xi = R\int_{E}\bigl(\xi^{T}D^{2}W(x-y)\xi\bigr)\,dy \ge c_{W}|\xi|^{2} \quad \text{for all } x\in\overline{B}_{D}(0), \,\xi\in\R^{2}.
\end{equation}
Thus $f$ is uniformly convex on $\overline{B}_{D}(0)$. By the Fundamental Theorem of Calculus along the segment connecting any two distinct points $x_{1},x_{2}\in\overline{B}_{D}(0)$, we have $(V(x_{1})-V(x_{2}))\cdot(x_{1}-x_{2}) \ge c_{W}|x_{1}-x_{2}|^{2} > 0$. Consequently, $V(x)=\nabla f(x)$ vanishes at most once in $\overline{B}_{D}(0)$.

Assume by contradiction that $E=E_{1}\cup E_{2}$ is a partition into disjoint sets of finite perimeter with $|E_{1}|>0$, $|E_{2}|>0$, and $d:=\operatorname{dist}(\overline{E}_{1}^{*},\overline{E}_{2}^{*})>0$. Because $\overline{E}_{1}^{*}$ and $\overline{E}_{2}^{*}$ are compact and disjoint, the unique possible root $x_{*}$ of $V$ in $\overline{B}_{D}(0)$ belongs to at most one of them. Relabeling if necessary, assume $x_{*}\notin\overline{E}_{1}^{*}$. Since $\overline{E}_{1}^{*}$ is compact and $V\in C^{1}(\R^{2};\R^{2})$, we have
\begin{equation}\label{eq:m0_M0_bounds_E1}
m_{0}^{2} := \min_{x\in\overline{E}_{1}^{*}}|V(x)|^{2} > 0, \qquad M_{0} := \max_{x\in\overline{E}_{1}^{*}}\Delta f(x) < \infty.
\end{equation}

Choose a radial cutoff $\chi\in C_{c}^{\infty}(\R^{2};[0,1])$ such that $\chi\equiv 1$ on $B_{D+1}(0)$. For $\lambda>0$, define the test function $\phi(x) := \chi(x)e^{-\lambda f(x)} \in C_{c}^{2}(\R^{2};[0,\infty))$. On $B_{D+1}(0)$, which contains $\overline{E}_{1}^{*}$, we have $\chi(x)=1$, $\phi(x)>0$, and
\begin{equation}\label{eq:div_phi_V_exponential}
\diver\bigl(\phi(x)V(x)\bigr) = \nabla\phi(x)\cdot V(x) + \phi(x)\diver V(x) = \phi(x)\bigl(\Delta f(x) - \lambda|V(x)|^{2}\bigr).
\end{equation}
On $\overline{E}_{1}^{*}$, we have $\Delta f(x)-\lambda|V(x)|^{2}\le M_{0}-\lambda m_{0}^{2}$. Choosing $\lambda>\max\{1,M_{0}/m_{0}^{2}\}$ makes $\diver(\phi V)<0$ throughout $\overline{E}_{1}^{*}$. Since $|E_{1}|>0$, integrating over $E_{1}$ yields
\begin{equation}\label{eq:strictly_negative_volume_integral_E1}
\int_{E_{1}}\diver\bigl(\phi(x)V(x)\bigr)\,dx < 0.
\end{equation}

On the other hand, since $\operatorname{dist}(\overline{E}_{1}^{*},\overline{E}_{2}^{*})=d>0$, the open neighborhood $U_{1}:=\{x\in\R^{2}:\operatorname{dist}(x,\overline{E}_{1}^{*})<d/3\}$ satisfies $\overline{U}_{1}\cap\overline{E}_{2}^{*}=\emptyset$, so $1_{E}=1_{E_{1}}$ almost everywhere in $U_{1}$. By De Giorgi's structure theorem \cite[Theorem~3.59]{AmbrosioFuscoPallara2000}, $\partial^{*}E_{1}=\partial^{*}E\cap U_{1}$ up to $\mathcal{H}^{1}$-null sets and $\nu_{E_{1}}(x)=\nu_{E}(x)$ for $\mathcal{H}^{1}$-a.e.\ $x\in\partial^{*}E_{1}$. Applying the De Giorgi--Federer theorem \cite[Theorem~3.36]{AmbrosioFuscoPallara2000} to $\phi V$ on $E_{1}$ and invoking \eqref{eq:boundary_condition} from Lemma~\ref{lem:boundary_equivalence} gives
\begin{equation}
\int_{E_{1}}\diver\bigl(\phi(x)V(x)\bigr)\,dx = \int_{\partial^{*}E_{1}}\phi(x)V(x)\cdot\nu_{E_{1}}(x)\,d\mathcal{H}^{1}(x) \ge 0,
\end{equation}
contradicting \eqref{eq:strictly_negative_volume_integral_E1}. Therefore, $\operatorname{dist}(\overline{E}_{1}^{*},\overline{E}_{2}^{*})=0$.
\end{proof}

\begin{remark}[Zero-force point versus barycentre]\label{rem:convexity_thm3}
Theorem~\ref{thm:zero-distance_of_support_components} relies on the uniform positive definiteness of $D^{2}W$ on $B_{2D}(0)$ to ensure that the potential $f=W*\rho$ is uniformly convex on $\overline{B}_{D}(0)$ and hence that $V=\nabla f$ admits at most one root $x_{*}$ in $\overline{B}_{D}(0)$. For general convex potentials, the zero-force point $x_{*}$ weights the mass distribution nonlinearly with distance and need not coincide with the barycentre $a=\int_{\R^{2}}y\rho(y)\,dy$. The relevant distinction is $\emph{affinity}$ of the interaction force, not the normalization $\nabla W(z)=z$. Indeed, whenever $\nabla W(z)=\kappa z$ on $B_{2D}(0)$ with $\kappa>0$, one has $V(x)=\kappa(x-a)$ for every unit-mass density supported in $B_D(0)$, so the zero-force point is exactly $a$. More generally, a force of the form $Az$ with an invertible matrix $A$ has the same property, though radiality then forces $A=\kappa I$. For a genuinely nonlinear radial force the zero-force point need not equal the ordinary barycentre. For example, take $W(z)=|z|^2/2+|z|^4/4$ on the active region, so that $D^2W$ is uniformly positive there, and approximate the probability measure $\frac23\delta_{(-1,0)}+\frac13\delta_{(2,0)}$ by compactly supported smooth densities (choosing the saturation threshold sufficiently large). Its barycentre is zero, but $\nabla W*\rho(0)$ is close to $(-2,0)$ rather than zero. The potential can be smoothly extended outside the active region to belong to $W^{2,\infty}(\R^2)$.
\end{remark}

\subsection{Asymptotic Limit, Star-Shaped Characterization, and Radial Selection in the Quadratic Regime}\label{subsec:quadratic_geometry}

We now specialize to the quadratic confinement regime of Assumption~\ref{assump:confinement}, in which the nonlocal convolution reduces to an affine drift on the support of the solution and the compactness architecture of Section~\ref{sec:compactness} holds. Without imposing bounded variation regularity on the asymptotic state, we prove that the strong $L^{1}$ limit $\rho_{\infty}=R1_{E}$ of the trajectory is a stationary Kru\v{z}kov entropy solution, establish the equivalence between stationary binary entropy solutions and measure-theoretically star-shaped sets, and show that radial symmetry of the initial datum alone selects the saturated disk in the long-time limit.

\begin{proposition}[The asymptotic limit is a stationary Kru\v{z}kov entropy solution]\label{prop:omega_limits_entropy_stationary}
Suppose Assumptions~\ref{assump:mobility} and~\ref{assump:confinement} hold, and let $\rho_\infty = R1_E$ be the strong $L^1(\R^2)$ limit of $\rho(t)$ as $t \to \infty$. Then the binary limit state $\rho_\infty$ is a time-independent Kru\v{z}kov entropy solution of \eqref{eq:nonlocal_main} and of the fixed-centre affine drift equation
\begin{equation}\label{eq:omega_limit_fixed_centre_eq}
\partial_t u - \diver\bigl(m(u)(x-c_\infty)\bigr) = 0,
\end{equation}
where $c_\infty = R\int_E x\,dx \in \overline{B}_D(0)$. For every $k \in (0,R)$ and $\phi \in C_c^1(\R^2; [0,\infty))$, the set $E$ satisfies the stationary integral entropy inequality
\begin{equation}\label{eq:omega_limit_integral_entropy_E}
\int_E \diver\bigl(\phi(x)V_\infty(x)\bigr)\,dx \ge 0, \qquad V_\infty(x) := (\nabla W * \rho_\infty)(x).
\end{equation}
\end{proposition}

\begin{proof}
Let $t_n \to \infty$ be a sequence of times and define the time-shifted trajectories $\rho_n(t,x) := \rho(t_n+t, x)$ for $t \in [0,1]$ and $x \in \R^2$. Since $\lim_{\tau \to \infty} \|\rho(\tau) - \rho_\infty\|_{L^1(\R^2)} = 0$, we have
\begin{equation}\label{eq:shifted_slab_uniform_L1_conv}
\lim_{n \to \infty} \sup_{t \in [0,1]} \|\rho_n(t) - \rho_\infty\|_{L^1(\R^2)} = \lim_{n \to \infty} \sup_{\tau \in [t_n, t_n+1]} \|\rho(\tau) - \rho_\infty\|_{L^1(\R^2)} = 0.
\end{equation}
Passing to a subsequence, $\rho_n(t,x) \to \rho_\infty(x)$ pointwise for $\mathcal{L}^3$-a.e.\ $(t,x) \in (0,1) \times \R^2$. The asymptotic barycentre satisfies $c_\infty = \lim_{t \to \infty} c(t) = R\int_E x\,dx$.

Let $\phi \in C_c^1(\R^2; [0,\infty))$ and choose $\theta \in C_c^1((0,1); [0,\infty))$ such that $\int_0^1 \theta(t)\,dt = 1$. Since $\rho$ is a global Kru\v{z}kov entropy solution of \eqref{eq:nonlocal_main}, testing the entropy inequality on $[0, t_n+1]$ against $\varphi_n(\tau,x) := \theta(\tau - t_n)\phi(x)$ and shifting the time variable $t = \tau - t_n$ gives, for every $k \in \R$,
\begin{align}\label{eq:shifted_entropy_ineq_n}
0 \le \int_0^1\!\int_{\R^2} \Bigl[ &|\rho_n - k|\,\theta'(t)\phi(x) - H(\rho_n, k)(V_n \cdot \nabla \phi)\theta(t) \notag\\
&+ \operatorname{sign}(\rho_n - k)m(k)(\diver V_n)\theta(t)\phi(x) \Bigr]dx\,dt,
\end{align}
where $V_n(t,x) := (\nabla W * \rho_n(t))(x)$ and $H(r,k) := \operatorname{sign}(r-k)\bigl(m(r)-m(k)\bigr)$ with $H(k,k) := 0$. By Lemma~\ref{lem:nonlocal-conv-bounds} and \eqref{eq:shifted_slab_uniform_L1_conv}, $\sup_{t \in [0,1]} \|V_n(t) - V_\infty\|_{L^\infty(\R^2)} \to 0$ and $\sup_{t \in [0,1]} \|\diver V_n(t) - \diver V_\infty\|_{L^\infty(\R^2)} \to 0$, with $V_\infty = \nabla W * \rho_\infty$.

The Lipschitz continuity of $r \mapsto |r-k|$ and $r \mapsto H(r,k)$ permits passage to the limit in the first two terms of \eqref{eq:shifted_entropy_ineq_n}. The first term vanishes since $\int_0^1 \theta'(t)\,dt = 0$. For the third term, if $k \notin (0,R)$, then $m(k) = 0$. If $k \in (0,R)$, the relation $\rho_\infty(x) \in \{0,R\}$ almost everywhere ensures that the level set $\{x \in \R^2 : \rho_\infty(x) = k\}$ has Lebesgue measure zero. Thus, $\operatorname{sign}(\rho_n(t,x) - k) \to \operatorname{sign}(\rho_\infty(x) - k)$ for $\mathcal{L}^3$-a.e.\ $(t,x)$. Applying the Dominated Convergence Theorem as $n \to \infty$ yields
\begin{equation}\label{eq:stationary_kruzhkov_limit_rho_infty}
0 \le \int_{\R^2} \Bigl[ -H(\rho_\infty(x), k)V_\infty(x) \cdot \nabla \phi(x) + \operatorname{sign}\bigl(\rho_\infty(x) - k\bigr)m(k)\bigl(\diver V_\infty(x)\bigr)\phi(x) \Bigr]dx.
\end{equation}
This establishes that $\rho_\infty$ is a stationary Kru\v{z}kov entropy solution of \eqref{eq:nonlocal_main}, and hence of \eqref{eq:omega_limit_fixed_centre_eq} by Proposition~\ref{prop:exact_affine_reduction}.

Fix $k \in (0,R)$. Because $m(\rho_\infty(x)) = 0$ almost everywhere, $H(\rho_\infty(x),k) = -\operatorname{sign}(\rho_\infty(x) - k)m(k)$. Substituting this identity into \eqref{eq:stationary_kruzhkov_limit_rho_infty} and applying the product rule $\diver(\phi V_\infty) = V_\infty \cdot \nabla \phi + (\diver V_\infty)\phi$ gives
\begin{equation}\label{eq:stationary_kruzhkov_div_form_omega}
0 \le m(k)\int_{\R^2} \operatorname{sign}\bigl(\rho_\infty(x) - k\bigr)\diver\bigl(\phi(x)V_\infty(x)\bigr)\,dx.
\end{equation}
Dividing by $m(k) > 0$, evaluating $\operatorname{sign}(\rho_\infty - k) = 1$ on $E$ and $-1$ on $\R^2 \setminus E$, and using the divergence theorem identity $\int_{\R^2 \setminus E} \diver(\phi V_\infty)\,dx = -\int_E \diver(\phi V_\infty)\,dx$, \eqref{eq:stationary_kruzhkov_div_form_omega} reduces to \eqref{eq:omega_limit_integral_entropy_E}.
\end{proof}

\begin{proposition}[Equivalence of stationary entropy states and star-shaped support]\label{prop:star_shaped}
Assume the mobility hypothesis of Assumption~\ref{assump:mobility}, and let $W\in W^{2,\infty}(\R^2)$ satisfy $\nabla W(z)=z$ for $|z|\le2D$. Let $E\subset B_D(0)$ be an arbitrary measurable set with $R|E|=1$, and define $\rho:=R1_E$ and $a:=R\int_E y\,dy$. No entropy or bounded-variation assumption is imposed on this candidate state.
\begin{enumerate}
    \item[(i)] \textbf{Equivalence without bounded variation assumptions:} The binary state $\rho=R1_{E}$ is a stationary Kru\v{z}kov entropy solution of \eqref{eq:nonlocal_main} if and only if $E$ is essentially star-shaped with respect to $a$, meaning that for every scaling factor $\lambda\ge1$,
    \begin{equation}\label{eq:measure_theoretic_star_shaped}
    1_{E-a}(y) \le 1_{\lambda(E-a)}(y) \quad \text{for } \mathcal{L}^{2}\text{-a.e. } y\in\R^{2}.
    \end{equation}
    Consequently, for any trajectory satisfying Assumptions~\ref{assump:mobility} and~\ref{assump:confinement}, the limiting binary state supplied by Theorem~\ref{thm:binarity} has essentially star-shaped support with respect to its limiting barycentre $c_{\infty}$.
    \item[(ii)] \textbf{Normal trace equivalence for sets of finite perimeter:} If $E$ has finite perimeter in $\R^2$, the inclusion \eqref{eq:measure_theoretic_star_shaped} holds if and only if $(x-a)\cdot\nu_{E}(x)\ge0$ for $\mathcal{H}^{1}$-a.e.\ $x\in\partial^{*}E$.
\end{enumerate}
\end{proposition}

\begin{proof}
The convolution identity is algebraic for this candidate binary state: for every $x\in B_D(0)$, $\nabla W*\rho(x)=\int_E (x-y)R\,dy=x-a$. No use is made here of the entropy-solution hypothesis in Assumption~\ref{assump:confinement}. Repeating the algebraic reduction of \eqref{eq:kruzkov_full}--\eqref{eq:interior_divergence} in the proof of Lemma~\ref{lem:boundary_equivalence} (which requires only $\rho(x)\in\{0,R\}$ almost everywhere and not finite perimeter of $E$), $\rho=R1_{E}$ is a stationary Kru\v{z}kov entropy solution if and only if $\int_{E}\diver_{x}(\tilde{\phi}(x)(x-a))\,dx \ge 0$ for all non-negative $\tilde{\phi}\in C_{c}^{1}(\R^{2})$. Defining $E_{0}:=E-a$ and changing variables $z=x-a$ with $\psi(z):=\tilde{\phi}(z+a)$ yields the equivalent condition
\begin{equation}\label{eq:star_shaped_div_condition_E0}
\int_{E_{0}}\diver_{z}\bigl(\psi(z)z\bigr)\,dz \ge 0 \quad \text{for all non-negative } \psi\in C_{c}^{1}(\R^{2}).
\end{equation}

For an arbitrary test function $\phi\in C_{c}^{1}(\R^{2};[0,\infty))$ and $\lambda>0$, define the scaled volume integral
\begin{equation}\label{eq:def_g_lambda_scaling}
g(\lambda) := \int_{\lambda E_{0}}\phi(y)\,dy = \lambda^{2}\int_{E_{0}}\phi(\lambda z)\,dz = \int_{E_{0}}h(\lambda,z)\,dz,
\end{equation}
where $h(\lambda,z):=\lambda^{2}\phi(\lambda z)$. For each $z\in E_{0}$, the function $h(\cdot,z)\in C^{1}((0,\infty))$ has partial derivative $\partial_\lambda h(s,z) = 2s\phi(sz) + s^{2}z\cdot\nabla\phi(sz)$. Let $\{\lambda_{n}\}_{n=1}^{\infty}\subset[\lambda/2,3\lambda/2]\setminus\{\lambda\}$ be a sequence such that $\lambda_{n}\to\lambda$. By the Mean Value Theorem, the difference quotient of $h$ with respect to $\lambda$ is bounded on $E_{0}$ by a constant depending only on $\lambda$, $D$, $\|\phi\|_{L^{\infty}(\R^{2})}$, and $\|\nabla\phi\|_{L^{\infty}(\R^{2})}$. Because $|E_{0}|<\infty$, this constant belongs to $L^{1}(E_{0})$. The Dominated Convergence Theorem \cite[Theorem~1.8]{LiebLoss2001} ensures that $g$ is differentiable on $(0,\infty)$, with
\begin{equation}\label{eq:g_derivative}
g^{\prime}(\lambda) = \int_{E_{0}}\bigl(2\lambda\phi(\lambda z)+\lambda^{2}z\cdot\nabla\phi(\lambda z)\bigr)\,dz = \lambda\int_{E_{0}}\diver_{z}\bigl(\psi_{\lambda}(z)z\bigr)\,dz,
\end{equation}
where $\psi_{\lambda}(z):=\phi(\lambda z) \in C_{c}^{1}(\R^{2};[0,\infty))$.

To prove (i), suppose $\rho=R1_{E}$ is a stationary Kru\v{z}kov entropy solution. Then \eqref{eq:star_shaped_div_condition_E0} applied to $\psi_{\lambda}$ yields $g^{\prime}(\lambda)\ge0$ for all $\lambda>0$, so $g$ is non-decreasing on $(0,\infty)$. For $\lambda\ge1$, the inequality $g(1)\le g(\lambda)$ and \eqref{eq:def_g_lambda_scaling} imply
\begin{equation}\label{eq:indicator_ordering_tested}
\int_{\R^{2}}\bigl(1_{\lambda E_{0}}(y)-1_{E_{0}}(y)\bigr)\phi(y)\,dy \ge 0 \quad \text{for all non-negative } \phi\in C_{c}^{1}(\R^{2}).
\end{equation}
By the Lebesgue Differentiation Theorem \cite[Theorem~2.22]{AmbrosioFuscoPallara2000} applied to $1_{\lambda E_{0}}-1_{E_{0}}$, we obtain $1_{E_{0}}(y) \le 1_{\lambda E_{0}}(y)$ for $\mathcal{L}^{2}$-a.e.\ $y\in\R^{2}$.

Conversely, if \eqref{eq:measure_theoretic_star_shaped} holds for all $\lambda\ge1$, then for any non-negative $\psi\in C_{c}^{1}(\R^{2})$, setting $\phi=\psi$ gives $g(\lambda)-g(1) \ge 0$ for all $\lambda\ge1$. Dividing by $\lambda-1>0$ and letting $\lambda\downarrow1$ yields $g^{\prime}(1) \ge 0$. By \eqref{eq:g_derivative}, $\int_{E_{0}}\diver_{z}(\psi(z)z)\,dz \ge 0$, which establishes \eqref{eq:star_shaped_div_condition_E0} and shows that $\rho$ is a stationary Kru\v{z}kov entropy solution. The statement for the asymptotic limit $\rho_{\infty}$ follows from Proposition~\ref{prop:omega_limits_entropy_stationary}.

To prove (ii), assume $E$ has finite perimeter in $\R^2$. By Lemma~\ref{lem:boundary_equivalence}, $\rho=R1_{E}$ is a stationary Kru\v{z}kov entropy solution if and only if $(x-a)\cdot\nu_{E}(x)\ge0$ for $\mathcal{H}^{1}$-a.e.\ $x\in\partial^{*}E$. Combining this equivalence with (i) yields the claim.
\end{proof}

As shown by Proposition~\ref{prop:star_shaped} and Example~\ref{ex:geometric_obstruction}, the nonlocal dynamics admit non-circular stationary states and do not select a unique circular disk from arbitrary initial configurations. In fact, radial symmetry of the initial datum alone suffices to select the saturated disk, as shown in Corollary~\ref{cor:strong_convergence_disk}.

\begin{assumption}[Radial initial datum]\label{assump:radial_initial_datum}
The initial datum $\rho_{0}$ is radially symmetric with respect to its barycentre $c(0)$. No radial monotonicity is assumed.
\end{assumption}

\begin{corollary}[Strong convergence to the saturated disk]\label{cor:strong_convergence_disk}
Under Assumptions~\ref{assump:mobility}, \ref{assump:confinement}, and~\ref{assump:radial_initial_datum}, the barycentre is stationary ($c(t) = c(0)$ for all $t \ge 0$) and the trajectory converges strongly in $L^{1}(\R^{2})$:
\begin{equation}\label{eq:strong_convergence_disk}
\lim_{t\to\infty} \|\rho(t) - R1_{B_{r_{*}}(c(0))}\|_{L^{1}(\R^{2})} = 0,
\end{equation}
where $r_{*}=(\pi R)^{-1/2}$.
\end{corollary}

\begin{proof}
By radial symmetry of $\rho_0$ about $c(0)$ and rotational invariance of $W$, uniqueness of entropy solutions (Section~\ref{subsec:stability}) implies that $\rho(t)$ remains radial about $c(0)$ at every time. Its barycentre is therefore $c(t)=c(0)$. By Theorem~\ref{thm:binarity}, the full trajectory converges in $L^1(\R^2)$ to a binary density $\rho_\infty=R1_E$ whose barycentre is $c(0)$. The strong limit is radial about the same point.

On the other hand, Proposition~\ref{prop:star_shaped} implies that $E$ is essentially star-shaped with respect to $c(0)$, i.e. $1_{E-c(0)}\le 1_{\lambda(E-c(0))}$ almost everywhere for every $\lambda\ge1$. A radial measurable set with this property is, up to a null set, a ball centred at the origin: indeed, in polar coordinates the characteristic function depends only on the radius, and the dilation ordering forces that radial profile to be non-increasing almost everywhere (one may first use a countable dense family of dilations and then continuity of dilations in $L^1_{\mathrm{loc}}$). Hence $E=B_r(c(0))$ up to a null set for some $r>0$. Mass conservation gives $R\pi r^2=1$, so $r=r_*=(\pi R)^{-1/2}$, proving \eqref{eq:strong_convergence_disk}.
\end{proof}

\section*{Acknowledgements} 

The authors are members of GNAMPA - Istituto Nazionale di Alta Matematica (INdAM). The authors acknowledge the support of the INdAM - GNAMPA 2026 Project ''Modelli aggregazione-diffusione con mobilit`a  non lineare: buona positura, comportamento asintotico`` (CUP: E53C25002010001). E.~Radici acknowledges the support of the Italian ''National Centre for HPC, Big Data and Quantum Computing`` - Spoke 5 ''Environment and Natural Disasters``, and of the InterMaths Network \url{www.intermaths.eu}.

\end{document}